# ARITHMÉTIQUE DES COURBES ELLIPTIQUES À RÉDUCTION SUPERSINGULIÈRE EN $p$

BERNADETTE PERRIN-RIOU

Résumé : *Nous faisons le point sur la conjecture principale pour une courbe elliptique sur $\mathbb{Q}$ ayant bonne rduction supersingulire en p et en donnons quelques consquences. Puis nous dfinissons la notion de $\lambda$ invariant et de $\mu$ invariant dans cette situation, gnralisant un travail de Kurihara et en dduisons la forme de l'ordre du groupe de Shafarevich-Tate le long de la $\mathbb{Z}_p$-extension cyclotomique. Par des exemples, nous donnons quelques arguments qui, en alliant théorèmes et calculs numériques, permettent de calculer l'ordre de la composante p-primaire du groupe de Shafarevich-Tate dans des cas non connus jusqu'à présent (groupe de Shafarevich-Tate non trivial, courbes de rang $\geq 1$).*

Dans l'article "Parabolic points and zeta functions of modular curves" en 1972 ([12]), Manin donne entre autres des formules explicites pour les valeurs des fonctions $L$ associées à une forme modulaire de poids 2 en termes de symboles modulaires. Les reliant conjecturalement à la conjecture de Birch et Swinnerton-Dyer, il exprime le comportement conjectural du rang de la courbe elliptique associée $E$ le long de la $\mathbb{Z}_p$-extension cyclotomique $\mathbb{Q}_\infty$. Lorsque $E$ a bonne réduction ordinaire, Mazur montre dans [13] une partie de ces conjectures en utilisant la théorie d'Iwasawa etla technique des $\Gamma$-modules et énonce ce qu'on appelle désormais "les conjectures principales". Pour le cas où $E$ a bonne réduction supersingulière en $p$, citons Manin : *les nombres premiers supersinguliers ont résisté jusqu'à maintenant faute d'une technique de $\Gamma$-modules.*

Ces techniques existent maintenant : interpolation des valeurs spéciales $p$-adiquement depuis déjà longtemps (Manin-Vishik, Amice-Vélu, Mazur-Tate-Teitelbaum) et plus récemment théorie d'Iwasawa complète : énoncé des conjectures principales, calcul des valeurs spéciales en termes du groupe de Shafarevich-Tate. Par exemple, dans [16], sous une hypothèse de non-nullité d'une fonction $L$ $p$-adique, il est montré que $E(\mathbb{Q}_\infty)$ est de rang fini.[1] Le théorème de Kato dont la démonstration utilise un système d'Euler-Kolyvagin permet de montrer cette non-nullité et d'obtenir une partie de la conjecture principale. En utilisant le théorème de Kato, Kurihara [11] a récemment montré que si $L(E,1)/\Omega_E$ est une unité, $E(\mathbb{Q}_\infty)$ est fini et $\text{III}(E/\mathbb{Q}_n)(p)$ est fini pour toute extension finie $\mathbb{Q}_n/\mathbb{Q}$ contenue dans $\mathbb{Q}_\infty$ et a calculé l'ordre de $\text{III}(E/\mathbb{Q}_n)(p)$.

Dans cet article, nous donnons une généralisation du théorème de Kurihara. En alliant alors ces théorèmes à des calculs numériques, nous sommes en mesure de donner des exemples où l'on peut montrer la conjecture principale complète, où il est possible de calculer l'ordre de la composante $p$-primaire du groupe de Shafarevich-Tate et de **montrer** qu'il est égal à l'ordre conjectural. Il est possible aussi de calculer de manière non conjecturale la manière dont il varie le long de la $\mathbb{Z}_p$-extension cyclotomique. Précisons bien qu'**il s'agit vraiment de l'ordre et**

---

*Date*: 12 juillet 2001.

[1]Plus exactement, pour que l'énoncé de *loc. cit.* ne soit pas conjectural, l'hypothèse mise était que $E(\mathbb{Q})$ et le groupe de Tate-Shafarevich sont finis, mais l'hypothèse qu'une certaine fonction $L$ $p$-adique est non nulle est suffisante, hypothèse maintenant démontrée par Kato.

**non de l'ordre conjectural**. Les calculs numériques en question ne font intervenir que des calculs de symboles modulaires et donnent des renseignements sur les invariants arithmétiques de la courbe elliptique sur $\mathbb{Q}_\infty$. On revient ainsi aux sources des conjectures principales qui relient des données analytiques à des données arithmétiques. Beaucoup de cas particuliers ne démontrent pas une conjecture. Cependant, les arguments mis en œuvre sont, me semblent-ils, intéressants car ils montrent ce qu'il est possible de tirer de la conjecture à partir de quelques calculs de symboles modulaires.

Les calculs ont été faits en utilisant le logiciel `GP/Pari` [1], une version préliminaire d'un programme `GP` de Joseph L. Wetherell (que l'on peut trouver à l'adresse [24]) modifiée avec l'aide de Dominique Bernardi pour calculer les symboles modulaires et les valeurs de fonctions $L$ $p$-adiques, les courbes calculées par J. Cremona et les programmes `mwrank` et `ratpoints` pour calculer des points dans le groupe de Mordell-Weil de $E$ ([7]), le programme `GP` de Tom Womack pour calculer l'ordre conjectural du groupe de Shafarevich-Tate. Je remercie les auteurs de ces programmes et données en libre accès. Quelques programmes seront mis à l'adresse [26]. Ils ont tourné à la fois sur des ordinateurs Macintosh et sur les serveurs du GDR Médicis de l'école Polytechnique ([25]).

Je remercie avec grand plaisir Dominique Bernardi et Karim Belabas pour l'aide qu'ils m'ont apportée.

## 1. Fonction $L$ $p$-adique analytique et congruences provenant de la distribution modulaire

La construction de la fonction $L$ $p$-adique associée à une courbe elliptique modulaire repose sur l'existence et les propriétés des symboles modulaires. Rappelons rapidement quelques notations. Les références sont [12], [15], [7], [24].

1.1. Soit $E$ une courbe elliptique définie sur $\mathbb{Q}$. Soit $D(E)$ le quotient de l'espace vectoriel des formes différentielles sur $E$ définies sur $\mathbb{Q}$ par le sous-espace vectoriel des formes différentielles exactes. C'est un $\mathbb{Q}$-espace vectoriel de dimension 2 contenant la droite $\operatorname{Fil}^0 D(E)$ engendrée par l'image des formes différentielles invariantes. Soit $\omega_E$ une forme différentielle de Néron sur $\mathbb{Z}$. Si
$$y^2 + a_1 xy + a_3 y = x^3 + a_2 x^2 + a_4 x + a_6$$
est un modèle minimal de Weierstrass, on peut prendre $\omega_E = \frac{dx}{2y+a_1 x+a_3}$. Les classes de $\omega_E$ et de $\eta = x\omega_E$ forment une base de $D(E)$. Ainsi, en tant que $\mathbb{Q}$-espace vectoriel, $D(E) = H^1_{dR}(E/\mathbb{Q})$. En tant que module filtré, il s'agit plutôt de $H^1_{dR}(E/\mathbb{Q})[-1]$. Le $\mathbb{Q}$-espace vectoriel $D(E)$ est muni d'une forme bilinéaire alternée $[\cdot,\cdot]_{D(E)}$ vérifiant $[\omega_E, \eta]_{dR} = 1$. Soit $N_E$ le conducteur de $E$. Notons $f$ la forme modulaire associée à $E$, choisissons une paramétrisation $\pi : X_0(N_E) \to E$ de degré minimal, $c_\pi$ la constante associée : $\omega_E = c_\pi \omega_f$ avec $\pi^* \omega_f = f(z) 2 i \pi dz$.

Rappelons que $D_p(E) = \mathbb{Q}_p \otimes D(E)$ admet un endomorphisme $\varphi$ dit de Frobenius. On renvoie à la note [3] pour une description concrète. Le polynôme caractéristique de $\varphi$ tel qu'il a été choisi est alors $X^2 - p^{-1} a_p X + p^{-1}$. On pose $L(E/\mathbb{Q}_p, s) = (1 - a_p p^{-s} + p^{1-2s})^{-1}$.

Soit $M_E$ le $\mathbb{Z}_p$-réseau engendré par $\omega_E$ et $\eta$ dans $D_p(E)$. L'endomorphisme $p\varphi$ laisse stable $M_E$. D'autre part, $\varphi \omega_E \in M_E$. Comme $p$ est supposé supersingulier, il est facile de voir que $[\varphi \omega_E, \omega_E]_{D_p(E)}$ est une unité (si $\begin{pmatrix} p\gamma & \alpha \\ p\delta & \beta \end{pmatrix}$ est la matrice de $p\varphi$ dans la base $(\omega_E, \eta)$ avec $\alpha, \beta, \gamma, \delta \in \mathbb{Z}_p$, on a $p\gamma + \beta = a_p$, $\gamma\beta - \alpha\delta = 1$; comme $a_p$ est divisible par $p$, il en est de même de $\beta$; donc $[\varphi \omega_E, \omega_E]_{D_p(E)} = \delta$ est une unité). On en déduit que $M_E = \mathbb{Z}_p \omega_E \oplus \mathbb{Z}_p \varphi \omega_E$ (noter le changement de notation par rapport à [3]). Remarquons dès maintenant que lorsque $E$ a bonne réduction supersingulière en $p$, $p$ ne divise pas $c_\pi$ (A.2). Ainsi, si $M_f = \mathbb{Z}\omega_f + \mathbb{Z}_p\varphi\omega_f$, $M_E = c_\pi M_f = M_f$.

Soit $c_0$ le nombre de composantes connexes de $E(\mathbb{R})$. Soit $m$ un entier tel que $\pi_*\{r\} \in H_1(E, m^{-1}\mathbb{Z})$ où $\{r\}$ est l'image dans $X_0(N_E)(\mathbb{C})$ d'un chemin joignant $\infty$ à $r$ (son existence est due à un théorème de Manin).

Soient $\gamma^+$ et $\gamma^-$ des bases des $\mathbb{Z}$-modules $H_1(E,\mathbb{Z})^+$ et $H_1(E,\mathbb{Z})^-$ telles que $\int_{\gamma^+}\omega_E \in \mathbb{R}^+$ et $\int_{\gamma^-}\omega_E \in i\mathbb{R}^+$. Si $r$ est un rationnel, on note $x^\pm(r)$ les composantes de $-\pi^*\{-r\}$ dans la base $(\gamma^+, \gamma^-)$. Ce sont des éléments de $m^{-1}\mathbb{Z}$. Lorsque le dénominateur de $r$ est premier à $N_E$, $x^\pm(r)$ est de plus entier en $p$ (la pointe $0$ et la pointe $r$ sont alors équivalentes sous le sous-groupe de congruence $\Gamma_0(N_E)$, $(p-1-a_p)\pi_*\{0\}$ s'exprime en termes des $\pi_*\{\frac{a}{p}\} - \pi_*\{\infty\}$ (avec $a \in \mathbb{Z}$ premier à $p$) qui sont entiers et $p-1-a_p$ est une unité en $p$). On pose $c_1 = c_\pi c_0$ et $C = 2c_\pi m$. Pour $p$ impair supersingulier pour $E$, ce sont des entiers premiers à $p$.

1.2. Le symbole modulaire vérifie des propriétés de congruences qui peuvent se traduire comme l'existence d'une distribution $p$-adique d'ordre $1/2$ (Vishik, Amice-Vélu). Ainsi, pour $a$ premier à $p$, les

$$\mu(a + p^n\mathbb{Z}_p) = c_1^{-1}(x^+(a/p^n)\varphi^n(\omega_E) - x^+(a/p^{n-1})\varphi^{n+1}(\omega_E)) \in D_p(E) .$$

définissent une distribution sur $\mathbb{Z}_p^*$ et vérifient $\mu(a + p^n\mathbb{Z}_p) \in C^{-1}p^{-(\lfloor n/2 \rfloor+1)}M_E$.

On pose $L_p(E)(\rho) = \int_{\mathbb{Z}_p^*} \rho(\chi^{-1}a)d\mu$ un des avatars de la fonction $L$ $p$-adiques pour $\rho$ caractère continu de $G_\infty$ à valeurs dans $\mathbb{C}_p^*$. Nous aurons besoin d'autres versions. Soit $\mathcal{H}$ l'anneau des fonctions analytiques sur le disque unité $\{|x| < 1\}$ vérifiant une propriété de croissance (par exemple $\sup_n p^{nh}||f||_{p^{-1/(p-1)p^n}} < \infty$ pour un $h \in \mathbb{R}$), $\mathcal{H}(\Gamma)$ l'image de $\mathcal{H}$ par $x \mapsto (\gamma - 1)$ avec $\gamma$ générateur de $\Gamma = \mathrm{Gal}(\mathbb{Q}(\mu_{p^\infty})/\mathbb{Q}(\mu_p))$ et $\mathcal{H}(G_\infty) = \mathbb{Z}_p[\mathrm{Gal}(\mathbb{Q}(\mu_p)/\mathbb{Q})] \otimes \mathcal{H}(\Gamma)$. On a alors $L_p(E)(\rho) = \rho(\hat{L}_p(E))$ avec $\hat{L}_p(E) \in \mathcal{H}(G_\infty) \otimes D_p(E)$, ou encore si $\gamma$ est un générateur de $\mathrm{Gal}(\mathbb{Q}(\mu_{p^\infty}/\mathbb{Q}(\mu_p))$ ou ce qui revient au même de $1 + 2p\mathbb{Z}_p$ (à travers le caractère cyclotomique $\chi$) et $e_i$ l'idempotent de $\Delta = \mathrm{Gal}(\mathbb{Q}(\mu_p)/\mathbb{Q})$ associé au caractère $\chi^i_{|\mathrm{Gal}(\mathbb{Q}(\mu_p)/\mathbb{Q})}$

$$\hat{L}_p(E) = \sum_i \hat{L}_{p,(i)}(E)(\gamma - 1)e_i$$

avec $\hat{L}_{p,(i)}(E) \in \mathcal{H} \otimes D_p(E)$, fonction analytique de $x$. Ainsi, si Teich est le caractère de Teichmüller,

(1.2.1) $$\hat{L}_{p,(i)}(E) = \int_{\mathbb{Z}_p^*} \mathrm{Teich}^i(u)(1+x)^{\frac{\log_p u}{\log_p \chi(\gamma)}}d\mu(u) .$$

On a alors $\int_{\mathbb{Z}_p^*} x^k d\mu = L_p(E)(\chi^k) = \hat{L}_{p,(k)}(E)(\langle\chi(\gamma)\rangle^k - 1)$. La dérivée d'ordre $k$ en $\mathbf{1}$ (resp. en caractère $\delta$ de Dirichlet de conducteur une puissance de $p$) se calcule par la formule

$$L_p^{(k)}(E)(\mathbf{1}) = \int_{\mathbb{Z}_p^*} \log_p^k a \; d\mu , \quad L_p^{(k)}(E)(\delta) = \int_{\mathbb{Z}_p^*} \delta(a) \log_p^k a \; d\mu .$$

Ainsi, si $\delta$ est un caractère d'ordre $p^n$ et de conducteur $p^{n+1}$ que l'on voit soit comme un caractère de Dirichlet, soit comme un caractère de $\mathrm{Gal}(\mathbb{Q}(\mu_{p^{n+1}})/\mathbb{Q}(\mu_p)) = \mathrm{Gal}(\mathbb{Q}_n/\mathbb{Q})$ avec $\mathbb{Q}_n = \mathbb{Q}(\mu_{p^{n+1}})^\Delta$, on a

$$L_p^{(k)}(E)(\delta) = (\log_p \chi(\gamma))^{-k} \hat{L}_{p,(0)}^{(k)}(E)(\delta(\gamma) - 1)$$

que l'on note aussi $\delta(L_p^{(k)}(E))$, ce qui ne devrait pas porter à confusion.

La fonction $L_p(E)$ vérifie l'équation fonctionnelle
$$L_p(E)(\rho^{-1}) = \epsilon_{anal} L_p(E)(\rho)$$
où $\epsilon_{anal}$ est le signe de l'équation fonctionnelle complexe : on a donc $\epsilon_{anal} = (-1)^{r_{anal}^\infty} = (-1)^{r_{anal}^p}$ où $r_{anal}^\infty$ est l'ordre d'annulation de $L(E,s)$ en 1 et $r_{anal} = r_{anal}^p$ l'ordre d'annulation de $L_p(E)$ en **1**.

Cet article reposant en partie sur des calculs numériques, indiquons dès maintenant quelques outils de calculs de ces fonctions. Nous prendrons alors $\chi(\gamma) = 1+p$ (pour $p$ impair). Considérons d'abord les polynômes d'interpolation des fonctions $\hat{L}_{p,(i)}(E)$ modulo $\omega_n(x)$ avec $\omega_n(x) = (1+x)^{p^n} - 1$. Plus précisément, si $\xi_n(x) = \omega_n(x)/\omega_{n-1}(x)$, nous montrerons que si $j$ est une classe modulo $p-1$, il existe une et une seule famille de polynômes $P_n^{(j)}$ pour $n \geq 0$, de degré $< p^n$ vérifiant

$$\hat{L}_{p,(j)}(E) \equiv P_n^{(j)} \varphi^{n+1} \omega_E - \xi_n P_{n-1}^{(j)} \varphi^{n+2} \omega_E \bmod \omega_n(x) D_p(E)$$

pour $n \geq 1$. Explicitement, on a

$$P_n^{(j)} = c_1^{-1} \sum_{\substack{a=0 \\ (a,p)=1}}^{p^{n+1}} \mathrm{Teich}^j(a) x^+ (\frac{a}{p^{n+1}})(1+x)^{r_n(a)}$$

où $r_n(a)$ est déterminé par $r_n(a) \equiv \frac{\log_p a}{\log_p 1+p} \bmod p^n$ et $0 \leq r_n(a) < p^n$. Ce polynôme est à coefficients dans $C^{-1}\mathbb{Z}_p = \mathbb{Z}_p$.

Si $a$ et $b$ sont deux vecteurs colinéaires d'un $\mathbb{Q}_p$-espace vectoriel : $a = \lambda b$, on écrit $a \sim b$ si $\lambda$ est une unité de $\mathbb{Z}_p$.

**1.3. Conjecture de Birch et Swinnerton-Dyer $p$-adique.** La conjecture de Birch et Swinnerton-Dyer $p$-adique dans le cas supersingulier a la forme suivante ([3]) :

**Conjecture.** *Soit $r$ le rang de $E(\mathbb{Q})$ et $\mathrm{III}(E/\mathbb{Q})$ le groupe de Shafarevich-Tate de $E/\mathbb{Q}$ : alors,*

$$(1-\varphi)^{-1}(1-p^{-1}\varphi^{-1}) \frac{\mathbf{1}(L_p^{(r)}(E))}{r!} = L(E/\mathbb{Q}_p, 1)^{-1} \frac{\mathrm{Tam}(E)}{\sharp E(\mathbb{Q})_{tor}^2} \sharp\mathrm{III}(E/\mathbb{Q}) R^{(r)}(E) \ .$$

Ici $\mathrm{Tam}(E)$ est le produit des nombres de Tamagawa locaux ; $R^{(r)}(E)$ est un élément de $D_p(E)$ dont le calcul demande la connaissance de $E(\mathbb{Q})$ (régulateur $p$-adique du groupe de Mordell-Weil $E(\mathbb{Q})$). On peut le décrire de la manière suivante. Rappelons d'abord qu'une fois choisie la forme différentielle de Néron $\omega_E$, on associe à tout élément $\nu$ de $D_p(E)$ une forme quadratique $h_\nu$ sur $E(\mathbb{Q})$, appelée hauteur $p$-adique et on note $\langle \cdot, \cdot \rangle_\nu$ la forme bilinéaire associée (voir [14], [3] pour une description en termes de fonctions $\sigma$). Ainsi, $h_{\omega_E}(P) = -(\log_{\omega_E} P)^2$ où $\log_{\omega_E}$ est le logarithme sur $E(\mathbb{Q})$ associé à la forme différentielle invariante $\omega$. Supposons $r \geq 1$. Le noyau de la forme bilinéaire $\langle \cdot, \cdot \rangle_{\omega_E}$ est de rang $r-1$, c'est le noyau de localisation $\mathbb{Z}_p \otimes_\mathbb{Z} E(\mathbb{Q}) \to \mathbb{Z}_p \otimes E(\mathbb{Q}_p)$ que l'on note $(\mathbb{Z}_p \otimes_\mathbb{Z} E(\mathbb{Q}))_0$. Les $h_\nu$ varient linéairement en $\nu \in D_p(E)$. On en déduit que la restriction de $\langle \cdot, \cdot \rangle_\nu$ à $(\mathbb{Z}_p \otimes_\mathbb{Z} E(\mathbb{Q}))_0$ est indépendante de $\nu \notin \mathrm{Fil}^0 D_p(E) = \mathbb{Q}_p \omega_E$ : on la note $\langle\langle \cdot, \cdot \rangle\rangle$. L'orthogonal de $(\mathbb{Z}_p \otimes_\mathbb{Z} E(\mathbb{Q}))_0$ est indépendant de $\nu$. Cette forme bilinéaire est conjecturée être non dégénérée sur $(\mathbb{Z}_p \otimes_\mathbb{Z} E(\mathbb{Q}))_0$. Il existe un unique élément $R^{(r)}(E) \in D_p(E)$ tel que

$$[R^{(r)}(E), \nu]_{D_p(E)} = \frac{\det\left((\langle P_i, P_j \rangle_\nu)\right)}{[E(\mathbb{Q})/E(\mathbb{Q})_{tors} : \sum_{i=1}^r \mathbb{Z} P_i]}$$

pour $\nu \notin \mathrm{Fil}^0 D_p(E)$ et $(P_i)$ un système libre de rang $r$ de $E(\mathbb{Q})$. On peut le décrire explicitement par

(1) si $r=0$, $R^{(0)}(E) = \omega_E$ ;

(2) si $r = 1$,
$$R^{(1)}(E) = \frac{h_\nu(P)\omega_E - h_{\omega_E}(P)\nu}{[\omega_E, \nu]_{D_p(E)}}$$

avec $P$ un générateur de $E(\mathbb{Q})$ modulo torsion et $h_\nu$ la hauteur $p$-adique associée à $\nu$; remarquons que $h_{\omega_E}(P)$ n'est jamais nul car $P$ n'est pas de torsion;

(3) si $r \geq 2$,
$$R^{(r)}(E) = \frac{\det\left(\left(\langle P_i, P_j\rangle_\nu\right)\right)\nu' - \det\left(\left(\langle P_i, P_j\rangle_{\nu'}\right)\right)\nu}{[E(\mathbb{Q})/E(\mathbb{Q})_{tors} : \sum_{i=1}^r \mathbb{Z}P_i][\nu', \nu]_{D_p(E)}} \ .$$
$$= \frac{\disc\left(\langle \cdot, \cdot\rangle_\nu \nu' - \langle \cdot, \cdot\rangle_{\nu'}\nu\right)}{[\nu', \nu]_{D_p(E)}}$$

pour $(\nu, \nu')$ une base de $D_p(E)$ avec $\nu$ et $\nu'$ n'appartenant pas à $\Fil^0 D_p(E)$. À une unité $p$-adique près, $R^{(r)}(E)$ est aussi le produit du discriminant de $\langle\langle \cdot, \cdot\rangle\rangle$ sur $(\mathbb{Z}_p \otimes_\mathbb{Z} E(\mathbb{Q}))_0$ par $(h_\nu(P_u)\omega_E - h_{\omega_E}(P_u)\nu)/[\omega_E, \nu]$ si $P_u$ est une base de l'orthogonal de $(\mathbb{Z}_p \otimes_\mathbb{Z} E(\mathbb{Q}))_0$ dans $\mathbb{Z}_p \otimes_\mathbb{Z} E(\mathbb{Q})$ :
$$R^{(r)}(E) \sim \disc\langle\langle \cdot, \cdot\rangle\rangle \frac{h_\nu(P_u)\omega_E - h_{\omega_E}(P_u)\nu}{[\omega_E, \nu]}$$

et on a
$$[R^{(r)}(E), \omega_E]_{D_p(E)} \sim \disc\langle\langle \cdot, \cdot\rangle\rangle \log^2_{\omega_E}(P_u) \ .$$

## 2. Fonction $L$ $p$-adique arithmétique

2.1. **Définition du module arithmétique.** Soit $T_p(E)$ le module de Tate des points de $p^\infty$-torsion de $E$ et $V_p(E) = \mathbb{Q}_p \otimes T_p(E)$. On définit le régulateur $p$-adique [18]
$$\mathcal{L}_E : \varprojlim_n H^1(\mathbb{Q}_p(\mu_{p^n}), T_p(E)) \to \mathcal{H}(G_\infty) \otimes D_p(E) \ .$$

Il est d'ordre $\leq 0$ au sens suivant : on dit que $f \in \mathbb{Q}_p[[x]] \otimes D_p(E)$ analytique sur le disque unité $\{|x|_p < 1\}$ de $\mathbb{C}_p$ est d'ordre $\leq 0$ si pour un (ou tout) $\rho < 1$, les $||(1 \otimes \varphi)^{-n} f||_{\rho^{1/p^n}}$ sont bornés avec $n$.

Soit $H^1_\infty(\mathbb{Q}, T_p(E))$ la limite projective des groupes de cohomologie galoisienne à ramification limitée $H^1(G_{S,\mathbb{Q}(\mu_{p^n})}, T_p(E))$ et $H^2_{\infty,p}(\mathbb{Q}, T_p(E))$ la limite projective des noyaux $H^2(G_{S,\mathbb{Q}(\mu_{p^n})}, T_p(E)) \to \oplus_{v \in S} H^2(\mathbb{Q}(\mu_{p^n})_v, T_p(E))$ ($S$ est un ensemble de places contenant les places divisant $p$ et les places où $E$ a mauvaise réduction et $G_{S,K}$ le groupe de Galois sur $K$ de la plus grande extension de $K$ non ramifiée en dehors des places au dessus de $S$). Définissons alors $\mathbb{I}_{arith}$ comme le sous-$\mathbb{Z}_p[[G_\infty]]$-module de $\mathcal{H}(G_\infty) \otimes D_p(E)$ engendré par l'image de
$$\left(\det_{\mathbb{Z}_p[[G_\infty]]} H^2_{\infty,p}(\mathbb{Q}, T_p(E))\right)^{-1} \otimes H^1_\infty(\mathbb{Q}, T_p(E))$$

par $\mathcal{L}_E$, ce qui signifie pratiquement que si $F_2$ est une série caractéristique du module de $\mathbb{Z}_p[[G_\infty]]$-module de torsion $H^2_{\infty,p}(\mathbb{Q}, T_p(E))$ et si $\mathfrak{c}$ est un élément de $H^1_\infty(\mathbb{Q}, T_p(E))$ tel que $H^1_\infty(\mathbb{Q}, T_p(E))/\mathbb{Z}_p[[G_\infty]]\mathfrak{c}$ soit de $\mathbb{Z}_p[[\Gamma]]$-torsion (et de série caractéristique $F_\mathfrak{c}$), on a
$$\mathbb{I}_{arith}(E/\mathbb{Q}) = \mathbb{Z}_p[[G_\infty]] F_\mathfrak{c}^{-1} F_2 \mathcal{L}_E(\mathfrak{c}) \ .$$

Cette définition est conforme à l'idée de l'utilisation d'un régulateur pour mesurer un objet compliqué. On note $I_{arith}(E/\mathbb{Q})$ un générateur de $\mathbb{I}_{arith}(E/\mathbb{Q})$.

## 2.2. Formule arithmétique de Birch et Swinnerton-Dyer.
Soit $Sel(E, p^n)$ le groupe de Selmer relatif à la multiplication par $p^n$. On définit

$$Sel_p(E/\mathbb{Q}) = \varinjlim_n Sel(E, p^n) = H^1_f(\mathbb{Q}, V_p(E)/T_p(E))$$

$$\check{S}_p(E) = \varprojlim_n Sel(E, p^n) = H^1_f(\mathbb{Q}, T_p(E))$$

(les deuxièmes notations sont celles de Bloch-Kato). Rappelons que l'on a une suite exacte

$$0 \to \mathbb{Z}_p \otimes E(\mathbb{Q}) \to \check{S}_p(E) \to T_p(\mathbf{III}(E/\mathbb{Q})) \to 0 \ .$$

où $T_p(\mathbf{III}(E/\mathbb{Q}))$ est le module de Tate du groupe de Shafarevich-Tate de $E/\mathbb{Q}$. Les formes bilinéaires introduites précédemment se prolongent naturellement à $\check{S}_p(E)$, de même que $\langle\langle \cdot, \cdot \rangle\rangle$ au noyau de localisation $\check{S}_p(E)_0$. Soit $\mathbf{III}(T_p(E)/\mathbb{Q})(p)$ le quotient de $Sel_p(E/\mathbb{Q})$ par sa partie divisible maximale, ce qui est aussi le quotient de $\mathbf{III}(E/\mathbb{Q})(p)$ par sa partie divisible maximale et est égal à $\mathbf{III}(E/\mathbb{Q})(p)$ lorsque ce dernier est fini. On a alors le théorème ([18]) :

**2.2.1. Théorème.** $\mathbb{I}_{arith}(E/\mathbb{Q})$ a un zéro en **1** de multiplicité supérieure ou égale à $r_{arith} = \mathrm{rg}_{\mathbb{Z}_p} \check{S}_p(E)$. Il est égal à $r_{arith}$ si et seulement si la forme bilinéaire $\langle\langle \cdot, \cdot \rangle\rangle$ est non dégénérée sur $\check{S}_p(E)_0$. On a alors

$$(1-\varphi)^{-1}(1-p^{-1}\varphi^{-1})\mathbf{1}(I^*_{arith}(E/\mathbb{Q}))$$
$$\sim L(E/\mathbb{Q}_p, 1)^{-1} \frac{\mathrm{Tam}(E)}{\sharp E(\mathbb{Q})^2_{tor}} \sharp\mathbf{III}(T_p(E)/\mathbb{Q})(p) R^{(r)}(E) \in D_p(E)$$

où $R^{(r)}(E)$ est défini comme en 1.3 en remplaçant $\mathbb{Z}_p \otimes E(\mathbb{Q})$ par $\check{S}_p(E)$.

Ici, $\mathbf{1}(I^*_{arith}(E/\mathbb{Q})) = \dfrac{\mathbf{1}(I^{(r_{arith})}_{arith}(E/\mathbb{Q}))}{r_{arith}!}$.

**Remarque.** Ce théorème est démontré dans [18] et dans un cadre plus général dans [20]. Avec la notation de [20, chap. 3], le sous-espace $N$ engendré par $\nu \in D_p(E)$ est régulier si et seulement si $\langle\langle \cdot, \cdot \rangle\rangle$ est non dégénérée; les suites $s_N$ et $s_{f,N}$ sont réduites à

$(s_N) \qquad\qquad\qquad 0 \to N \to D_p(E)/\mathrm{Fil}^0 D_p(E) \to 0$

$(s_{f,N}) \qquad\qquad\qquad 0 \to \check{S}_p(E)^* \to \check{S}_p(E) \to 0$

dès que $N \neq \mathrm{Fil}^0 D_p(E)$ et à

$(s_N) \quad 0 \to N \to N \xrightarrow{0} D_p(E)/\mathrm{Fil}^0 D_p(E) \to D_p(E)/\mathrm{Fil}^0 D_p(E) \to 0$

$(s_{fN}) \qquad\qquad\qquad 0 \to \check{S}_p(E)^*_0 \to \check{S}_p(E)_0 \to 0$

lorsque $N = \mathrm{Fil}^0 D_p(E)$ (la forme bilinéaire $\langle \cdot, \cdot \rangle_N$ est dans ce cas égale à $\langle\langle \cdot, \cdot \rangle\rangle$).

Le théorème 2.2.1 se généralise à $K = \mathbb{Q}(\mu_{p^n})$ (et même à une extension finie de $\mathbb{Q}$, composée d'une extension de $\mathbb{Q}$ non ramifiée en $p$ et d'une sous-extension finie de $\mathbb{Q}(\mu_{p^\infty})$). Soit $\Delta_K$ le groupe de Galois de $K/\mathbb{Q}$ et $\hat{\Delta}_K$ le groupe des caractères de $\Delta_K$. Si $\delta$ est un caractère de $\Delta$, soit $p^{n_\delta}$ la $p$-partie de son conducteur. On note $G_{\infty,K}$ le groupe de Galois de $K(\mu_{p^\infty})/K$. On attache à un élément $v$ de $D_p(E)$ et à $K$ l'élément suivant

$$\tilde{v}_K = \otimes_{\delta \in \hat{\Delta}_K} \varphi^{n_\delta} v \otimes (1-\delta(p)\varphi)(1-p^{-1}\overline{\delta}(p)\varphi^{-1})^{-1} v$$

de $\kappa \otimes \otimes^{[K:\mathbb{Q}_p]}_{\mathbb{Q}_p} D_p(E)$ avec $\kappa$ le corps de coefficients des caractères de $\Delta_K$. Remarquons que $\tilde{\omega}_{\mathbb{Q}_p} = (1-\varphi)(1-p^{-1}\varphi^{-1})^{-1}\omega$, que si $\delta$ est un caractère de Dirichlet non ramifié en p (i.e. $\delta(p) \neq 0$, $n(\delta) = 0$), la contribution de $\delta$ est $(1 - \delta(p)\varphi)(1 - p^{-1}\delta(p)^{-1}\varphi^{-1})^{-1}$. Lorsque $\delta$ est ramifié (donc $\delta(p) = 0$), sa contribution est $\varphi^{n(\delta)}$. Si $\chi_K$ est le caractère cyclotomique de $G_{\infty,K}$, on a $\chi_K(\tau^{[K:\mathbb{Q}_p]}) = \chi(\tau)^{[K:\mathbb{Q}_p]}$. On pose alors

$$I_{arith}(E/K)(\chi_K^s) = \otimes_\delta I_{arith}(E/K)(\delta\chi^s) \in Funct(\hat{G}_{\infty,K}, \mathbb{C}_p) \otimes^{[K:\mathbb{Q}_p]} D_p(E) \ .$$

Notons $r_K$ le rang de $\check{S}_p(E/K)$ et $I_{arith}^{(r_K)}(E/K)(\mathbf{1}_K)$ la dérivée $r_K$-ième en $\mathbf{1}_K$ (caractère trivial de $G_{\infty,K}$ par rapport à $\chi_K$). Soit $v$ un élément de $D_p(E)$. Si $N = Kv$ est le $K$-espace vectoriel engendré par $v$, vu comme $\mathbb{Q}_p$-espace vectoriel, on peut lui associer comme dans [20] (voir aussi [6]) une forme bilinéaire $\langle \cdot, \cdot \rangle_{Kv}$ sur $\check{S}_p(E/K)$ à valeurs dans $\mathbb{Q}_p$ (et même dans $\mathbb{Z}_p$). On note $\langle\langle \cdot, \cdot \rangle\rangle_K$ sa restriction au noyau de localisation $\check{S}_p(E/K)_0$ en $p$.

**2.2.2. Théorème.** 1) *On a*

$$I_{arith}(E/K)(\mathbf{1}_K) = \otimes_{\delta \in \hat{\Delta}_K} \delta(I_{arith}) \in \mathbb{Q}_p \tilde{\omega}_{E,K} \ .$$

2) *Si $\delta(I_{arith})$ est non nul, alors le rang de $E(K)^{(\delta)}$ est nul et $\text{III}(K)^{(\delta)}$ est fini.*

3) *En particulier, si $I_{arith}(E/K)(\mathbf{1})$ est non nul, $E(K)$ et $\text{III}(E/K)(p)$ sont finis. On a alors*

$$I_{arith}(E/K)(\mathbf{1}_K) \sim L(E/K_p, 1)^{-1} \frac{\text{Tam}(E/K) \sharp \text{III}(E/K)}{\sharp E(K)^2} \tilde{\omega}_{E,K} \ .$$

4) *Plus généralement, si $\check{S}_p(E/K)$ est de rang $r_K \geq 1$, et $v$ un élément de $D_p(E)$ qui n'appartient pas à $\text{Fil}^0 D_p(E)$,*

$$[\frac{I_{arith}^{(r_K)}(E/K)(\mathbf{1}_K)}{r_K!}, \tilde{v}_K]_{D(E)} / [\omega_E, v]_{D_p(E)}^{[K:\mathbb{Q}]}$$
$$\sim L(E/K_p, 1)^{-1} \text{Tam}(E/K) \sharp \text{III}(T_p(E)/K)(p) \text{disc}\langle \cdot, \cdot \rangle_{K \otimes v}$$

La démonstration est semblable à la démonstration du théorème 3.5.2 de [20] en prenant $N = Kv$. Nous ne redonnons pas la preuve de l'assertion sur l'ordre du zéro, mais seulement les étapes du calcul du coefficient dominant. La modification des facteurs d'Euler qui sont ici cachés dans $\tilde{v}_K$ provient de ce que les formules donnant les valeurs spéciales de $\mathcal{L}_E$ en un caractère $\delta$ de conducteur $p^n$ avec $n \geq 0$ font intervenir $\varphi^n$.

On reprend les notations du §2.1 et posons $K = \mathbb{Q}(\mu_{p^n})$. On pose $G_{\infty,K} = \text{Gal}(K_\infty/K)$ et $\Delta_K = \text{Gal}(K/\mathbb{Q})$. On note $K_p = \mathbb{Q}_p \otimes_\mathbb{Q} K = \prod_{v|p} K_v$, $\mathcal{O}_{K_p} = \mathbb{Q}_p \otimes_\mathbb{Q} \mathcal{O}_K$. Pour $\delta$ caractère de $\Delta_K = \text{Gal}(K/\mathbb{Q})$ et $\mathcal{O}_\delta$ l'anneau engendré par les valeurs de $\delta$, on note $e_\delta = \sum_{\sigma \in \Delta_K} \delta^{-1}(\sigma)\sigma$. Soit $H^1(K, T_p(E))_0$ l'image de $H^1_\infty(K, T_p(E))_{G_{\infty,K}}$ dans $H^1(K, T_p(E))$. Notons $\Sigma_f$ l'ensemble des caractères de $\Delta$ tels que $e_\delta H^1(K, V_p(E))_0$ soit contenu dans $e_\delta H^1_f(K, V_p(E)) = e_\delta(\mathbb{Q}_p \otimes \check{S}_p(E))$ et $\Sigma_{/f}$ l'ensemble des caractères de $\Delta$ tels que $e_\delta H^1(K, V_p(E))_0$ ne soit pas contenu dans $e_\delta H^1_f(K_p, V_p(E))$.

Si $M_1$ et $M_2$ sont deux $\mathbb{Z}_p$-modules de même rang et $\mathbb{Q}_p \otimes M_1 \xrightarrow{f} \mathbb{Q}_p \otimes M_2$ est un isomorphisme, on note $[M_2 : M_1] = [M_1 \xrightarrow{f} M_2]$ l'indice généralisé égal à $\frac{\sharp (M_2)_{tors}}{\sharp (M_1)_{tors}} [M_2/(M_2)_{tors} : M_1/(M_1)_{tors}]$.

On peut supposer que $F_{\mathfrak{c}}$ n'a pas de zéro en un caractère $\delta$. On a alors le résultat suivant sur $F_2$ ([18], même démonstration) :

### 2.2.3. Proposition. *On a*

$$\prod_\delta \delta(F_c^{-1}F_2^*) \sim [H^2_{\infty,p}(K,T_p(E))_{G_{\infty,K}} : H^2_{\infty,p}(K,T_p(E))^{G_{\infty,K}}]$$

$$\times [H^1_\infty(K,T_p(E))_{G_{\infty,K}} : \mathbb{Z}_p[\Delta_K]P_K(\mathfrak{c})]$$

$$\sim \frac{L(E/K_p,1)^{-1}|d_K|_p \operatorname{Tam}(E/K)\sharp\operatorname{III}(T_p(E)/K)}{\sharp E(K)^2_{tors}} \prod_{v|p} \sharp E(K_v)_{tors}$$

$$\times \operatorname{disc}\langle\langle\cdot,\cdot\rangle\rangle_K \mathcal{N}(P_K(\mathfrak{c}))^{-1} \prod_{\delta \in \Sigma_f} [e_\delta \mathcal{O}_{K_p}\omega_E^* : e_\delta \log_E \check{S}_p(E/K)_u]^2$$

*où $P_K(\mathfrak{c})$ est la projection de $\mathfrak{c}$ dans $H^1(K,T_p(E))$, où $\mathcal{N}(P_K(\mathfrak{c}))$ avec*

$$\mathcal{N}(P_K(\mathfrak{c})) = \prod_{\delta \in \Sigma_{/f}} [e_\delta \mathcal{O}_{K_p}\omega_E : e_\delta \mathbb{Z}_p[\Delta_K] \exp^*_{K_p} P_K(\mathfrak{c})]$$

$$\prod_{\delta \in \Sigma_f} [e_\delta \mathcal{O}_{K_p}\omega_E^* : e_\delta \mathbb{Z}_p[\Delta_K] \log_{K_p} P_K(\mathfrak{c})]$$

*et où $\check{S}_p(E/K)_u$ est l'intersection de $\check{S}_p(E/K)$ et de l'image de $H^1_\infty(\mathbb{Q},T_p(E))$.*

La première égalité est classique. La seconde se démontre à partir de la suite exacte de Poitou-Tate comme dans [18] ou [20] et en utilisant le lemme suivant qui modifie le lemme 3.6.6 de [20] et fait intervenir la différente $\mathcal{D}_K$ et le discriminant $d_K$ de $K$.

### 2.2.4. Lemme.

$$\operatorname{Tam}_p(E/K) = [H^1_f(K_p,T_p(E)) : \mathcal{O}_{K_p}\omega_E^*]L(E/K_p,1)$$
$$= \sharp E(K_p)_{tors}[\mathcal{D}_{K_p}^{-1}\omega_E : H^1_{/f}(K_p,T_p(E))]L(E/K_p,1)$$
$$= \sharp E(K_p)_{tors}|d_K|_p^{-1}[\mathcal{O}_{K_p}\omega_E : H^1_{/f}(K_p,T_p(E))]L(E/K_p,1)$$

*où $\omega_E^*$ est un élément de $D_p(E)$ tel que $[\omega_E,\omega_E^*] = 1$ dans la dualité naturelle, par exemple $\eta$ ou $v/[\omega_E,v]$.*

La première affirmation est la définition à la Bloch-Kato du nombre de Tamagawa ([5]). La deuxième se déduit des dualités :

$$[\mathcal{O}_{K_p}\omega_E : H^1_{/f}(K_p;T_p(E))] = \sharp E(K_p)^{-1}_{tors}[\mathcal{D}_{K_p}^{-1}\omega_E^* : H^1_f(K_p,T_p(E))]$$
$$= \sharp E(K_p)^{-1}_{tors}|d_K|_p[\mathcal{O}_{K_p}\omega_E^* : H^1_f(K_p,T_p(E))]$$
$$= \sharp E(K_p)^{-1}_{tors}|d_K|_p \operatorname{Tam}_p(E/K)L_p(E/K_p,1)^{-1}$$

On a $\operatorname{Tam}_p(E/K) = [E(K_p) : E^0(K_p)]$ où $E^0(K_p)$ est l'ensemble des points de $E(K_p)$ ne se réduisant pas sur le point singulier (s'il existe). Ici cela vaut donc 1.

La théorie du régulateur $p$-adique d'Iwasawa permet d'interpréter $\mathcal{N}(P_K(\mathfrak{c}))$ en termes de $\mathcal{L}_E(\mathfrak{c})$. Soit $G$ un élément de $\mathcal{H} \otimes D_p(E)$ tel que $(1-\varphi)G = \mathcal{L}_E(\mathfrak{c})\cdot(1+x)$. La proposition suivante est très importante :

### 2.2.5. Proposition. *Soit $\delta$ un caractère non trivial de $\Delta_K$. Alors,*

$$(p\varphi)^{-n}e_\delta G(\zeta_n - 1) = e_\delta \exp^*_{K_p}(P_K(\mathfrak{c}))$$

*et lorsque $e_\delta G(\zeta_n - 1) = 0$,*

$$D((p\varphi)^{-n}e_\delta(G))(\zeta_n - 1) = \frac{\left(\log^2_{\omega_E} e_\delta P_K(\mathfrak{c})\right) v - h_v(e_\delta P_K(\mathfrak{c}))\omega_E}{\log_{\omega_E} e_\delta P_K(\mathfrak{c})}$$

Ici, $D = (1+x)d/dx$, rappelons que $(1-p\varphi)D(G) = \mathcal{L}'_E(\mathfrak{c}) \cdot (1+x)$. Des formules analogues sont vraies pour un caractère non ramifié. Cette proposition est une conséquence essentielle de [18]. On commence par regarder le cas où le conducteur de $\delta$ est le conducteur de $K$. Il s'agit essentiellement des propositions 2.1.4, 2.2.2 et 2.3.5 (ou plutôt sa version pour un caractère ramifié) de [18]. Il suffit de vérifier que les valeurs contre $\nu$ par l'accouplement $[\cdot,\cdot]_{D_p(E)}$ (étendu par linéarité) des membres de droite et de gauche sont égales pour $\nu = \omega_E$ et $v \notin \text{Fil}^0 D_p(E)$; pour $\nu = \omega_E$, c'est la proposition 2.2.2 *loc.cit*; pour $\nu = v$, si $\beta = \beta_{\mathbb{Q}_p} v$ est la projection de $D_p(E)$ sur $\text{Fil}^0 D_p(E)$ parallèlement à $\mathbb{Q}_p v$, on a $[\beta x, v]_{D_p(E)} = [x,v]_{D_p(E)}$; comme $\beta$ est à valeurs dans $\text{Fil}^0$, on peut remplacer le second membre par $\log P_K(\mathfrak{c}) \equiv \log_{\omega_E} P_K(\mathfrak{c}) v/[v, \omega_E]_{D_p(E)} \mod \text{Fil}^0 D_p(E)$ pour obtenir la formule désirée. Lorsque $\delta$ se factorise par un quotient $\Delta_L$ de $\Delta_K$ avec $L$ de $p$-conducteur $p^{n_\delta}$, on utilise les formules suivantes pour se ramener au cas précédent :

$$e_\delta((p\varphi)^{-n}G(\zeta_n - 1)) = e_{\delta,L}(Tr_{K/K_{n_\delta}}((p\varphi)^{-n}G(\zeta_n - 1))) = e_{\delta,L}((p\varphi)^{-n_\delta} G(\zeta_{n_\delta} - 1))$$

avec $e_{\delta,L} = \sum_{\sigma \in \Delta_L} \delta(\sigma)^{-1} \sigma$.

Revenons à la démonstration principale. Posons $s_K = \zeta_n + \zeta_{n-1} + \ldots + \zeta_1$. On a alors $\mathcal{O}_K = \mathcal{A}_K s_K$ avec $\mathcal{A}_K$ l'ordre associé à $\mathcal{O}_K$. Remarquons à ce propos que $s_K$ est obtenu par la même recette que $G(\zeta_n - 1)$ : l'équation $(1 - \varphi)S = (1+x) - 1$ implique que $S(\zeta_n - 1) - \sum_{i=1}^{n} \zeta_i \in \mathbb{Z}_p$, c'est-à-dire que $S(\zeta_n - 1) - s_K \in \mathbb{Z}_p$. [2]On a $\mathcal{A}_K s_K = \mathcal{A}_K S(\zeta_n - 1)$. De même,

$$\varphi^{-n} G(\zeta_n - 1) - \sum_{i=1}^{n} \varphi^{-i} \mathcal{L}_E(\zeta_i - 1) \in D_p(E)$$

On a pour $\delta \in \Sigma_{/f}$,

$$p^{-n} \frac{[\delta(\varphi^{-n} \mathcal{L}_E(\mathfrak{c})), v]_{D_p(E)}}{[\omega_E, v]_{D_p(E)}} \sim \frac{[e_\delta((p\varphi)^{-n} G(\zeta_n - 1)), v]_{D_p(E)}}{[\omega_E, v]_{D_p(E)} e_{\ delta}(s_K)}$$
$$\sim [e_\delta \mathbb{Z}_p[\Delta_K] s_K \omega_E : e_\delta \mathbb{Z}_p[\Delta_K](p\varphi)^{-n} G(\zeta_n - 1)]$$
$$\sim [e_\delta \mathbb{Z}_p[\Delta_K] s_K \omega_E : e_\delta \mathbb{Z}_p[\Delta_K] \exp^* P_K(\mathfrak{c})]$$

et pour $\delta \in \Sigma_f$,

$$p^{-n} \frac{[\delta(\varphi^{-n} \mathcal{L}'_E(\mathfrak{c})), v]_{D_p(E)}}{[\omega_E, v]_{D_p(E)}} \sim \frac{[e_\delta(D((p\varphi)^{-n}G)(\zeta_n - 1)), v]_{D_p(E)}}{[\omega_E, v]_{D_p(E)} e_\delta(s_K)}$$
$$\sim \frac{h_v(e_\delta P_K(\mathfrak{c}))}{e_\delta(s_K) \log_{\omega_E} e_\delta P_K(\mathfrak{c})}$$
$$\sim [e_\delta \mathbb{Z}_p[\Delta_K] s_K \omega^* : e_\delta \mathbb{Z}_p[\Delta_K] \log_E P_K(\mathfrak{c})] \frac{h_v(e_\delta P_K(\mathfrak{c}))}{\log^2_{\omega_E} e_\delta P_K(\mathfrak{c})}$$
$$\sim \frac{[e_\delta \mathbb{Z}_p[\Delta_K] s_K \omega^* : e_\delta \mathbb{Z}_p[\Delta_K] \log_E P_{u,\delta}]}{[e_\delta \mathbb{Z}_p[\Delta_K] \omega^* : e_\delta \log_E e_\delta \check{S}_p(E/K)_u)]^2} h_v(P_{u,\delta})$$

pour $P_{u,\delta}$ un générateur de $e_\delta \check{S}_p(E/K)_u$). En prenant le produit sur tous les caractères $\delta$ (y compris les caractères non ramifiés pour lesquels nous n'avons pas écrit la formule) et en utilisant le fait dû essentiellement à Leopoldt que (voir par exemple [8], la manière de l'utiliser m'a été rappelée par [2])

$$|d_K|_p^{-1}[\mathcal{A}_K : \mathbb{Z}_p[\Delta_K]] = p^{n(p^{n-1}(p-1)-1)} \ ,$$

(le discriminant est à une unité près $\prod_\xi f_\xi$, l'indice est $\prod_\xi [K : \mathbb{Q}(\xi)]$ où $f_\xi$ est le conducteur de $\xi$ et $\mathbb{Q}(\xi)$ le corps de ses valeurs et on remarque alors que $f_\xi[K$ :

---

[2]Dans le cas où $K = K_0(\mu_{p^n})$ avec $K_0$ extension non ramifiée en $p$, il faut remplacer $(1+x)-1$ par $c((1+x) - 1)$ avec $c$ une base de $\mathcal{O}_{K_0}$ sur $\mathbb{Z}_p$

$\mathbb{Q}(\xi)] \sim p^n$ pour $\xi$ caractère non trivial), on obtient en notant $s_K$ le rang de $\check{S}_p(E/K) \cap H^1(K, T_p(E))_0$,

$$\frac{[\mathbf{1}(\mathcal{L}^{(s_K)}(\mathfrak{c}/K), \tilde{v}_K]_{D_p(E)}}{[\omega_E, v]^{[K:\mathbb{Q}_p]}} \sim |d_K|_p^{-1} \mathcal{N}(P_K(\mathfrak{c})) \prod_{\delta \in \Sigma_f} \frac{h_v(P_{u,\delta})}{[e_\delta \mathbb{Z}_p[\Delta_K]\omega^* : e_\delta \log_E e_\delta \check{S}_p(E/K)_u]^2} .$$

En multipliant cette formule par celle de la proposition 2.2.3, on en déduit le théorème.

2.2.6. **Remarque.** Si $\mathfrak{c} \in H^1_\infty(\mathbb{Q}, T_p(E))$ a une image non nulle $P_\mathbb{Q}(\mathfrak{c})$ dans $\mathbb{Q} \otimes H^1(\mathbb{Q}, T_p(E))$, alors il n'est pas divisible par $\gamma - 1$ dans $H^1_\infty(\mathbb{Q}, T_p(E))$ et $\mathcal{L}_E(\mathfrak{c})$ ne s'annule pas en $\mathbf{1}$ ou a un zéro d'ordre 1 en $\mathbf{1}$. Dans le premier cas, $P_\mathbb{Q}(\mathfrak{c})$ n'appartient pas à $\check{S}_p(E/\mathbb{Q}) = H^1_f(\mathbb{Q}, T_p(E))$ et $(1-\varphi)^{-1}(1-p^{-1}\varphi^{-1})\mathbf{1}(\mathcal{L}_E(\mathfrak{c})) \in \text{Fil}^0 D_p(E)$. Dans le second cas, $P_\mathbb{Q}(\mathfrak{c})$ appartient à $\check{S}_p(E/\mathbb{Q})$ et $(1-\varphi)^{-1}(1-p^{-1}\varphi^{-1})\mathbf{1}(\mathcal{L}'_E(\mathfrak{c}))$ n'appartient pas à $\text{Fil}^0 D_p(E)$. En effet un tel point ne peut pas être dans le noyau de localisation en $p$ ([16, lemme 4.5.1], la démonstration est exactement la même). Le même raisonnement s'applique en remplaçant $\mathbf{1}$ par un caractère $\delta$ de conducteur $p^{n+1}$ et l'endomorphisme d'Euler par $(p\varphi)^{n+1}$. Ainsi, l'annulation ou non de la composante modulo $\text{Fil}^0 D_p(E)$ est un élément important.

2.3. **Équation fonctionnelle.** Il est montré dans [20] que le module $\mathbb{I}_{arith}(E/\mathbb{Q})$ est invariant par l'involution $\iota : \tau \mapsto \tau^{-1}$ (ou $\rho \mapsto \rho^{-1}$). On peut alors construire un élément $c$ de $H^1(\langle \iota \rangle, \mathbb{Z}_p[[\Gamma]]^*)$ donné par la classe du cocycle déterminé par $a_\iota \in \mathbb{Z}_p[[\Gamma]]^*$ avec $e_0 I^\iota = a_\iota e_0 I$ si $I$ est un générateur de $\mathbb{I}_{arith}(E/\mathbb{Q})$. On vérifie facilement que $H^1(\langle \iota \rangle, \Lambda^*)$ est d'ordre 2 et a comme éléments les classes de 1 et de $-1$. On en déduit qu'il existe un unique élément $\epsilon_{arith} \in \{\pm 1\}$ et un générateur $I_{arith}(E/\mathbb{Q})$ de $\mathbb{I}_{arith}(E/\mathbb{Q})$ tels que

$$I_{arith}(E/\mathbb{Q})(\rho^{-1}) = \epsilon_{arith} I_{arith}(E/\mathbb{Q})(\rho) .$$

Cette construction est due à Greenberg. Notons $r'_{arith}$ la multiplicité du zéro en $\mathbf{1}$ de $\mathbb{I}_{arith}(E/\mathbb{Q})$.

**Proposition.** *On a $\epsilon_{arith} = (-1)^{r'_{arith}}$.*

On a $r_{arith} \leq r'_{arith}$. Il est clair que $r_{arith} = 0$ si et seulement si $r'_{arith} = 0$. Ainsi, si $r'_{arith} = 1$, on a $r_{arith} = 1$. Réciproquement, si $r_{arith} = 1$ et si $r'_{arith} > 1$, $\check{S}_p(E)$ est contenu dans le noyau de localisation en $p$ (regarder la composante modulo $\text{Fil}^0 D_p(E)$). Cela n'est pas possible si $E(\mathbb{Q})$ est infini. Ainsi, si $E(\mathbb{Q})$ est de rang 1 et $\mathrm{III}(E/\mathbb{Q})(p)$ fini, $r_{arith} = r'_{arith} = 1$.

**Remarque.** Dans le cas ordinaire, Greenberg [9] montre par un très joli argument dû à Guo que

$$r'_{arith} \equiv r_{arith} \mod 2 .$$

3. Conjecture principale et théorème de Kato

3.1. **Conjecture principale.** La conjecture principale pour $E$ et $p$ telle qu'elle est énoncée dans [18] dit :

3.1.1. **Conjecture.** *$L_p(E)$ est un générateur du $\mathbb{Z}_p[[G_\infty]]$-module $\mathbb{I}^{arith}(E)$ :*

$$\mathbb{I}^{arith}(E) = \mathbb{Z}_p[[G_\infty]] L_p(E)$$

Le théorème fondamental de Kato qui s'appuie sur la technique des systèmes d'Euler introduite par Kolyvagin est une avancée fondamentale pour la démonstration de cette conjecture. Il y a deux versions du résultat de Kato. La première est énoncée dans [22] :

**3.1.2. Théorème** (Kato). *Il existe un entier $r$ tel que $p^r L_p(E) = g I^{arith}$ avec $g \in \mathbb{Z}_p[[G_\infty]]$.*

Soit $\rho_p$ la représentation $p$-adique de $G_\mathbb{Q}$ donnant l'action sur les points de $p$-torsion : $\rho_p : G_\mathbb{Q} \to GL_2(\mathbb{Z}/p\mathbb{Z})$.

**3.1.3. Théorème** (Kato). *Si $E$ a multiplication complexe ou si $\rho_p(G_\mathbb{Q}) = GL_2(\mathbb{Z}/p\mathbb{Z})$, alors $L_p(E) = g I^{arith}$ avec $g \in \mathbb{Z}_p[[G_\infty]]$.*

Le passage de la formulation de Kato à celle-là est faite dans [18]. Pour déterminer si $\rho_p(G_\mathbb{Q}) = GL_2(\mathbb{Z}/p\mathbb{Z})$, les critères donnés par Serre dans [23] sont extrêmement commodes. Nous supposons toujours que $E$ a bonne réduction supersingulière en $p$, ce qui diminue le nombre de cas à envisager. Dans les exemples que nous étudierons, la surjectivité peut être montrée dans le cas sans multiplication complexe à l'aide des critères suivants (A.1)

**3.1.4. Proposition** (Serre). *Supposons que $E$ a réduction supersingulière en $p$.*
  1) *Pour $p = 3$, $\rho_3(G_\mathbb{Q}) = GL_2(\mathbb{Z}/3\mathbb{Z})$ si et seulement si $\Delta$ n'est pas un cube.*
  2) *Pour $p \geq 5$, si $N_E$ est sans facteurs carrés, $\rho_p(G_\mathbb{Q}) = GL_2(\mathbb{Z}/p\mathbb{Z})$.*
  3) *S'il existe un nombre premier $\ell$ divisant strictement $N_E$ tel que $\mathrm{ord}_\ell(j_E) \not\equiv 0 \bmod p$, alors $\rho_p(G_\mathbb{Q}) = GL_2(\mathbb{Z}/p\mathbb{Z})$.*
  4) *S'il existe un nombre premier $\ell$ tel que $a_l \not\equiv 0 \bmod p$, $a_l^2 \not\equiv 4\ell \bmod p$, $(\frac{a_l^2 - 4\ell}{p}) = 1$, alors $\rho_p(G_\mathbb{Q}) = GL_2(\mathbb{Z}/p\mathbb{Z})$.*

Sous l'hypothèse que $E$ a bonne réduction supersingulière en $p$, si $\rho_p(G_\mathbb{Q}) \neq GL_2(\mathbb{Z}/p\mathbb{Z})$, il est nécessairement égal au normalisateur d'un sous-groupe de Cartan non déployé dans $GL_2(\mathbb{Z}/p\mathbb{Z})$.

**Remarque.** La première courbe dans la liste de Cremona sans multiplication complexe, ayant bonne réduction supersingulière en 3 et telle que $\rho_3(G_\mathbb{Q}) \neq GL_2(\mathbb{Z}/3\mathbb{Z})$ est la courbe "1952C" de conducteur 1952, d'équation $y^2 = x^3 - 332x + 2752$, d'invariant $j = -988047936/226981 = (-996/61)^3$ et de discriminant $-929714176 = (-16 \times 61)^3$. L'image de $G_\mathbb{Q}$ dans $PGL_2(\mathbb{Z}/3\mathbb{Z}) \cong S_4$ est le groupe de Galois du polynôme $3x^4 - 1992x^2 - 21927936$ qui est le groupe diédral d'ordre 8. Il en est de même de la courbe "1952D" d'équation $y^2 = x^3 - 332x - 2752$. Ces deux courbes deviennent isomorphes sur $\mathbb{Q}(\sqrt{-1})$. Pour la courbe "1044A" d'équation $y^2 = x^3 - 3105x - 139239$, en prenant $l \leq 100000$, la proposition ne permet pas de montrer que $\rho_5$ est surjective. En utilisant le logiciel de calcul Magma, on vérifie que la clôture galoisienne du corps engendré par la $x$-coordonnée d'un point de 5-torsion (de degré 12) est d'ordre 24 et donc que $\rho_5$ n'est pas surjective.

Nous allons maintenant donner quelques conséquences du théorème de Kato, toujours dans le cas de bonne réduction supersingulière.

**3.2. Sans hypothèse sur le rang.** Notons $r'_{arith}$ l'ordre du zéro dans $I^{arith}$ et $r_{anal}$ l'ordre du zéro dans $L_p(E)$ en **1**.

**3.2.1. Proposition** (Kato). *On a les inégalités*

$$\mathrm{rg}_{\mathbb{Z}_p} \check{S}_p(E) \leq r'_{arith} \leq r_{anal}$$

D'après le théorème 2.2.1, on a $\mathrm{rg}_{\mathbb{Z}_p} \check{S}_p(E)) \leq r'_{arith}$ avec égalité dans le cas de non dégénérescence de $\langle\langle \cdot, \cdot \rangle\rangle$ sur $\check{S}_p(E)_0 = H^1_f(\mathbb{Q}, T_p(E))_0$.

**3.2.2. Proposition.** *Supposons que la conjecture de Birch et Swinnerton-Dyer $p$-adique soit vraie pour $E/\mathbb{Q}$ (resp. pour $E/\mathbb{Q}(\mu_p)$) à une unité $p$-adique près. Alors, la conjecture principale est vraie pour $E/\mathbb{Q}$ (resp. pour $E/\mathbb{Q}(\mu_p)$) et $p$.*

3.3. **Cas où** $r_{anal} = 0$. Lorsque $L(E,1) \neq 0$, c'est-à-dire $\mathbf{1}(L_p(E)) \neq 0$, $E(\mathbb{Q})$ et $\mathbf{III}(E/\mathbb{Q})(p)$ sont finis. Kato montre en fait que $\mathbf{III}(E/\mathbb{Q})$ est fini (ce résultat avait été obtenu auparavant par Kolyvagin). De plus, $(1-\varphi)^{-1}(1-p^{-1}\varphi^{-1})\mathbf{1}(L_p(E))$ appartient à $\mathrm{Fil}^0 D_p(E)$ et vaut $\frac{L_{\{p\}}(E,0)}{\Omega_E}\omega_E$.

3.3.1. **Proposition** (Kato). *Si $\rho_p$ est surjective et si $\mathbf{1}(L_p(E)) \neq 0$,*

$$\mathrm{ord}_p \mathbf{III}(E/\mathbb{Q})(p) \leq \mathrm{ord}_p \frac{L(E,1)}{\Omega_E} - \mathrm{ord}_p \frac{\mathrm{Tam}(E)}{\sharp E(\mathbb{Q})_{tor}^2} \ .$$

*Lorsqu'il y a égalité, la conjecture principale est vraie pour le caractère trivial de $\mathrm{Gal}(\mathbb{Q}(\mu_p)/\mathbb{Q})$, c'est-à-dire*

$$e_0 I^{arith}(E) = \Lambda e_0 L_p(E) \ .$$

**Exemple.** Pour les courbes $X_0(17)^{(D)}$ avec $D$ égal à une de valeurs suivantes, le second membre dans la proposition 3.3.1 est égal à 2 pour $p = 3$ : $-947$, $-923$, $-907$, $-827$, $-823$, $-691$, $-635$, $-503$, $-488$, $-479$, $-419$, $-347$, $-311$, $-283$, $-199$, $-167$, $-139$, $253$, $373$, $457$, $557$, $701$, $749$, $917$, $953$.

3.4. **Zéro d'ordre $\geq 1$.** Si $f$ est une fonction analytique sur $G_\infty$ et $\rho$ un caractère de $G_\infty$, notons $f^*(\rho)$ ou $\rho(f^*)$ le coefficient dominant en $\rho$ dans le développement de $f(\rho\langle\chi\rangle^s)$. On pose

$$L'_{p,\omega_E}(E,\mathbf{1}) = [(1-\varphi)^{-1}(1-p^{-1}\varphi^{-1})\mathbf{1}(L'_p(E)), \omega_E]_{D_p(E)}$$
$$= [\mathbf{1}(L'_p(E)), \tilde{\omega}_E]_{D_p(E)}$$
$$L^{(r)}_{p,\omega_E}(E,\mathbf{1}) = [\mathbf{1}(L^{(r)}_p(E)), \tilde{\omega}_E]_{D_p(E)}$$
$$L^*_{p,\omega_E}(E,\mathbf{1}) = [\mathbf{1}(L^{(*)}_p(E)), \tilde{\omega}_E]_{D_p(E)} \ .$$

3.4.1. **Proposition.** *On suppose $r_{anal} \geq 1$. Alors, si $L^*_{p,\omega_E}(E,\mathbf{1}) \neq 0$, le rang de $\check{S}_p(E)$ est supérieur ou égal à 1.*

La condition $L^*_{p,\omega_E}(E,\mathbf{1}) \neq 0$ signifie que $(1-\varphi)^{-1}(1-p^{-1}\varphi^{-1})\mathbf{1}(L^*_p(E))$ n'appartient pas à $\mathrm{Fil}^0 D_p(E)$. Les hypothèses sont vérifiables numériquement.

*Démonstration.* Supposons que $\mathbf{1}(I_{arith}) \neq 0$. Comme $(1-\varphi)^{-1}(1-p^{-1}\varphi^{-1})\mathbf{1}(I_{arith})$ est proportionnel à $\omega_E$, il en est de même de $(1-\varphi)^{-1}(1-p^{-1}\varphi^{-1})\mathbf{1}(L^*_p(E))$. Ce qui est contradictoire avec l'hypothèse, donc $r'_{arith}$ et $r_{arith}$ sont supérieurs à 1. Ainsi, $\mathbf{1}(I_{arith}) = 0$ et le rang de $\check{S}_p(E)$ est supérieur ou égal à 1. □

Cette proposition se généralise à un caractère non trivial. Si $K$ est une extension abélienne de $\mathbb{Q}$, on pose

$$L^*_{p,\omega_E}(E/K,\mathbf{1}) = [\mathbf{1}(L^*_p(E/K)), \tilde{\omega}_{E,K}]_{D_p(E)}$$

et si $\delta$ est un caractère de conducteur $p^{n_\delta+1}$ non trivial,

$$L^*_{p,\omega_E}(E,\delta) = [(p\varphi)^{-(n_\delta+1)}\delta(L^*_p(E)), \omega_E]_{D_p(E)}$$

**Proposition.** *Supposons que $\delta(L_p(E))$ est nul pour un caractère $\delta$ d'ordre $p^m$, de conducteur $p^{m+1}$. Alors, si $L^*_{p,\omega_E}(E,\delta)$ est non nul,*

$$\mathrm{rg}_{\mathbb{Z}_p} \check{S}_p(E/\mathbb{Q}_m) - \mathrm{rg}_{\mathbb{Z}_p} \check{S}_p(E/\mathbb{Q}_{m-1}) \geq (p-1)p^{m-1} \ .$$

**Exemple.** Soit $E$ la courbe $y^2+y = x^3-4x+2$ de conducteur 1909 (c'est la courbe 1909$A$). Le nombre premier $p = 3$ est supersingulier. On a $L(E,1) = 0$, le signe de l'équation fonctionnelle est 1. Donc, $r_{anal} \geq 2$. Le calcul montre que $L^{(2)}_{p,\omega_E}(E,\mathbf{1})$ est non nul modulo $9M_E$. Donc, le rang de $\check{S}_p(E)$ est 1 ou 2.

3.4.2. **Conjecture.** *Si $L(E,1) = 0$ (c'est-à-dire si $r_{anal} \geq 1$), $L^*_{p,\omega_E}(E,\mathbf{1}) \neq 0$. Si $\delta$ est un caractère d'ordre $p^m$, de conducteur $p^{m+1}$ et si $\delta(L_p(E)) = 0$, alors $L^*_{p,\omega_E}(E,\delta) \neq 0$.*

D'une part, tous les exemples numériques calculés confirment cette conjecture, d'autre part nous avons montré que si $E(\mathbb{Q})$ est infini, $L^*_{p,\omega_E}(E,\mathbf{1}) \neq 0$. De plus, si la conjecture est fausse par exemple pour le caractère trivial, $\text{III}(E/\mathbb{Q})(p)$ serait infini ; en effet, cela impliquerait (§2.2.6) que $\check{S}_p(E/\mathbb{Q})$ serait contenu dans le noyau de localisation en $p$ et donc que $\check{S}_p(E/\mathbb{Q})$ serait de rang $\geq 0$. Comme l'intersection de $E(\mathbb{Q})$ et du noyau de localisation est de torsion, $\text{III}(E/\mathbb{Q})(p)$ serait infini. L'énoncé correspondant à des caractères non triviaux peut s'écrire sous la forme :

3.4.3. **Conjecture.** *Soit $n > 0$. Supposons que pour $r > 0$, $\hat{L}^{(k)}_{p,(0)}(E) \equiv 0 \mod \xi_n$ pour $k < r$ et $\hat{L}^{(r)}_{p,(0)}(E) \not\equiv 0 \mod \xi_n$. Alors $[\hat{L}^{(r)}_{p,(0)}(E), \varphi^{-n-1}\omega_E]_{D_p(E)} \not\equiv 0 \mod \xi_n$.*

Remarquons que la conjecture est un énoncé sur une propriété des symboles modulaires.

On déduit encore du théorème de Kato et des résultats du paragraphe 2.2 le fait suivant :

3.4.4. **Proposition.** *Les pentes $\lambda_{arith}$ de $(1 - p^{-1}\varphi^{-1})(1-\varphi)^{-1}\mathbf{1}(I^*_{arith})$ et $\lambda_{anal}$ de $(1 - p^{-1}\varphi^{-1})(1-\varphi)^{-1}\mathbf{1}(L^*_p(E))$ dans une base $(\nu, \omega_E)$ de $D_p(E)$ sont égales. Si $E(\mathbb{Q})$ est infini, il existe un point $P$ de $\check{S}_p(E/\mathbb{Q})$ tel que*

$$\lambda_{anal} = -\frac{h_\nu(P)}{\log^2_{\omega_E} P} \ .$$

En effet, ces pentes sont aussi égales à celle de $\mathcal{L}_E(\mathfrak{c})$ pour $\mathfrak{c}$ un élément non nul de $H^1_\infty(\mathbb{Q}, T_p(E))$ et engendrant un module d'indice premier à $\gamma - 1$. Lorsque $E(\mathbb{Q})$ est infini, $\mathbf{1}(\mathcal{L}_E(\mathfrak{c})) = 0$ et sa dérivée a comme pente $-\frac{h_\nu(P)}{\log^2_{\omega_E} P}$ pour $P = P_\mathbb{Q}(\mathfrak{c}) \in \check{S}_p(E/\mathbb{Q})$.

On poeut voir cette proposition comme un moyen numérique de vérifier que le rang est $\geq 2$ connaissant un point d'ordre infini. En effet, lorsque $E(\mathbb{Q})$ est de rang 1, la valeur de $-\frac{h_\nu(P)}{\log^2_{\omega_E} P}$ est indépendant de $P \in E(\mathbb{Q})$. Si par chance on a trouvé un point $P \in E(\mathbb{Q})$ tel que

$$\lambda_{anal} \neq -\frac{h_\nu(P)}{\log^2_{\omega_E} P} \ ,$$

$\check{S}_p(E/\mathbb{Q})$ est de rang $\geq 2$. Cela peut se généraliser à un rang plus grand. Supposons trouvés $r$ points indépendants de $E(\mathbb{Q})$. Si le régulateur $p$-adique (vectoriel) de ces points (au sens de 1.3) n'est pas proportionnel à $(1 - p^{-1}\varphi^{-1})(1-\varphi)^{-1}\mathbf{1}(L^*_p(E))$, le rang de $\check{S}_p(E)$ est strictement plus grand que $r$.

Remarquons cependant qu'une base naturelle dans laquelle s'exprime le régulateur $p$-adique est $(\omega_E, \eta)$ de $D_p(E)$ alors que les valeurs de la fonctions $L$ $p$-adiques s'expriment naturellement dans la base $\omega_E, \varphi\omega_E$. On a ainsi besoin de calculer $\varphi\omega_E$ dans la base $(\eta, \omega_E)$ (programme `frobenius`, [26])

3.5. **Zéro d'ordre 1.**

3.5.1. **Proposition.** *Supposons que $L(E,1) = 0$ (et donc $\mathbf{1}(L_p(E)) = 0$). Si $L'_{p,\omega_E}(E,\mathbf{1}) \neq 0$, (en particulier, $\mathbf{1}(L'_p(E)) \neq 0$), $H^1_f(\mathbb{Q}, V_p(E))$ est de dimension 1.*

*Démonstration.* Le fait que $\mathbf{1}(L'_p(E)) \neq 0$ implique que $\mathbf{1}(I'_{arith}) \neq 0$. On en déduit que le rang de $\tilde{S}_p(E)$ est inférieur ou égal à 1. On a vu précédemment que ce rang est supérieur ou égal à 1. D'où la proposition. $\square$

Sous l'hypothèse que $L'(E,1) \neq 0$, Kolyvagin a montré beaucoup plus : le rang de $E(\mathbb{Q})$ est 1 et $\text{III}(E/\mathbb{Q})$ est fini.

**3.5.2. Proposition.** 1) *Supposons que $\mathbf{1}(L_p(E)) = 0$. Alors, le point $P_{Kato}$ de $H^1(\mathbb{Q}, T_p(E))$ construit par Kato appartient à $\check{S}_p(E/\mathbb{Q})$ et $L_p(E)$ a un zéro d'ordre 1 en $\mathbf{1}$ si et seulement si $P_{Kato}$ n'est pas de torsion.*

*2) Supposons que $\mathbf{1}(L_p(E)) = 0$ et $L'_{p,\omega_E}(E, \mathbf{1}) \neq 0$. Si le rang de $E(\mathbb{Q})$ est non nul, le point $P_{Kato}$ de $H^1(\mathbb{Q}, T_p(E))$ construit par Kato appartient à $\mathbb{Z}_p \otimes_{\mathbb{Z}} E(\mathbb{Q})$ et est non nul. Si la conjecture $p$-adique de Birch et Swinnerton-Dyer pour $E/\mathbb{Q}$ est vraie, $P_{Kato}$ n'appartient pas à $\mathbb{Q} \otimes_{\mathbb{Z}} E(\mathbb{Q})$.*

*Démonstration.* Lorsque $\mathbf{1}(L_p(E)) = 0$, l'exponentielle duale du point de Kato $P_{Kato}$ est nulle, ce qui signifie que $P_{Kato}$ appartient à $\check{S}_p(E) = H^1_f(\mathbb{Q}, T_p(E))$ et vérifie $(1-\varphi)^{-1}(1-p^{-1}\varphi^{-1})\mathbf{1}(L'_p(E)) \equiv \log P_{Kato} \mod \text{Fil}^0 D_p(E)$ d'après [18]. Si $P_{Kato}$ n'est pas de torsion, $\mathbf{1}(L'_p(E))$ est non nul. Réciproquement, si $P_{Kato}$ est de torsion, comme $P_{Kato}$ est la projection de $\mathfrak{c} \in H^1_\infty(\mathbb{Q}, T_p(E))$ et que $\mathbb{Q}_p \otimes H^1_\infty(\mathbb{Q}, T_p(E))_{G_\infty}$ s'injecte dans $H^1_f(\mathbb{Q}, V_p(E))$, $\mathfrak{c} = (\gamma - 1)\mathfrak{c}'$ et $L_p(E) = \mathcal{L}_E(\mathfrak{c})$ est nul en $\mathbf{1}$.

Supposons maintenant que $\mathbf{1}(L'_p(E))$ est non nul ; si le rang de $E(\mathbb{Q})$ est non nul, il est nécessairement égal à 1 et $\text{III}(E/\mathbb{Q})(p)$ est fini. On en déduit que $\check{S}_p(E) = \mathbb{Z}_p \otimes_{\mathbb{Z}} E(\mathbb{Q})$ et donc que $P_{Kato}$ appartient à $\mathbb{Z}_p \otimes_{\mathbb{Z}} E(\mathbb{Q})$. Soit $P$ un générateur de $E(\mathbb{Q})$ (modulo torsion). Soit $m$ un entier $p$-adique tel que $P_{Kato} \equiv mP$. Sous la conjecture $p$-adique de Birch et Swinnerton-Dyer, on a alors

$$\log_E P_{Kato} \equiv \frac{\text{Tam}(E)}{\sharp E(\mathbb{Q})^2_{tor}} \sharp \text{III}(E/\mathbb{Q}) \log_{\omega_E}(P)^2 \omega_E^* \mod \text{Fil}^0 D_p(E)$$

où $\omega_E^*$ est une base duale de $\omega_E$ (on a par définition $\log_E P_{Kato} \equiv \log_{\omega_E} P_{Kato} \omega_E^* \mod \text{Fil}^0 D_p(E)$). D'où,

$$m = \frac{\text{Tam}(E)}{\sharp E(\mathbb{Q})^2_{tor}} \sharp \text{III}(E/\mathbb{Q}) \log_{\omega_E}(P) .$$

Mais $\frac{\text{Tam}(E)}{\sharp E(\mathbb{Q})^2_{tor}} \sharp \text{III}(E/\mathbb{Q})$ est un rationnel. Il ne reste plus qu'à vérifier que, d'après [4], si $P \in E(\mathbb{Q})$, $\log_{\omega_E} P$ n'appartient pas à $\mathbb{Q}$. $\square$

**Remarque.** Comme dans le cas des points de Heegner, si $\mathbf{1}(L_p(E)) = 0$, $L_p(E)$ a un zéro d'ordre 1 si et seulement si $P_{Kato}$ n'est pas de torsion. Mais, ici le point de Kato n'est pas dans $E(\mathbb{Q})$, mais seulement dans $\mathbb{Q}_p \otimes E(\mathbb{Q})$. Par contre, si les conjectures de Bloch-Kato $p$-adiques sont vraies pour $E$ en un entier $k > 1$, les points de Kato correspondant semblent être des multiples entiers d'un point motivique.

**3.5.3. Proposition.** *Supposons $r_{anal} = r_{arith} = 1$ et $E(\mathbb{Q})$ de rang 1. Alors,*

*Autrement dit, on a la majoration (rappelons que $\sharp E(\mathbb{Q})$ est premier à $p$)*
$\text{ord}_p(\text{III}(E/\mathbb{Q})(p)) \leq \text{ord}_p \text{Tam}(E)$
$$+ \text{ord}_p \left( L(E/\mathbb{Q}_p, 1)^{-1} \frac{(1-\varphi)^{-1}(1-p^{-1}\varphi^{-1})\mathbf{1}(L'_p(E))}{R^{(1)}(E)} \right)$$

*(le quotient de deux vecteurs liés est par définition leur coefficient de proportionnalité). On a $L(E/\mathbb{Q}_p, 1)^{-1} = \frac{p+1-a_p}{p}$. Pour le calcul, il suffit de regarder la composante modulo $\omega_E$, ce qui permet de remplacer le calcul de hauteur $p$-adique par un simple calcul de logarithme $p$-adique sur la courbe elliptique : avec*

$$L'_{p,\omega_E}(E, \mathbf{1}) = [(1-\varphi)^{-1}(1-p^{-1}\varphi^{-1})\mathbf{1}(L'_p(E)), \omega_E]_{D_p(E)} ,$$

*comme $[\omega_E, \varphi\omega_E]$ est une unité p-adique, on obtient que*

$$\mathrm{ord}_p(\mathrm{III}(E/\mathbb{Q})(p)) \leq \mathrm{ord}_p \mathrm{Tam}(E) + \mathrm{ord}_p \frac{L'_{p,\omega_E}(E,\mathbf{1})}{p(\frac{\log_{\omega_E}((p+1-a_p)P)}{p})^2}$$

*ou encore*

$$\mathrm{ord}_p(\mathrm{III}(E/\mathbb{Q})(p)) \leq \mathrm{ord}_p \mathrm{Tam}(E) + \frac{Z_2}{(\frac{\log_{\omega_E}((p+1-a_p)P)}{p})^2}$$

*si*

$$(1-\varphi)^{-1}(1-p^{-1}\varphi^{-1})\mathbf{1}(L'_p(E)) = \frac{Z_1\omega_E - Z_2 p\varphi\omega_E}{p+1-a_p}$$

avec $\mathbf{1}(L'_p(E)) = X\omega_E - Yp\varphi\omega_E$, $Z_1 = (p-1-2a_p+a_p^2/p)X + (2p-a_p/p)Y$, $Z_2 = (p-1)Y - (2-a_p/p)X$.

On peut reprendre les calculs qui ont été faits dans [3]. Par exemple, la 3-composante du groupe de Shafarevich-Tate de la courbe $X_0(17)^{(5)}$ d'équation $y^2 = x^3 - 15x^2 - 200x - 110000$ ou de la courbe $X_0(17)^{(-4)}$ d'équation $y^2 = x^3 + 12x^2 - 128x + 56320$ (dont le groupe de Mordell-Weil est de rang 1) est triviale. La 19-composante de la courbe $X_0(11)^{(8)}$ d'équation $y^2 = x^3 - 32x^2 - 10240x - 647168$ est triviale. De plus, pour ces courbes et nombres premiers correspondants, la conjecture principale p-adique est vraie.

**Exemple.** Soit la courbe $E = 43A$ d'équation $y^2 + y = x^3 + x^2$. Le nombre premier 7 est supersingulier. Le groupe de Mordell-Weil $E(\mathbb{Q})$ est de rang 1 engendré par $P = (0,0)$. Le logarithme 7-adique de $8P = (\frac{11}{49}, \frac{20}{343})$ est congru à 28 mod $7^2$. Comme

$$\begin{aligned}L(E/\mathbb{Q}_7,1)(1-\varphi)^{-1}&(1-7^{-1}\varphi^{-1})\mathbf{1}(L'_7(E)) \\ &= (5\times 7 + 6\times 7^2 + 4\times 7^3 + 4\times 7^4 + O(7^5))\omega_E \\ &\quad - 7(3\times 7 + 4\times 7^2 + 3\times 7^3 + 5\times 7^4 + O(7^5))\varphi\omega_E ,\end{aligned}$$

on en déduit que $\mathrm{III}(E/\mathbb{Q})(7)$ est trivial et que la conjecture 7-adique principale est vraie. Plus précisément, on trouve par le calcul la même valeur pour le régulateur de $P$. Comme $\mathrm{Tam}(E/\mathbb{Q}) = 1$ et que $E(\mathbb{Q})$ n'a pas de torsion, ces calculs sont ainsi compatibles avec la conjecture de Birch et Swinnerton-Dyer p-adique.

**Remarque.** Malheureusement dans le cas supersingulier, on ne sait pas relier l'ordre du zéro de $L_p(E)$ avec celui de $L(E,s)$. Supposons que $E$ a bonne réduction ordinaire en $p$. Lorsqu'une certaine composante $L_p^\pi(E)$ de la fonction $L$ p-adique vectorielle (la fonction $L$ p-adique de Mazur-Swinnerton-Dyer) a un zéro simple en $\mathbf{1}$, il est démontré que $L(E,s)$ a aussi un zéro simple en 1. La réciproque est vraie lorsqu'on sait montrer que la hauteur p-adique associée ne s'annule pas sur $E(\mathbb{Q})$. On obtient alors la divisibilité

$$\sharp\mathrm{III}(E/\mathbb{Q})(p) \mid \frac{L'(E/\mathbb{Q},1)}{\Omega_E h_\infty(P)} .$$

où $P$ est un générateur de $E(\mathbb{Q})$ modulo torsion et $h_\infty(P)$ la hauteur de Néron-Tate (cela utilise pêle-mêle le théorème de Gross-Zagier et son analogue p-adique, le théorème de Kato dans le cas ordinaire et les calculs des valeurs spéciales de la fonction $L$ arithmétique analogues à ceux de 2.2.1).

# 4. $\lambda$-INVARIANTS ET $\mu$-INVARIANTS

4.1. *Nous reprenons ici les idées de Kurihara [11] pour définir les $\lambda$-invariants et les $\mu$-invariants d'une classe d'éléments de $\mathcal{H}(G_\infty) \otimes D_p(E)$. Nous avons récemment appris qu'un travail analogue a été fait par Pollack [21] dans le cas où $a_p = 0$.*

Posons $\omega_n(x) = (1+x)^{p^n} - 1$ et $\xi_n(x) = \omega_n(x)/\omega_{n-1}(x)$ le polynôme cyclotomique de degré $p^n$.

**4.1.1. Lemme.** *Soit $F \in \mathcal{H} \otimes D_p(E)$ tel que $\varphi^{-n-1} F(\zeta_n - 1) \in \mathbb{C}_p \otimes \mathrm{Fil}^0 D_p(E)$ pour tout entier $n \geq 1$. Alors, il existe une et une seule famille de polynômes $P_n$ pour $n \geq 0$, de degré $< p^n$ vérifiant*
$$F \equiv P_n \varphi^{n+1} \omega_E - \xi_n P_{n-1} \varphi^{n+2} \omega_E \mod \omega_n(x) D_p(E)$$
*pour $n \geq 1$. On a alors la relation pour $n \geq 2$*
$$P_n - a_p P_{n-1} + \xi_{n-1} P_{n-2} \equiv 0 \mod \omega_{n-1}(x) \ .$$
*Si de plus*
$$(1 - p^{-1} \varphi^{-1})(1 - \varphi)^{-1} F(0) \in \mathrm{Fil}^0 D_p(E)$$
*(resp. $F(0) \in \varphi \mathrm{Fil}^0 D_p(E)$), on a*
$$(a_p - 2) P_1 \equiv ((a_p - 2) a_p - (p-1)) P_0 \mod \omega_0(x)$$
*(resp. $P_1 \equiv a_p P_0 \mod \omega_0(x)$). Enfin, si $F$ est d'ordre $\leq 0$, il existe une constante $M$ telle que $P_n \in M^{-1} \mathbb{Z}_p[x]$ pour tout $n$.*

*Démonstration.* Pour $n \geq 0$, écrivons $F \equiv P_n \varphi^{n+1} \omega_E - Q_n \varphi^{n+2} \omega_E \mod \omega_n(x)$ avec $P_n$ et $Q_n$ des polynômes de degré $< p^n$ à coefficients dans $\mathbb{Q}_p$. Pour $n \geq 1$, $\varphi^{-n-1} F(\zeta_n - 1) \in \mathbb{C}_p \omega_E$. Donc, $Q_n(\zeta_n - 1) = 0$ c'est-à-dire que $Q_n$ est divisible par $\xi_n$ : $Q_n = \xi_n R_{n-1}$ pour $n \geq 1$. Ici, $R_{n-1}$ est défini modulo $\omega_{n-1}(x)$ (mais la condition sur les degrés le détermine uniquement). On a la relation de récurrence pour $n \geq 2$
$$P_n \varphi^{n+1} \omega_E - \xi_n R_{n-1} \varphi^{n+2} \omega_E \equiv P_{n-1} \varphi^n \omega_E - \xi_{n-1} R_{n-2} \varphi^{n+1} \omega_E \mod \omega_{n-1}(x) \ .$$
En utilisant $\varphi^2 \omega_E - p^{-1} a_p \varphi \omega_E + p^{-1} \omega_E = 0$,
$$P_n \varphi^{n+1} \omega_E - \xi_n R_{n-1} \varphi^{n+2} \omega_E = P_n \varphi^{n+1}(\omega_E) - \xi_n R_{n-1} \varphi^n (p^{-1} a_p \varphi \omega_E - p^{-1} \omega_E)$$
$$= (P_n - \xi_n R_{n-1} p^{-1} a_p) \varphi^{n+1}(\omega_E) + p^{-1} \xi_n R_{n-1} \varphi^n(\omega_E)$$
$$\equiv R_{n-1} \varphi^n(\omega_E) + (P_n - R_{n-1} a_p) \varphi^{n+1}(\omega_E) \mod \omega_{n-1}(x)$$
car $\xi_n \equiv p \mod \omega_{n-1}(x)$. On en déduit que
$$P_n - a_p R_{n-1} \equiv -\xi_{n-1} R_{n-2} \mod \omega_{n-1}(x)$$
$$P_{n-1} \equiv R_{n-1} \mod \omega_{n-1}(x) \ .$$
Donc, pour $n \geq 1$, $R_n = P_n$ et
$$P_{n+1} - a_p R_n + \xi_n R_{n-1} \equiv 0 \mod \omega_n(x) \ .$$
La relation de récurrence pour $n = 1$ donne
$$P_0 \omega_E - Q_0 \varphi \omega_E \equiv P_1 \varphi \omega_E - Q_1 \varphi^2 \omega_E \mod \omega_0(x)$$
$$\equiv (P_1 - Q_1 p^{-1} a_p) \varphi(\omega_E) + p^{-1} Q_1 \omega_E \mod \omega_0(x)$$
d'où
$$P_1 - Q_1 p^{-1} a_p + Q_0 \equiv 0 \mod \omega_0(x)$$
$$P_0 \equiv p^{-1} Q_1 \mod \omega_0(x) \equiv p^{-1} \xi_1 R_0 \mod \omega_0(x)$$
$$\equiv R_0 \mod \omega_0(x) \ .$$

On en déduit finalement que

$$P_{n+1} - a_p P_n + \xi_n P_{n-1} \equiv 0 \bmod \omega_n(x) \text{ pour } n \geq 1$$
$$P_1 \equiv P_0 a_p - Q_0 \bmod \omega_0(x) .$$

Lorsque $F(0) \in \varphi \operatorname{Fil}^0 D_p(E)$, $Q(0)$ est nul et on a donc $P_1 \equiv P_0 a_p \bmod \omega_0(x)$. Lorsque $(1 - p^{-1}\varphi^{-1})(1-\varphi)^{-1} F(0) \in \varphi \operatorname{Fil}^0 D_p(E)$, La relation sur $F(0)$ implique que $(p-1)P_0 = (a_p - 2)Q_0$ et donc

$$(a_p - 2)P_1 \equiv ((a_p - 2)a_p - (p-1)) P_0 \bmod \omega_0(x) .$$

$\square$

*Si $P$ est un polynôme non nul de $\mathbb{Q}_p[x]$, par le théorème de Weierstrass, il peut s'écrire de manière unique sous la forme $P = u p^\mu (x^\lambda + pR)$ où $u$ est une unité de $\mathbb{Z}_p[[x]]$ et où $R \in \mathbb{Z}_p[x]$ est de degré $< \lambda$. On appelle $\lambda = \lambda(P)$ et $\mu = \mu(P)$ les $\lambda$ et $\mu$ invariants de $P$. De manière équivalente, on peut écrire $P = v p^\mu (x^\lambda + pQ + x^{\lambda+1} Q^1)$ avec $v \in \mathbb{Z}_p^*$ et $Q$ et $Q^1 \in \mathbb{Z}_p[x]$. Par exemple, $\lambda(\omega_n) = p^{n-1}(p-1)$, $\mu(\omega_n(x)) = 0$, $\lambda(\xi_n) = p^{n-2}(p-1)$ et $\mu(\xi_n) = 0$. Si $F$ est comme dans le lemme 4.1.1 et $P_n$ la suite de polynômes de degré $< p^n$ tels que*

$$F \equiv P_n \varphi^{n+1} \omega_E - \xi_n P_{n-1} \varphi^{n+2} \omega_E \bmod \omega_n(x) D_p(E)$$

*on note $\mu_+(F)$ (resp. $\mu_-(F)$) le minimum des $\mu(P_n)$ pour $n$ pair (resp. impair).*

*Si $n$ est un entier et $K_n = \mathbb{Q}(\mu_{p^n})$, on écrit*

$$\tilde{\omega}_{K_n} = \tilde{\omega}_{K/K_{n_0}} \otimes \tilde{\omega}_{K_{n_0}}$$

*(notations du paragraphe 2.2).*

4.1.2. **Proposition.** *Soit $F \in \mathcal{H} \otimes D_p(E)$, non nul, d'ordre $\leq 0$ et vérifiant*

$$F(\zeta_n - 1) \in \varphi^{n+1} \operatorname{Fil}^0 D_p(E)$$

*pour toute racine de l'unité $\zeta_n$ d'ordre $p^n$ non triviale et soit $\mu_\pm = \mu_\pm(F)$. Alors,*

*1) $F$ a un nombre fini de zéros de la forme $\zeta - 1$ avec $\zeta$ racine de l'unité d'ordre une puissance de $p$ : il existe un entier $n_0$ tel que $F(\zeta_n - 1) \neq 0$ pour $n \geq n_0$.*

*2) Si $a_p = 0$ ou si $\mu_+ = \mu_-$, il existe des rationnels $\lambda_+$, $\lambda_-$ et $\nu$ tels que pour $n \geq n_0$*

$$\operatorname{ord}_p \left( \frac{\prod_{\zeta \in \mu_{p^n} - \mu_{p^{n_0}}} F(\zeta - 1)}{\tilde{\omega}_{K/K_{n_0}}} \right)$$
$$= \frac{p^{2\lfloor \frac{n}{2} \rfloor + 1} - p}{p+1} \mu_+ + \frac{p^{2\lfloor \frac{n+1}{2} \rfloor} - 1}{p+1} \mu_- + \frac{p}{p^2 - 1}(p^n - 1) + \lambda_+ \lfloor \frac{n}{2} \rfloor + \lambda_- \lfloor \frac{n+1}{2} \rfloor - \nu .$$

*Une formule analogue existe dans le cas où $a_p \neq 0$ et $\mu_+ \neq \mu_-$. Nous ne l'avons pas écrite car aucun cas de ce type n'a été rencontré numériquement.*

**Remarque.** Si

$$A_n = \frac{p^{2\lfloor \frac{n}{2} \rfloor + 1} - p}{p+1} \mu_+ + \frac{p^{2\lfloor \frac{n+1}{2} \rfloor} - 1}{p+1} \mu_- + \frac{p}{p^2 - 1}(p^n - 1) + \lambda_+ \lfloor \frac{n}{2} \rfloor + \lambda_- \lfloor \frac{n+1}{2} \rfloor$$

et si $F(\zeta - 1)$ ne s'annule jamais, on a

$$\operatorname{ord}_p \left( \frac{\prod_{\zeta \in \mu_{p^n}} F(\zeta - 1)}{\tilde{\omega}_{K/K_{n_0}}} \right)$$
$$= A_n - A_{n_0} + \operatorname{ord}_p \left( \frac{\prod_{\zeta \in \mu_{p^{n_0}}} F(\zeta - 1)}{\tilde{\omega}_{K_{n_0}}} \right) .$$

On peut aussi écrire $A_n$ sous la forme

$$A_n = \begin{cases} \left(\frac{p\mu_+ + \mu_-}{p+1} + \frac{p}{p^2-1}\right)(p^n - 1) + \frac{\lambda_+ + \lambda_-}{2}n & \text{si } n \text{ est pair} \\ \left(\frac{\mu_+ + p\mu_-}{p+1} + \frac{p}{p^2-1}\right)(p^n - 1) + (\lambda_+ + \lambda_-)\lfloor\frac{n}{2}\rfloor + \lambda_- & \text{si } n \text{ est impair.} \end{cases}$$

Nous verrons dans la démonstration que

$$\begin{cases} \lambda_+ \equiv -\frac{1}{p+1} \bmod \mathbb{Z} \\ \lambda_- \equiv -\frac{p}{p+1} \bmod \mathbb{Z} \end{cases}.$$

Ainsi, $\lambda_+ + \lambda_-$ est un entier, $\lambda_+ - \lambda_- \equiv \frac{p-1}{p+1} \bmod \mathbb{Z}$ et $A_n$ est bien un entier.

Un cas particulier important est celui où $\mu_+ = \mu_- = 0$, on a alors

$$A_n = \frac{p}{p^2-1}(p^n - 1) + \lambda_+ \lfloor\frac{n}{2}\rfloor + \lambda_- \lfloor\frac{n+1}{2}\rfloor.$$

Nous donnerons la définition de $\lambda_+$ et $\lambda_-$ au cours de la démonstration

4.1.1. *Démonstration.* Soit $P_n$ la suite de polynômes de degré $< p^n$ tels que

$$F \equiv P_n \varphi^{n+1} \omega_E - \xi_n P_{n-1} \varphi^{n+2} \omega_E \bmod \omega_n(x) D_p(E)$$

On a les relations

(4.1.1) $\qquad P_n - a_p P_{n-1} + \xi_{n-1} P_{n-2} \equiv 0 \bmod \omega_{n-1}(x) D_p(E)$

pour $n \geq 2$. Lorsque $a_p = 0$ (resp. $a_p \neq 0$), dès que $P_n \neq 0$ pour un entier $n$, les polynômes $P_m$ sont non nuls pour $m \geq n$ (resp. pour $m$ assez grand) de même parité. Soit $n$ tel que $P_n \neq 0$. Écrivons $P_n = u_n p^{\mu_n}(x^{\lambda_n} + Q_n)$ avec $Q_n \in (p, x^{\lambda_n})\mathbb{Z}_p[x]$ et $u_n$ une unité de $\mathbb{Z}_p$. Comme le degré de $P_n$ est $< p^n$, $\lambda_n < p^n$. Dans la suite, $\epsilon$ désigne indifféremment un signe $\pm$ ou une classe modulo 2, c'est-à-dire une parité. Quitte à multiplier tous les $P_n$ par un entier indépendant de $n$, on peut supposer que $\mu_n \geq 0$ (c'est la condition que $F$ est d'ordre $\leq 0$ qui assure que cela est possible).

4.1.3. **Lemme.** *La suite des $\mu_n$ pour $n$ de parité $\epsilon$ est stationnaire. Si $a_p = 0$ ou si $\mu_+ = \mu_-$, il existe un entier $n_\epsilon \equiv \epsilon \bmod 2$ tel que pour $n \geq n_\epsilon$, $n \equiv \epsilon \bmod 2$,*

$$\lambda_n = \lambda_{n_\epsilon} + \frac{p^n - p^{n_\epsilon}}{p+1} = \lambda_\epsilon + \frac{p^n}{p+1};$$

*de plus, pour tout entier $n$, $\lambda_n \geq M_n \geq p^{n-2}(p-1)$ avec*

$$M_n = \begin{cases} \frac{p^n - 1}{p+1} & \text{si } n \text{ est pair} \\ \frac{p^n - p}{p+1} & \text{si } n \text{ est impair} \end{cases}.$$

*Si $a_p \neq 0$ et $\mu_+ \neq \mu_-$, pour $n$ assez grand, on a*

$$\lambda_{n+1} = \lambda_n = \lambda_\epsilon + \frac{p^n}{p+1}$$

*pour $n$ de parité $\epsilon$ avec $\mu_\epsilon < \mu_{-\epsilon}$ et $\lambda_n < p^{n-1}(p-1)$, $\lambda_{n+1} < p^n(p-1)$.*

On note $\mu_\epsilon$ la limite des $\mu_n$. La limite des $\mu_n$ est aussi sa borne inférieure. Lorsque $a_p = 0$ ou $\mu_+ = \mu_-$, la suite des $\mu_n$ (pour $n$ de parité fixe $\epsilon$) est décroissante, puis stationnaire et on peut prendre pour $n_\epsilon$ le plus petit entier tel que $\mu_n$ est minimal pour $n$ de parité $\epsilon$. Pour $n > n_\epsilon$ de parité $\epsilon$, $\lambda_n < p^{n-1}$ et $\lambda_n - \frac{p^n}{p+1}$ est indépendant de $n$. On note $\lambda_\epsilon$ la limite des $\lambda_n - \frac{p^n}{p+1}$ pour $n$ de parité $\epsilon$. On a

$$\mu_\pm = \mu_n \quad , \quad \lambda_\pm = \lambda_n - \frac{p^n}{p+1}$$

pour $n \geq n_\epsilon$ de parité $\epsilon$. Une variante entière de $\lambda_\pm$ est de poser

$$\tilde{\lambda}_\pm = \lim_{\substack{k \to \infty \\ (-1)^k = \pm}} \lambda_k - M_k = \lambda_n - M_n$$

pour $n$ tel que $\mu_n$ est minimal. On a donc $\lambda_+ = \tilde{\lambda}_+ - \frac{1}{p+1}$ et $\lambda_- = \tilde{\lambda}_- - \frac{p}{p+1}$.

*Démonstration.* Supposons d'abord $a_p = 0$. Fixons $\epsilon$, les entiers $n$ considérés sont supposés $\equiv \epsilon \mod 2$. La relation (4.1.1) s'écrit

$$u_n p^{\mu_{n+2}}(x^{\lambda_{n+2}} + pQ_{n+2} + x^{\lambda_{n+2}+1}Q^1_{n+2})$$
$$+ u_{n-2}(x^{p^n(p-1)} + pZ_{n+1})p^{\mu_n}(x^{\lambda_n} + pQ_n + x^{\lambda_n+1}Q^1_n))$$
$$= u'_n(x^{p^{n+1}} + p\delta_{n+1})p^{\nu_{n+1}}(x^{g_{n+1}} + pR_{n+1} + x^{g_{n+1}+1}R^1_{n+1})$$

avec $u'_n$ une unité et $R_{n+1}$ et $R^1_{n+1}$ des polynômes à coefficients dans $\mathbb{Z}_p$, $\xi_{n+1} = x^{p^n(p-1)} + pZ_{n+1}$. Si $\mu_n > \mu_{n+2}$, on a $\lambda_{n+2} \geq p^{n+1} \geq \frac{p^{n+2}-1}{p+1}$. Supposons $\mu_n \leq \mu_{n+2}$. Comme $\lambda_n < p^n$,

$$\lambda_n + p^n(p-1) < p^{n+1} .$$

L'identité (4.1.1) implique alors que $\lambda_{n+2} = \lambda_n + p^n(p-1) < p^{n+1}$ et que $\mu_{n+2} = \mu_n$. Donc la suite des $\mu_n$ est décroissante, puis stationnaire. Soit $n_\epsilon$ le plus petit des entiers $n$ tels que $\mu_n$ soit minimal pour $n$ de parité $\epsilon$. On a alors $\mu_{n_\epsilon+2} = \mu_{n_\epsilon}$ et pour tout entier $n \geq n_\epsilon$ de parité $\epsilon$, $\mu_n = \mu_{n_\epsilon}$. Ainsi, pour tout entier $n \geq n_\epsilon$, $n \equiv \epsilon \mod 2$, on a

$$\lambda_{n+2} = \lambda_n + p^n(p-1) .$$

On en déduit que pour $n \geq n_\epsilon$

$$\lambda_n = \lambda_{n_\epsilon} + \frac{p^n - p^{n_\epsilon}}{p+1} = \lambda_\epsilon + \frac{p^n}{p+1} .$$

L'inégalité $\lambda_n \geq M_n$ se déduit directement de la relation $P_{n+2} \equiv -\xi_{n+1}P_n \mod \omega_n$ et d'un calcul de degré que nous reverrons dans le paragraphe 4.1.2.

Supposons maintenant que $a_p$ est non nul (on a donc $p = 3$ et $a_3 = \pm 3$). Soit $\mu_\epsilon$ le minimum des $\mu_n$ pour $n$ de parité $\epsilon$. Quitte à diviser les polynômes $P_n$ par une puissance de $p$, on peut supposer que pour une des parités, disons $\epsilon$, $\mu_{-\epsilon} \geq \mu_\epsilon = 0$. Soit $n$ un entier de parité $\epsilon$ tel que $\mu_n = 0$. Alors la relation (4.1.1) appliquée à $n+2$ implique que $\mu_{n+2} = \mu_\epsilon = 0$ puisque $\mu(a_p P_{n+1}) \geq 1$. Donc $\mu_m = 0$ pour tout entier $m \geq n$ de parité $\epsilon$ et on a encore $\lambda_{m+2} = \lambda_m + \lambda(\xi_{m+1}) = \lambda_m + p^m(p-1)$ et $\lambda_n = \lambda_\epsilon + \frac{p^n}{p+1}$. Pour la parité $-\epsilon$, soit $n$ de parité $-\epsilon$ tel que $\mu_n = \mu_{-\epsilon}$ et tel que $\mu_{n+1} = 0$. Trois cas peuvent se produire.

Cas 1 : $\mu_{-\epsilon} > \text{ord}_p(a_p) = 1$. La relation (4.1.1) appliquée à $n+2$ implique que $\mu_{n+2} = \text{ord}_p a_p < \mu_{-\epsilon}$, ce qui contredit la définition de $\mu_{-\epsilon}$.

Cas 2 : $\mu_{-\epsilon} < \text{ord}_p(a_p) = 1$. Alors, $\mu(P_{n+2}) = \mu(P_n) = \mu_{-\epsilon}$ et on a $\mu_m = \mu_{-\epsilon}$ pour tout entier $m \geq n$ de même parité. Dans ce cas, comme dans le cas $a_p = 0$, $\lambda_{m+2} = \lambda_m + p^m(p-1)$.

Cas 3 : $\mu_{-\epsilon} = \text{ord}_p(a_p) = 1$. Si $\mu_{n+2} > \mu_{-\epsilon}$, on a

$$\lambda_n = \lambda_{n+1} - p^n(p-1) = \lambda_\epsilon + \frac{p^{n+1}}{p+1} - p^n(p-1)$$
$$= \lambda_\epsilon - \frac{p^n(p^2 - p - 1)}{p+1}$$

Ce qui n'est pas possible pour $n$ assez grand. Donc, pour $n$ assez grand de parité $-\epsilon$, $\mu_n = \mu_{-\epsilon} = \text{ord}_p(a_p)$. Comme $\lambda(\xi_{n+1}) = p^n(p-1) > \frac{p^{n+1}}{p+1}$ et que $\lambda_{n+1} = \lambda_\epsilon + \frac{p^{n+1}}{p+1}$, on déduit de l'égalité 4.1.1 pour $n+2$ que $\lambda_{n+2} = \lambda_{n+1} = \lambda_\epsilon + \frac{p^{n+1}}{p+1}$. En particulier, pour $n$ assez grand, $\lambda_{n+2} < p^{n+1}(p-1)$. □

**4.1.4. Lemme.** *Soit $P$ un polynôme vérifiant $\lambda(P) < p^{n-1}(p-1))$. Alors, $P(\zeta_n-1)$ est non nul pour $\zeta_n$ racine de l'unité d'ordre $p^n$ et on a alors*

$$\operatorname{ord}_p P(\zeta_n - 1) = \mu(P) + \frac{\lambda(P)}{p^{n-1}(p-1)} \ .$$

*Démonstration.* On a $P = u_n p^{\mu(P)}(x^{\lambda(P)} + pQ)$ avec $Q \in \mathbb{Z}_p[x]$, $u$ une unité de $\mathbb{Z}_p[[x]]$. Si $\zeta_n$ est une racine de l'unité d'ordre $p^n$,

$$\operatorname{ord}_p((\zeta_n - 1)^{\lambda(P)}) = \frac{\lambda(P)}{p^{n-1}(p-1)} < \frac{1}{p-1} < 1 \ .$$

Donc $P(\zeta_n - 1)$ ne s'annule pas et

$$\operatorname{ord}_p(P(\zeta_n - 1) = \mu(P) + \frac{\lambda(P)}{p^{n-1}(p-1)} \ .$$

□

4.1.2.  *Démontrons maintenant la proposition 4.1.2. La première assertion se déduit du lemme 4.1.4 (elle est aussi démontrée dans [19]). Prenons maintenant $a_p = 0$ ou $\mu_+ = \mu_-$ (le calcul est inspiré du calcul de Kurihara).*

Soit $n_0 = \sup(n_+, n_-)$. Il s'agit de calculer (avec $p^{n_\zeta}$ l'ordre de $\zeta$)

$$\sum_{\zeta \in \mu_{p^n} - \mu_{p^{n_0}}} \operatorname{ord}_p P_{n_\zeta}(\zeta - 1) = \sum_{j=n_0+1}^{n} \sum_{\zeta \in \mu_{p^j} - \mu_{p^{j-1}}} \operatorname{ord}_p P_j(\zeta - 1)$$

$$= \sum_{j=n_0+1}^{n} p^{j-1}(p-1) \left( \mu_{\epsilon(j)} + \frac{p}{p^2 - 1} + \frac{\lambda_{\epsilon(j)}}{p^{j-1}(p-1)} \right)$$

$$= A_n - A_{n_0}$$

avec

$$A_n = \sum_{j=1}^{n} \left( p^{j-1}(p-1) \left( \mu_{\epsilon(j)} + \frac{p}{p^2 - 1} \right) + \lambda_{\epsilon(j)} \right)$$

$$= \sum_{j=1}^{n} p^{j-1}(p-1)\mu_{\epsilon(j)} + \sum_{j=1}^{n} \frac{p^j}{p+1} + \sum_{j=1}^{n} \lambda_{\epsilon(j)} \ .$$

On a

$$\sum_{\substack{j=1 \\ j \ pair}}^{n} p^{j-1}(p-1)\tilde{\mu}_+ = \frac{p^{2\lfloor \frac{n}{2} \rfloor + 1} - p}{p+1} \mu_+$$

$$\sum_{\substack{j=1 \\ j \ impair}}^{n} p^{j-1}(p-1)\mu_- = \frac{p^{2\lfloor \frac{n+1}{2} \rfloor} - 1}{p+1} \mu_-$$

$$\sum_{j=1}^{n} \frac{p^j}{p+1} = \frac{p}{p^2 - 1}(p^n - 1)$$

$$\sum_{j=1}^{n} \lambda_{\epsilon(j)} = \lambda_+ \sum_{\substack{j \ pair}}^{n} 1 + \lambda_- \sum_{\substack{j \ impair}}^{n} 1$$

$$= \lambda_+ \lfloor \frac{n}{2} \rfloor + \lambda_- \lfloor \frac{n+1}{2} \rfloor \ .$$

*Donc*

$$A_n = \frac{p^{2\lfloor \frac{n}{2} \rfloor+1} - p}{p+1}\mu_+ + \frac{p^{2\lfloor \frac{n+1}{2} \rfloor} - 1}{p+1}\mu_- + \frac{p}{p^2-1}(p^n - 1)$$
$$+ \lambda_+ \lfloor \frac{n}{2} \rfloor + \lambda_- \lfloor \frac{n+1}{2} \rfloor \ .$$

*Donc*

$$A_n - A_{n_0} = \frac{p^{2\lfloor \frac{n}{2} \rfloor+1} - p^{2\lfloor \frac{n_0}{2} \rfloor+1}}{p+1}\mu_+ + \frac{p^{2\lfloor \frac{n+1}{2} \rfloor} - p^{2\lfloor \frac{n_0+1}{2} \rfloor}}{p+1}\mu_- + \frac{p}{p^2-1}(p^n - p^{n_0})$$
$$+ \lambda_+(\lfloor \frac{n}{2} \rfloor - \lfloor \frac{n_0}{2} \rfloor) + \lambda_-\left(\lfloor \frac{n+1}{2} \rfloor - \lfloor \frac{n_0+1}{2} \rfloor\right) \ .$$

*Ce qui termine la démonstration de la proposition 4.1.2 dans le cas $a_p = 0$ ou $\mu_+ = \mu_-$. Dans le cas contraire, lorsque $\mu_\epsilon < \mu_{-\epsilon}$, le calcul serait le même à condition de prendre 0 pour $\lambda_{-\epsilon}$ et $2\lambda_\epsilon$ pour $\lambda_\epsilon$.*

**Exemple.** Supposons que $F$ vérifie la condition supplémentaire

$$(1 - p^{-1}\varphi^{-1})(1 - \varphi)^{-1}F(0) \in \mathrm{Fil}^0 D_p(E) \ ,$$

que les polynômes $P_n$ sont à coefficients dans $\mathbb{Z}_p$ et que $P_0$ est une unité $p$-adique. En particulier $\mu_0 = 0$. La relation supplémentaire reliant $P_0$ et $P_1$ implique que $\mu_1 = 0$. Donc, $\mu_+ = \mu_- = 0$, $\lambda_1 = 0$, $\lambda_0 = 0$, on peut prendre $n_+ = 0$ et $n_- = 1$, donc $\tilde{\lambda}_+ = \tilde{\lambda}_- = 0$ et $\lambda_+ = -1/(p+1)$, $\lambda_- = -p/(p+1)$. D'où, $A_1 = 0$ et en multipliant par $F(0)$ (sous les hypothèses faite, $F(0)/(1-p^{-1}\varphi^{-1})^{-1}(1-\varphi)\omega_E$ est une unité), on obtient la formule

$$\mathrm{ord}_p\left(\frac{\prod_{\zeta \in \mu_{p^n}} F(\zeta-1)}{\tilde{\omega}_K}\right) = \frac{p(p^n-1)}{p^2-1} - \frac{1}{p+1}\lfloor \frac{n}{2} \rfloor - \frac{p}{p+1}\lfloor \frac{n+1}{2} \rfloor$$
$$= \lfloor \frac{p}{p^2-1}p^n - \frac{n}{2} \rfloor \ .$$

**Remarques.** 1) Nous avons vu au cours de la démonstration que si $\mu_{n_0} = \mu_+$, pour $n > n_0$, la fonction $F$ ne peut plus s'annuler en $\zeta - 1$ pour $\zeta$ racine de l'unité d'ordre $p^n$. Si de plus $\lambda_{n_0} < p^{n_0-1}(p-1)$, $F$ ne s'annule pas non plus en une racine de l'unité d'ordre $p^{n_0}$.

2) Supposons $a_p = 0$. Soit

$$M_n = \begin{cases} \frac{p^n - p}{p+1} & \text{si } n \text{ est impair} \\ \frac{p^n - 1}{p+1} & \text{si } n \text{ est pair} \end{cases} \ .$$

Pour tout entier $n$ de parité $\epsilon$ tel que $\mu_n = \mu_\epsilon$, on a $\lambda_n \geq M_n$.

3) $\lambda_+$ (resp. $\lambda_-$) se calcule à partir de $\lambda_{2n}$ (resp. $\lambda_{2n+1}$) dès que $\mu_{2n} = \mu_+$ (resp. $\mu_{2n+1} = \mu_-$).

**Définition.** Si $F$ est comme dans le lemme 4.1.1, on définit le $\lambda$-invariant de $F$ comme le couple $(\lambda_+, \lambda_-)$ et son $\mu$-invariant comme le couple $(\mu_+, \mu_-)$.

*Si $F_1$ et $F_2$ sont deux éléments de $\mathcal{H} \otimes D_p(E)$ d'ordre $\leq 0$ comme dans le lemme 4.1.1 tel que $F_2 = hF_1$ avec $\lambda_+(F_1) = \lambda_+(F_2)$ et $\mu_+(F_1) = \mu_+(F_2)$ et $h \in \mathbb{Z}_p[[x]]$. Alors $(F_1) = (F_2)$. On remarque en effet que si $h \in \mathbb{Z}_p[[x]]$, $\lambda(hF) = \lambda(F) + (\lambda(h), \lambda(h))$ et $\mu(hF) = \mu(F) + (\mu(h), \mu(h))$.*

## 5. Croissance du groupe de Shafarevich-Tate

5.1. *On peut maintenant mettre ensemble les paragraphes 2.2 et 4 pour obtenir le résultat suivant qui redémontre et généralise le résultat de Kurihara.*

5.1.1. **Théorème.** *Soit $E$ une courbe elliptique définie sur $\mathbb{Q}$ et $p$ un nombre premier tel que $E$ soit supersingulière en $p$.*

*1) Le rang de $E(K_\infty)$ est fini et le corang de $\mathrm{III}(E/K_n)$ est borné.*

*2) Supposons $E(K_n)$ et $\mathrm{III}(E/K_n)$ finis pour tout $n$. Alors il existe des entiers $\mu'_+$, $\mu'_-$, $\lambda'_+$, $\lambda'_-$ et un rationnel $\nu$ tels que*

$$\mathrm{ord}_p \sharp \mathrm{III}(E/K_n) = \frac{p^{2\lfloor \frac{n}{2} \rfloor + 1}}{p+1} \mu'_+ + \frac{p^{2\lfloor \frac{n+1}{2} \rfloor}}{p+1} \mu'_- + \frac{p}{p^2-1} p^n + \lambda'_+ \lfloor \frac{n}{2} \rfloor + \lambda'_- \lfloor \frac{n+1}{2} \rfloor + \nu \ .$$

*3) (Kurihara) Si $\frac{L(E,1)}{\Omega_E}$ est une unité en $p$, c'est-à-dire si $\mathbf{1}(L_p(E,1))$ n'est pas divisible par $p$ dans le réseau $M_E$ de $D_p(E)$, $E(K_n)$ et $\mathrm{III}(E/K_n)(p)$ sont finis pour tout $n$ et on a*

$$\mathrm{ord}_p \sharp \mathrm{III}(E/K_n)(p) = \lfloor \frac{p}{p^2-1} p^n - \frac{n}{2} \rfloor \ .$$

*4) Soit un entier $n_0$ tel que le rang de $E(K_n)$ et le corang de $\mathrm{III}(E/K_n)(p)$ soient stationnaires pour $n \geq n_0$ (qui existe d'après 1) et soit $s$ le rang de $\check{S}_p(E/\mathbb{Q}_\infty) = \check{S}_p(E/\mathbb{Q}_{n_0})$. Alors, il existe des entiers $\mu'_+$, $\mu'_-$, $\lambda'_+$, $\lambda'_-$ et un rationnel $\nu$ tels que pour $n \geq n_0$*

$\mathrm{ord}_p \sharp (\mathrm{III}(E/K_n)/div)$

$$= \frac{p^{2\lfloor \frac{n}{2} \rfloor + 1}}{p+1} \mu'_+ + \frac{p^{2\lfloor \frac{n+1}{2} \rfloor}}{p+1} \mu'_- + \frac{p}{p^2-1} p^n + (\lambda'_+ - s) \lfloor \frac{n}{2} \rfloor + (\lambda'_- - s) \lfloor \frac{n+1}{2} \rfloor + \nu \ .$$

*Ce théorème se déduit de la proposition 4.1.2 appliquée à la fonction $\hat{I}_{arith,(0)} \in \mathcal{H} \otimes D_p(E)$ (même définition que pour $L_p(E/\mathbb{Q})$, §1.2.1). Les propriétés de $I_{arith}$ impliquent en effet que $\hat{I}_{arith,(0)}$ vérifie les conditions du lemme 4.1.1. On a alors pour $\delta$ caractère d'ordre fini de $G_\infty$ :*

$$I_{arith}(\langle \chi \rangle^s \delta) = \hat{I}_{arith,(0)}(u^s \zeta - 1)$$

*avec $u = \langle \chi(\gamma) \rangle$ et $\zeta = \delta(\gamma)$. De plus, $I_{arith}(\langle \chi_{\mathbb{Q}_n} \rangle^s) = \prod_{\delta \in \hat{\Delta}_{\mathbb{Q}_n}} I_{arith}(\langle \chi \rangle^s \delta)$ où $\chi_K$ est le caractère cyclotomique attachée à $\mathbb{Q}_n$, c'est-à-dire la restriction de $\chi$ à $G_{\mathbb{Q}_n}$. Le rang de $\check{S}_p(E/\mathbb{Q}_n)$ est stationnaire pour $n$ assez grand. Comme $\check{S}_p(E/\mathbb{Q}_n)$ est facteur direct dans $H^1(\mathbb{Q}_n, T_p(E))$ et que le conoyau de l'application de restriction $H^1(\mathbb{Q}_n, T_p(E)) \to H^1(\mathbb{Q}_m, T_p(E))^{\mathrm{Gal}(\mathbb{Q}_m/\mathbb{Q}_n)}$ est $H^2(\mathbb{Q}_m/\mathbb{Q}_n, T_p(E)^{\mathrm{Gal}(\overline{\mathbb{Q}}/\mathbb{Q}_m)})$ et est donc nul, la limite inductive des $\check{S}_p(E/\mathbb{Q}_n)$ est un $\mathbb{Z}_p$-module de type fini. D'autre part, on a $\langle x, y \rangle_{\mathbb{Q}_m} = [\mathbb{Q}_m : \mathbb{Q}_n] \langle x, y \rangle_{\mathbb{Q}_n}$ pour $x$ et $y$ dans $\check{S}_p(E/\mathbb{Q}_n)$ et $m \geq n$. Donc, pour $n_0$ tel que $\check{S}_p(E/\mathbb{Q}_\infty) = \check{S}_p(E/\mathbb{Q}_{n_0})$*

$$\prod_{\delta \in \hat{\Delta}_n - \hat{\Delta}_{n_0}} I_{arith}(\langle \chi \rangle^s \delta) = \frac{\mathrm{Tam}(E/\mathbb{Q}_n) \sharp \mathrm{III}(T_p(E)/\mathbb{Q}_n)}{\mathrm{Tam}(E/\mathbb{Q}_0) \sharp \mathrm{III}(T_p(E)/\mathbb{Q}_{n_0})} p^{(n-n_0)s} \tilde{\omega}_{E, \mathbb{Q}_n/\mathbb{Q}_{n_0}}$$

*En utilisant la proposition 4.1.2, le fait que les $\mathrm{Tam}(E/\mathbb{Q}_n)$ deviennent stationnaires et que $\mathrm{III}(T_p(E)/\mathbb{Q}_n)$ est le quotient de $\mathrm{III}(E/\mathbb{Q}_n)$ par sa partie divisible, on en déduit le théorème.*

## 6. EXEMPLES NUMÉRIQUES

6.1. **Quelques généralités.** *Nous allons voir quelques critères numériques qui permettent de conclure à la conjecture de Birch et Swinnerton-Dyer en $p$ (toujours supersingulier) dans des exemples qui ne semblaient pas connus auparavant et à la conjecture principale. Nous avons fait des calculs systématiques sur les courbes elliptiques de conducteur $N_E < 150$ et leurs twists par des discriminants inférieurs à $500$ en valeur absolue et pour des nombres premiers $3, 5, 7$ lorsqu'ils sont supersinguliers systématiquement et occasionnellement pour d'autres nombres premiers* ($p = 11$, $p = 19$). Précisons que lorsqu'on parle de l'ordre du groupe de Shafarevich-Tate, il s'agit ici réellement du nombre d'éléments et non de l'ordre conjectural (appelé quelquefois ordre analytique).

Soit $F = L_p(E)$ ou un générateur de $\mathbb{I}_{arith}$. C'est un élément de $\mathcal{H}(G_\infty) \otimes D_p(E)$. On écrit $F = \sum_{\omega^i \in \hat{\Delta}} \hat{F}_i(\gamma - 1) = \sum_i \hat{F}_i(\gamma - 1)$ avec $F_i = e_{\omega^i}(F) \in \mathcal{H}$. Si $\omega^i \delta$ est de conducteur $p^{n+1}$ avec $\delta$ d'ordre $p^n$, on a alors $\omega^i \delta(F) = \delta(F_i) = \hat{F}_i(\zeta_n - 1)$ avec $\zeta_n = \delta(\gamma)$ une racine de l'unité d'ordre $p^n$.

Prenons $i = 0$. Alors, $\hat{F}_0$ est un élément de $\mathcal{H} \otimes D_p(E)$ vérifiant les propriétés du lemme 4.1.1. On note $P_n$ les polynômes associés à $L_p(E)$ et au caractère trivial de $\mathrm{Gal}(\mathbb{Q}(\mu_p)/\mathbb{Q})$ comme dans le lemme 4.1.1 et $Q_n$ les polynômes associés à un générateur de $\mathbb{I}_{arith}$. Les polynômes $P_n$ et $Q_n$ sont à coefficients dans $\mathbb{Z}_p$. On note $\lambda_n = \lambda(P_n)$ et $\mu_n = \mu(P_n)$ les invariants des polynômes $P_n$ et $\lambda'_n = \lambda(Q_n)$ et $\mu'_n = \mu(Q_n)$ les invariants de $Q_n$. Ces polynômes vérifient les relations pour $n \geq 2$

(6.1.1)
$$P_n - a_p P_{n-1} + \xi_{n-1} P_{n-2} \equiv 0 \bmod \omega_{n-1}(x)$$
$$Q_n - a_p Q_{n-1} + \xi_{n-1} Q_{n-2} \equiv 0 \bmod \omega_{n-1}(x)$$

et

(6.1.2)
$$(a_p - 2)P_1 \equiv ((a_p - 2)a_p - (p - 1)) P_0 \bmod \omega_0(x)$$
$$(a_p - 2)Q_1 \equiv ((a_p - 2)a_p - (p - 1)) Q_0 \bmod \omega_0(x) .$$

Par le théorème de Kato, on sait qu'il existe un élément $g \in \mathbb{Z}_p[[x]]$ tel que

(6.1.3)
$$P_n \equiv g Q_n \bmod \omega_n(x) .$$

Montrer la conjecture principale revient alors à montrer que $g(0)$ est une unité, c'est-à-dire que $\mu(g) = \lambda(g) = 0$.

Pour tous les exemples calculés dans la suite, l'image de $\rho_p$ est maximale. Nous ne le répéterons pas. Les donnes numriques obtenues sont mises dans un tableau avec le nom de la courbe, les coefficients $[a_1, a_2, a_3, a_4, a_6]$, le discriminant par lequel on twiste, les $\mu_i$, $\lambda_i$ intéressants, puis ventuellement les invariants $\tilde{\lambda}_+$ et $\tilde{\lambda}_-$.

Commençons par quelques résultats numériques concernant les $\mu$-invariants. Nous avons mis dans le tableau I les valeurs numériques pour la suite des $\mu_n$ trouvées par le calcul ($\mu_n = \infty$ signifie que $P_n = 0$, une case blanche signifie que $\mu_n$ est automatiquement nul). Dans le tableau II, seules des courbes vérifiant $a_3 = \pm 3$ interviennent pour $p = 3$. Remarquons entre autres que $\mu_2$ peut être supérieur à $\mu_1$ (par exemple $\mu_1 = 0$, $\mu_2 = 1$). Dans tous les exemples calculés, $\mu_3$ et $\mu_4$ sont nuls. Il semble raisonnable de faire la conjecture :

6.1.1. **Conjecture.** *Les invariants $\mu_+$ et $\mu_-$ de $L_p(E/\mathbb{Q})$ sont nuls.*

Remarquons que nous n'avons pas dans le cas supersingulier le problème d'"invariance par isogénie" : $E$ n'a pas de sous-groupes rationnels sur $\mathbb{Q}$ d'ordre $p$.

| $p$ | $\mu_0$ | $\mu_1$ | $\mu_2$ | $\mu_3$ | $\mu_4$ |
|---|---|---|---|---|---|
| 3 | 0 | 0 | | | |
| 3 | 1 | 0 | 0 | | |
| 3 | 1 | 0 | 1 | | 0 |
| 3 | 1 | 1 | 0 | 0 | |
| 3 | 1 | 1 | 1 | 0 | 0 |
| 3 | 2 | 0 | 0 | | |
| 3 | 2 | 0 | 1 | | 0 |
| 3 | 2 | 1 | 0 | 0 | |
| 3 | 2 | 1 | 1 | 0 | 0 |
| 3 | 2 | 2 | 0 | 0 | |
| 3 | 2 | 2 | 1 | 0 | 0 |
| 3 | 3 | 0 | 0 | | |
| 3 | 3 | 1 | 0 | 0 | |
| 3 | 3 | 2 | 0 | 0 | |
| 3 | 4 | 0 | 0 | | |
| 3 | $\infty$ | 0 | 0 | | |
| 3 | $\infty$ | 0 | 1 | | 0 |
| 3 | $\infty$ | 0 | $\infty$ | | 0 |
| 3 | $\infty$ | 1 | 0 | 0 | |
| 3 | $\infty$ | 1 | 1 | 0 | 0 |
| 3 | $\infty$ | 2 | 0 | 0 | |
| 3 | $\infty$ | 2 | 1 | 0 | 0 |
| 3 | $\infty$ | 3 | 0 | 0 | |

| $p$ | $\mu_0$ | $\mu_1$ | $\mu_2$ | $\mu_3$ | $\mu_4$ |
|---|---|---|---|---|---|
| 3 | $\infty$ | 3 | 1 | 0 | 0 |
| 3 | $\infty$ | 4 | 0 | 0 | |
| 3 | $\infty$ | $\infty$ | 0 | 0 | |
| 3 | $\infty$ | $\infty$ | 1 | 0 | 0 |
| 5 | 0 | 0 | | | |
| 5 | 1 | 0 | 0 | | |
| 5 | 1 | 1 | 0 | 0 | |
| 5 | 2 | 0 | 0 | | |
| 5 | 2 | 1 | 0 | 0 | |
| 5 | $\infty$ | 0 | | | |
| 5 | $\infty$ | 0 | 0 | | |
| 5 | $\infty$ | 1 | 0 | 0 | |
| 5 | $\infty$ | $\infty$ | 0 | 0 | |
| 7 | 0 | 0 | | | |
| 7 | 2 | 0 | 0 | | |
| 7 | $\infty$ | 0 | | | |
| 7 | $\infty$ | 1 | 0 | 0 | |
| 7 | $\infty$ | 1 | 0 | 0 | |
| 7 | $\infty$ | $\infty$ | 0 | 0 | |
| 11 | 0 | 0 | 0 | | |
| 11 | $\infty$ | 0 | 0 | | |
| 17 | 2 | 0 | 0 | | |

Table I

| $p$ | $\mu_0$ | $\mu_1$ | $\mu_2$ | $\mu_3$ | $\mu_4$ |
|---|---|---|---|---|---|
| 3 | 0 | 0 | 0 | 0 | 0 |
| 3 | 2 | 0 | 0 | 0 | 0 |
| 3 | 2 | 0 | 1 | 0 | 0 |
| 3 | 2 | 1 | 0 | 0 | 0 |
| 3 | 2 | 2 | 0 | 0 | 0 |
| 3 | 4 | 0 | 0 | 0 | 0 |
| 3 | $\infty$ | 0 | 1 | 0 | 0 |
| 3 | $\infty$ | 1 | 0 | 0 | 0 |
| 3 | $\infty$ | 2 | 0 | 0 | 0 |
| 3 | $\infty$ | 3 | 0 | 0 | 0 |

Table II ($a_3 \neq 0$)

*Donnons maintenant une conséquence de l'équation fonctionnelle. Si $\epsilon_{arith}$ est le signe de l'équation fonctionnelle de $\mathbb{I}_{arith}(E)$, les $\lambda'_n$ sont tous pairs pour $\epsilon_{arith} = 1$ et tous impairs pour $\epsilon_{arith} = -1$. Malheureusement, on ne sait pas montrer que $\epsilon_{arith} = \epsilon_{anal}$.*

*Rappelons enfin que lorsque $a_p = 0$ ou $\mu_+ = \mu_-$, les $\lambda'_k = \lambda(Q_k)$ ainsi que les $\lambda_k = \lambda(P_k)$ vérifient les inégalités $\lambda_k \geq M_k$, $\lambda'_k \geq M_k$ avec*

$$M_k = \begin{cases} \frac{p^k - p}{p+1} & \text{si } k \text{ est impair} \\ \frac{p^k - 1}{p+1} & \text{si } k \text{ est pair} \end{cases}.$$

*et que l'on pose $\tilde{\lambda}_+ = \lambda_k - M_k$, $\tilde{\lambda}'_+ = \lambda'_k - M_k$ (resp. $\tilde{\lambda}_- = \lambda_k - M_k$, $\tilde{\lambda}'_- = \lambda'_k - M_k$) pour $k$ pair (resp. impair) tel que $\mu_k = 0$. Donnons quelques valeurs des $M_k$.*

| $p$ | $M_0$ | $M_1$ | $M_2$ | $M_3$ | $M_4$ |
|---|---|---|---|---|---|
| 3 | 0 | 0 | 2 | 6 | 20 |
| 5 | 0 | 0 | 4 | 20 | 104 |
| 7 | 0 | 0 | 6 | 42 | 300 |
| 11 | 0 | 0 | 10 | 110 | 1220 |

*Et pour fixer les idées donnons une petite table des degrés des polynômes cyclotomiques de degré $p^n$ :*

| $p$ | 1 | 2 | 3 | 4 |
|---|---|---|---|---|
| 3 | 2 | 6 | 18 | 54 |
| 5 | 4 | 20 | 100 | 500 |
| 7 | 6 | 42 | 294 | 2058 |
| 11 | 10 | 110 | 1210 | 13310 |

*Nous supposerons dans la suite que si $a_p \neq 0$, $\mu_+ = \mu_-$ pour ne pas alourdir les énoncés (le cas contraire n'ayant jamais été rencontré numériquement).*

6.2. **Courbes de rang 0.**

*Exemple de Kurihara. Commençons par le cas étudié par Kurihara : $P_0$ est une unité (l'hypothèse $a_p = 0$ faite dans [11] n'est pas nécessaire). Grâce à la relation (6.1.2), il en est de même de $P_1$. Les $\mu_n$ sont tous nuls. On a $\lambda_0 = \lambda_1 = 0$ et $\tilde{\lambda}_+ = \tilde{\lambda}_- = 0$. Par le lemme 4.1.3, on a alors $\lambda_2 = p-1$, $\lambda_3 = p(p-1)$, $\lambda_4 = (p-1)(p^2+1)$, $\lambda_5 = p(p-1)(p^2+1)$ et $\lambda_+ = \frac{-1}{p+1}$, $\lambda_- = \frac{-p}{p+1}$. Il est rassurant de vérifier que c'est ce que donne le calcul numérique : par exemple, pour $E = X_0(17)$ et $p = 3$, on a $\lambda_2 = 2$, $\lambda_3 = 6$, $\lambda_4 = 20$, $\lambda_5 = 60$. La conjecture principale $e_0 \mathbb{I}_{arith}(E) = e_0 L_p(E/\mathbb{Q})\mathbb{Z}_p[[\Gamma]]$ est vraie car la congruence $P_0 \equiv g Q_0 \bmod x$ signifie que $g(0)$ est une unité; $E(\mathbb{Q}_\infty)$ est fini et les $\text{III}(\mathbb{Q}_n)(p)$ sont tous finis et de cardinal*

$$\text{ord}_p \sharp \text{III}(E/\mathbb{Q}_n) = \frac{p}{p^2-1}p^n - \frac{1}{p+1}\lfloor \frac{n}{2} \rfloor - \frac{p}{p+1}\lfloor \frac{n+1}{2} \rfloor = \lfloor \frac{p}{p^2-1}p^n - \frac{n}{2} \rfloor .$$

*Par exemple pour $p = 3$, on obtient*

$$\lfloor \frac{3^{n+1}}{8} - \frac{n}{2} \rfloor .$$

*Exemple de rang 0 et de groupe de Shafarevich-Tate non trivial. Prenons $E = X_0(17)^{(373)}$, $p = 3$. La valeur de $L(E/\mathbb{Q})/\Omega_E$ est $-36$. Les calculs donnent $\mu_0 = 2$, $\mu_1 = 0$, $\mu_2 = 0$, $\lambda_1 = 2$, $\lambda_2 = 6$. Donc*

$$\tilde{\lambda}_+ = 4 , \quad \tilde{\lambda}_- = 2$$

*On vérifie en même temps que $P_1$ et $P_2$ sont premiers respectivement à $\xi_1 = x^2 + 3x + 3$ et à $\xi_2$. Le rang de $E(\mathbb{Q})$ est 0. En utilisant la remarque 4.1.2, on obtient que $E(\mathbb{Q}_\infty)$ est de torsion et que $\text{III}(E/\mathbb{Q}_n)(3)$ est fini pour tout entier $n$. Passons maintenant à la conjecture principale. On a $\text{Tam}(E) = 16$ et $\sharp E(\mathbb{Q}) = 2$. Ainsi, $(1 - \varphi)^{-1}(1 - p^{-1}\varphi^{-1})\mathbf{1}(I_{arith}(E/\mathbb{Q}))/\omega_E = \sharp \text{III}(E/\mathbb{Q})(p) u$ divise 9 (avec $u$ unité de $\mathbb{Z}_3$). Comme on sait que $\text{III}(E/\mathbb{Q})$ est un carré, le cardinal de $\text{III}(E/\mathbb{Q})(3)$ est 9 ou 1. Nous allons montrer que le deuxième cas n'est pas possible. On conclut alors en remarquant que $P_0$ et $Q_0$ ont même valuation 3-adique et donc que $g(0)$ est une unité 3-adique puisque $P_0 = g(0)Q_0$, ce qui démontrera la conjecture principale $e_0 \mathbb{I}_{arith}(E/\mathbb{Q}) = e_0 L_p(E/\mathbb{Q})\mathbb{Z}_p[[\Gamma]]$.*

*Supposons que $\text{III}(E/\mathbb{Q})(3)$ est trivial. Dans ce cas, $Q_0$ est une unité. Comme en 6.2, les $\lambda'_n$ et $\mu'_n$ valent $\mu'_0 = \mu'_1 = 0$, $\lambda'_0 = 0$, $\lambda'_1 = 0$ et donc $\lambda'_2 = 2$, $\lambda'_3 = 6$, $\lambda'_4 = 20$, $\lambda'_5 = 60$ et $\tilde{\lambda}'_+ = \tilde{\lambda}'_- = 0$. Or les congruences $P_n \equiv gQ_n \bmod (1+x)^{3^n} - 1$ et le fait que les $\mu_n$ sont nuls et les $P_n$ non nuls (de degré $< 3^n$) pour $n \geq 1$ impliquent que $\lambda_n = \lambda(g) + \lambda'_n < 3^n$ pour tout entier $n \geq 1$. On en déduit en particulier que*

$$\lambda_2 - \lambda_1 = \lambda'_2 - \lambda'_1 = 2$$
$$\lambda_3 - \lambda_2 = \lambda'_3 - \lambda'_2 = 4 .$$

*ce qui est contradictoire avec les résultats numériques. Donc, $\text{III}(E/\mathbb{Q})(3)$ n'est pas trivial et d'ordre 9. On a donc ici montré par des calculs simples de symboles modulaires que*

$$\sharp \text{III}(X_0(17)^{(373)}/\mathbb{Q})(3) = 9 .$$

*De plus, $\lambda_- = 2 - 3/4 = 5/4$ et $\lambda_+ = 6 - 9/4 = 15/4$. Donc,*

$$\mathrm{ord}_3 \sharp \mathrm{III}(X_0(17)^{(373)}\mathbb{Q}_n)(3) = 2 + (A_n - A_0)$$

*avec*

$$A_n = \frac{3(3^n - 1)}{8} + \frac{15}{4}\lfloor \frac{n}{2} \rfloor + \frac{5}{4}\lfloor \frac{n+1}{2} \rfloor$$

*D'où,*

$$\mathrm{ord}_3 \sharp \mathrm{III}(X_0(17)^{(373)}/\mathbb{Q}_n)(3) = 2 + \frac{3(3^n - 1)}{8} + \frac{15}{4}\lfloor \frac{n}{2} \rfloor + \frac{5}{4}\lfloor \frac{n+1}{2} \rfloor \ .$$

*Remarquons que $L(E, 1)/\Omega_E$ est égal à 36, pour p ne divisant pas 6, les cas traités par Kurihara pour p supersingulier et par Mazur dans le cas où p est ordinaire montrent que la conjecture de Birch-Swinnerton-Dyer (à condition de vérifier que $\rho_p$ est surjective) est vraie à des puissances de 2, de 17 et de 373 près.*

*On peut généraliser cet exemple de la manière suivante.*

**6.2.1. Proposition.** *Supposons $\rho_p$ surjective, $\mu_0 \leq 2$ (en particulier, $L(E,1) \neq 0$) et $\tilde{\lambda}_+ \neq \tilde{\lambda}_-$. Alors, la conjecture principale est vraie et $\mathrm{III}(E/\mathbb{Q})(p)$ a l'ordre prédit par la conjecture de Birch et Swinnerton-Dyer. Si de plus, $\lambda_k - M_k < p^{k-1}(p-1)$ pour les entiers k tels que $\mu_{k-2} \neq 0$, alors $\mathrm{III}(E/\mathbb{Q}_n)(p)$ est fini pour tout entier n et $E(\mathbb{Q}_\infty) = 0$.*

*Ainsi, si $\mathrm{Tam}(E)/(\sharp E(\mathbb{Q}))^2$ est premier à p et si $\mu_0 = 2$, $\mathrm{III}(E/\mathbb{Q})(p)$ est d'ordre $p^2$. Si $\mu_0 = 1$, $\mathrm{III}(E/\mathbb{Q})(p)$ est trivial (par un théorème de Cassels, $\sharp \mathrm{III}(E/\mathbb{Q})(p)$ est un carré).*

*Démonstration.* On sait déjà que $E(\mathbb{Q})$ est fini ainsi que $\mathrm{III}(E/\mathbb{Q})(p)$. D'autre part $\sharp E(\mathbb{Q})_{tors}$ est d'ordre premier à p. Si la conjecture principale n'est pas vraie, comme $\sharp \mathrm{III}(E/\mathbb{Q})(p)$ est un carré, $Q_0$ est nécessairement une unité et $\tilde{\lambda}'_+ = \tilde{\lambda}'_- = 0$. On en déduit que $\tilde{\lambda}_+ = \tilde{\lambda}_-$, ce qui a été supposé faux. Donc $g$ est une unité. $\square$

**Exemple.** Soit $E = X_0(17)^{(-167)}$. Les données numériques sont

| $E_0$ | équation | $D$ | $p$ | $\mu_0$ | $\mu_1, \lambda_1$ | $\mu_2, \lambda_2$ | $\mu_3, \lambda_3$ | $\mu_4, \lambda_4$ |
|---|---|---|---|---|---|---|---|---|
| 17A | [1,-1,1,-1,-14] | -167 | 3 | 2 | 1,2 | 1,6 | 0,10 | 0,28 |

On a $\mu_0 = 2$, $\mu_1 = 1$, $\mu_2 = 1$, $\mu_n = 0$ pour $n \geq 3$, $\lambda_1 = 2$, $\lambda_2 = 6$, $\lambda_3 = 10$, $\lambda_4 = 28$. Donc $\tilde{\lambda}_+ = 28 - 20 = 8$ et $\tilde{\lambda}_- = 10 - 6 = 4$. Donc la conjecture principale est vraie. Le nombre de Tamagawa $\mathrm{Tam}(E)$ est égal à 4. Donc, $\mathrm{III}(E/\mathbb{Q})(3)$ est d'ordre 9. D'autre part, $P_1 = -12x^2 - 36x - 36 = -12\xi_1$. On en déduit que le rang de $\check{S}_3(E/\mathbb{Q}(\mu_9)^\Delta)$ est supérieur ou égal à $p - 1 = 2$ grâce au théorème 2.2.2 et au fait que la conjecture principale vient d'être vérifiée. Ce rang peut être 2 ou 4. Il faudrait vérifier numériquement que $L'_p(E) \not\equiv 0 \mod \xi_1$. Le polynôme

$$P_2 = 24x^8 + 216x^7 + 852x^6 + 1944x^5 + 2808x^4 + 2628x^3 + 1548x^2 + 540x + 108$$
$$= 12(x^2 + 3x + 3)(2x^6 + 12x^5 + 29x^4 + 39x^3 + 30x^2 + 12x + 3)$$

est premier avec $\xi_2$. Pour $n = 3$ et 4, la comparaison des degrés implique que $P_n$ est premier avec $\xi_n$ et donc pour $n \geq 5$ puisque $\mu_3 = \mu_4 = 0$. Ainsi, le rang de la limite inductive des $\check{S}_3(E/\mathbb{Q}(\mu_{3^n})^\Delta)$ est 2, $E(\mathbb{Q}_\infty)$ est de rang 2 ou $\mathrm{III}(E/\mathbb{Q}(\mu_9)^\Delta)(3)$ est infini. On a $\lambda_1 - M_1 = 2$, $\lambda_2 - M_2 = 4$, $\lambda_3 - M_3 = 4$, $\lambda_4 - M_4 = 8$, $\lambda_n - M_n = 4$ si n est un entier impair $\geq 3$ et $\lambda_n - M_n = 8$ si n est un entier pair $\geq 4$. On a alors

$$A_n - A_2 = \frac{3(3^n - 1)}{8} + \frac{31}{4}\lfloor \frac{n}{2} \rfloor + \frac{13}{4}\lfloor \frac{n+1}{2} \rfloor - 14$$

$$= \begin{cases} \frac{3^{n+1} - 115}{8} + \frac{11n}{2} \\ \frac{3^{n+1} - 133}{8} + \frac{11n}{2} \end{cases}$$

Ainsi,
$$\operatorname{ord}_p \sharp \mathbf{III}(E/\mathbb{Q}_n)(3) = \frac{3(3^n-1)}{8} + \frac{23}{4}\lfloor\frac{n}{2}\rfloor + \frac{5}{4}\lfloor\frac{n+1}{2}\rfloor + \nu$$

**Exemple.** Soit $E = X_0(17)^{(-187)}$ :

| $E_0$ | équation | $D$ | $p$ | $\mu_0$ | $\mu_1,\lambda_1$ | $\mu_2,\lambda_2$ | $\mu_3,\lambda_3$ | $\tilde{\lambda}_-$ | $\tilde{\lambda}_+$ |
|---|---|---|---|---|---|---|---|---|---|
| 17A | [1,-1,1,-1,-14] | -187 | 3 | 2 | 1,2 | 0,4 | 0,10 | 4 | 2 |

Le nombre de Tamagawa de $E$ est 8. On a $\mu_0 = 2$, $\mu_1 = 1$, $\mu_2 = 0$, $\mu_n = 0$ pour $n \geq 2$, $\lambda_1 = 2$, $\lambda_2 = 4$, $\lambda_3 = 10$. Donc, $\tilde{\lambda}_- = 4$ et $\tilde{\lambda}_+ = 2$. Comme $\tilde{\lambda}_+ \neq \tilde{\lambda}_+$, la conjecture principale est vraie et $\mathbf{III}(E/\mathbb{Q})(3)$ est d'ordre 9. Le polynôme $P_1$ est égal à $-24\xi_1$. Quand à $P_2$, il vaut

$$P_2 = -32x^8 - 284x^7 - 1140x^6 - 2680x^5 - 4012x^4 - 3948x^3 - 2568x^2 - 1080x - 216$$
$$= -4(x^2 + 3x + 3)(8x^6 + 47x^5 + 120x^4 + 169x^3 + 136x^2 + 72x + 18)$$

et est premier à $\xi_2$. On en déduit que les polynômes $P_n$ sont tous premiers à $\xi_n$ pour $n \geq 2$, que $E(\mathbb{Q}_\infty)$ est de rang $\geq 2$ ou $\mathbf{III}(E/\mathbb{Q}(\mu_9)^\Delta)(3)$ est infini, que le rang de la limite inductive des $\check{S}_3(E/\mathbb{Q}(\mu_{3^n})^\Delta)$ est 2 ou 4. Pour conclure pour le rang, il faut vérifier que $\hat{L}'_p(E)) \not\equiv 0 \bmod \xi_1$. Ce que je n'ai pas encore fait. Enfin,

$$A_n - A_1 = \frac{3(3^n-1)}{8} + \frac{7}{4}\lfloor\frac{n}{2}\rfloor + \frac{13}{4}\lfloor\frac{n+1}{2}\rfloor - 4$$

et

$$\operatorname{ord}_p \sharp \mathbf{III}(E/\mathbb{Q}_n)(3) = \frac{3(3^n-1)}{8} - \frac{1}{4}\lfloor\frac{n}{2}\rfloor - \frac{3}{4}\lfloor\frac{n+1}{2}\rfloor + \nu$$

**Exemple.** Soit $E = 40A^{(-379)}$ avec $40A$ d'équation $y^2 = x^3 - 7x - 6$ :

| $E_0$ | équation | $D$ | $p$ | $\mu_0$ | $\mu_1,\lambda_1$ | $\mu_2,\lambda_2$ | $\mu_3,\lambda_3$ | $\mu_4,\lambda_4$ |
|---|---|---|---|---|---|---|---|---|
| 40A | [0,0,0,-7,-6] | -379 | 3 | 2 | 2,0 | 1,4 | 0,10 | 0,28 |

Le nombre $\operatorname{Tam}(E)$ est égal à 32. On prend toujours $p = 3$. On a $\mu_0 = 2$, $\mu_1 = 1$, $\mu_2 = 1$, $\mu_n = 0$ pour $n \geq 3$ et $\lambda_1 = 0$, $\lambda_2 = 4$, $\lambda_3 = 10$, $\lambda_4 = 28$ et donc $\tilde{\lambda}_+ = 8$, $\tilde{\lambda}_- = 4$. La conjecture principale est vraie et $\mathbf{III}(E/\mathbb{Q})(3)$ est d'ordre 9. Pour des questions de degré, les polynômes $P_n$ sont tous premiers à $\xi_n$. Donc $E(\mathbb{Q}_\infty)$ est nul et $\mathbf{III}(E/\mathbb{Q}_n)(p)$ est fini pour tout entier $n$. On a

$$\operatorname{ord}_3 \mathbf{III}(E/\mathbb{Q}_n)(3) = A_n - A_2 + \operatorname{ord}_3 P_0 + \operatorname{ord}_3 \prod_{\zeta \in \mu_3 - \{1\}} P_1(\zeta - 1)$$
$$+ \operatorname{ord}_3 \prod_{\zeta \in \mu_9 - \mu_3} P_2(\zeta - 1) \ .$$

Ici, $P_1 = -144(x+1)$ et

$$P_2 = -48x^8 - 432x^7 - 1728x^6 - 4080x^5 - 6288x^4 - 6480x^3 - 4320x^2 - 1728x - 432$$
$$= -48\xi_1(x^6 + 6x^5 + 15x^4 + 22x^3 + 20x^2 + 9x + 3) \ .$$

Donc, $\operatorname{ord}_3 \prod_{\zeta \in \mu_3 - \{1\}} P_1(\zeta - 1) = 4$ et par le calcul

$$\operatorname{ord}_3 \prod_{\zeta \in \mu_9 - \mu_3} P_2(\zeta - 1) = 10 \ .$$

D'autre part,

$$A_n = \frac{3(3^n-1)}{8} + (\lambda_+ + \lambda_-)\lfloor\frac{n}{2}\rfloor + (\lambda_- \quad \text{si } n \text{ est impair})$$
$$= \frac{3(3^n-1)}{8} + \lambda_+\lfloor\frac{n}{2}\rfloor + \lambda_-\lfloor\frac{n+1}{2}\rfloor$$

avec $\lambda_+ = 28 - 3^4/4 = 31/4$ et $\lambda_- = 10 - 3^3/4 = 13/4$; d'où

$$A_n = \frac{3(3^n-1)}{8} + 11\lfloor\frac{n}{2}\rfloor(+\frac{13}{4} \quad \text{si } n \text{ est impair})$$
$$A_2 = 14 .$$

D'où pour $n \geq 3$

$$\text{ord}_3 \sharp \mathbf{III}(E/\mathbb{Q}_n)(3) = \frac{3(3^n-1)}{8} + \frac{31}{4}\lfloor\frac{n}{2}\rfloor + \frac{13}{4}\lfloor\frac{n+1}{2}\rfloor + 2 .$$

**Exemple.** $E = 62A^{(-296)}$ :

| $E_0$ | équation | $D$ | $p$ | $\mu_0$ | $\mu_1, \lambda_1$ | $\mu_2, \lambda_2$ | $\mu_3, \lambda_3$ | $\mu_4, \lambda_4$ | $\tilde{\lambda}_-$ | $\tilde{\lambda}_+$ |
|---|---|---|---|---|---|---|---|---|---|---|
| 62A | [1,-1,1,-1,1] | -296 | 3 | 2 | 1,2 | 1,4 | 0,12 | 0,28 | 6 | 8 |

Le rang de $E(\mathbb{Q}_\infty)$ est nul, les $\mathbf{III}(E/\mathbb{Q}_n)(3)$ sont finis. On a $\text{Tam}(E) = 8$, $\lambda_+ = 31/4$, $\lambda_- = 21/4$, $\tilde{\lambda}_+ = 8$, $\tilde{\lambda}_- = 6$. Enfin,

$$A_n - A_2 = \frac{3(3^n-1)}{8} + \frac{31}{4}\lfloor\frac{n}{2}\rfloor + \frac{21}{4}\lfloor\frac{n+1}{2}\rfloor - 16$$
$$= \begin{cases} \frac{3^{n+1}-131}{8} + \frac{13n}{2} \\ \frac{3^{n+1}-141}{8} + \frac{13n}{2} \end{cases} .$$

On a

$$\begin{cases} P_1 = 12(x^2 - 3x - 3) \\ P_2 = 3\xi_1(14x^6 + 75x^5 + 152x^4 + 166x^3 + 113x^2 + 36x + 12) \end{cases} .$$

On trouve alors que la contribution de $P_0$ (resp. $P_1$, $P_2$) à la valuation 3-adique de $\sharp\mathbf{III}(E/\mathbb{Q}_n)(3)$ pour $n \geq 3$ est 2 (resp. 4, 10).

On obtient donc que

$$\text{ord}_3 \sharp \mathbf{III}(E/\mathbb{Q}_n)(3) = \frac{3(3^n-1)}{8} + \frac{31}{4}\lfloor\frac{n}{2}\rfloor + \frac{21}{4}\lfloor\frac{n+1}{2}\rfloor .$$

**Exemple.** $E = 73A^{(-151)}$ :

| $E_0$ | équation | $D$ | $p$ | $\mu_0$ | $\mu_1, \lambda_1$ | $\mu_2, \lambda_2$ | $\mu_3, \lambda_3$ | $\mu_4, \lambda_4$ | $\tilde{\lambda}_-$ | $\tilde{\lambda}_+$ |
|---|---|---|---|---|---|---|---|---|---|---|
| 73A | [1, -1, 0, 4, -3] | -151 | 3 | 2 | 2,0 | 1,8 | 0,10 | 0,28 | 4 | 8 |

On a $\tilde{\lambda}_+ = 8$, $\tilde{\lambda}_- = 4$, $\lambda_+ = 31/4$, $\lambda_- = 13/4$. La conjecture principale est vraie. Le rang de $\check{S}_3(E/\mathbb{Q}_\infty)$ est supérieur à 6 car $P_2$ n'est pas premier à $\xi_2$ : on a $P_1 = -36$ et $P_2 = -12\xi_1\xi_2$. Donc, $\check{S}_3(\mathbb{Q}(\mu_9)^\Delta)$ est infini et de rang $\geq 6$. Si le rang de $\check{S}_3(\mathbb{Q}(\mu_9)^\Delta)$ était plus grand, et donc plus grand que 12, $\lambda_+$ et $\lambda_-$ seraient respectivement supérieurs à 12 et 6, ce qui est faux. Donc $\check{S}_3(\mathbb{Q}_\infty)$ est de rang 6.

**Exemple.** $E = 142C^{(397)}$ :

| $E_0$ | équation | $D$ | $p$ | $\mu_0$ | $\mu_1, \lambda_1$ | $\mu_2, \lambda_2$ | $\mu_3, \lambda_3$ | $\mu_4, \lambda_4$ |
|---|---|---|---|---|---|---|---|---|
| 142C | [1, -1, 0, -1, -3] | 397 | 3 | 1 | 1,0 | 1,2 | 0,10 | 0,30 |

Ici le nombre de Tamagawa est 12 et n'est donc pas premier à 3. Les invariants de $E$ pour $p = 3$ sont $\mu_0 = 1$, $\mu_1 = 1$, $\mu_2 = 1$, $\mu_n = 0$ pour $n \geq 3$, $\lambda_1 = 0$, $\lambda_2 = 2$, $\lambda_3 = 10$, $\lambda_4 = 30$. Donc, $\tilde{\lambda}_+ = 10$, $\tilde{\lambda}_- = 4$ et la conjecture principale est vraie. La composante 3-primaire de $\mathbf{III}(E/\mathbb{Q})$ est triviale. Le groupe $E(\mathbb{Q}_\infty)$ est nul. Les nombres de Tamagawa restent de valuation 3-adique 1. On obtient pour le membre de gauche

$$\begin{cases} 1 + 10 + \frac{3^{n+1}-131}{8} + \frac{13n}{2} \\ 1 + 10 + \frac{3^{n+1}-157}{8} + \frac{13n}{2} \end{cases} = \begin{cases} \frac{3^{n+1}-43}{8} + \frac{13n}{2} \\ \frac{3^{n+1}-69}{8} + \frac{13n}{2} \end{cases} .$$

L'ordre du groupe de Shafarevich-Tate est de valuation 3-adique

$$\begin{cases} \frac{3^{n+1}-51}{8} + \frac{13n}{2} \\ \frac{3^{n+1}-77}{8} + \frac{13n}{2} \end{cases}$$

*Donnons quelques exemples avec $p = 5, 7$. La croissance du groupe de Shafarevich-Tate s'en déduit comme précédemment.*

**Exemple.** $E = 106C^{(-312)}$ :

| $E_0$ | équation | $D$ | $p$ | $\mu_0$ | $\mu_1, \lambda_1$ | $\mu_2, \lambda_2$ | $\tilde{\lambda}_-$ | $\tilde{\lambda}_+$ |
|---|---|---|---|---|---|---|---|---|
| 106C | [1,0,0,-283,-2351] | -312 | 5 | 2 | 0,4 | 0,6 | 4 | 2 |

La conjecture principale 5-adique est vraie, $E(\mathbb{Q}_\infty) = 0$, $\mathrm{III}(E/\mathbb{Q})(5)$ est d'ordre 25, les $\mathrm{III}(E/\mathbb{Q}_n)(5)$ sont finis, $(\lambda_+, \lambda_-) = (11/6, 19/6)$.

**Exemple.** $E = 84A^{(-443)}$ :

| $E_0$ | équation | $D$ | $p$ | $\mu_0$ | $\mu_1, \lambda_1$ | $\mu_2, \lambda_2$ | $\tilde{\lambda}_-$ | $\tilde{\lambda}_+$ |
|---|---|---|---|---|---|---|---|---|
| 84A | [0,1,0,7,0] | -443 | 5 | 2 | 0,2 | 0,8 | 2 | 4 |

La conjecture principale 5-adique est vraie, $E(\mathbb{Q}_\infty) = 0$, $\mathrm{III}(E/\mathbb{Q})(5)$ est d'ordre 25, les $\mathrm{III}(E/\mathbb{Q}_n)(5)$ sont finis et $(\lambda_+, \lambda_-) = (23/6, 7/6)$.

**Exemple.** $E = 106B^{(-408)}$ :

| $E_0$ | équation | $D$ | $p$ | $\mu_0$ | $\mu_1, \lambda_1$ | $\mu_2, \lambda_2$ | $\tilde{\lambda}_-$ | $\tilde{\lambda}_+$ |
|---|---|---|---|---|---|---|---|---|
| 106B | [1,1,0,-7,5] | -408 | 7 | 2 | 0,2 | 0,10 | 2 | 4 |

La conjecture principale 7-adique est vraie, $E(\mathbb{Q}_\infty) = 0$, $\mathrm{III}(E/\mathbb{Q})(7)$ est d'ordre 49, les $\mathrm{III}(E/\mathbb{Q}_n)(7)$ sont finis et $(\lambda_+, \lambda_-) = (31/8, 9/8)$.

*La proposition suivante permet d'obtenir quelques cas où $\mu_0 \geq 3$.*

**6.2.2. Proposition.** *Supposons $a_p = 0$ et $L(E, 1) \neq 0$. Si $\tilde{\lambda}_+$ ou $\tilde{\lambda}_-$ est égal à 2 et si $\tilde{\lambda}_+ \neq \tilde{\lambda}_-$ ou si $\mathrm{ord}_p \mathrm{Tam}(E) \geq 1$, alors la conjecture principale est vraie et la composante p-primaire du groupe de Shafarevich-Tate a l'ordre prédit.*

*Démonstration.* Pour $k$ de parité convenable $\epsilon$ et assez grand, on a

$$0 \leq \lambda(g) \leq \lambda_k - M_k = \tilde{\lambda}_\epsilon = 2 \ .$$

D'autre part, $\mu(g) = 0$. Les équations fonctionnelles pour $L_p(E)$ et $\mathbb{I}_{arith}$ ont toutes deux le signe + puisqu'elles sont non nulles en **1**, donc $\lambda(g)$ est nécessairement pair. Si $\lambda(g) = 0$, $g$ est une unité et la conjecture principale est vraie. Sinon, $\lambda(g) = 2$. Dans ce cas $\tilde{\lambda}'_\epsilon = 0$ et $\lambda'_k = M_k$. Le fait qu'il soit égal à la valeur limite implique que $\lambda'_j$ est minimal pour tout entier $j$, ce qui n'est possible que si $\mu'_0 = 0$. Si $\mathrm{Tam}(E)$ est divisible par $p$, cela n'est pas possible. Sinon, on a nécessairement $\tilde{\lambda}_- = \tilde{\lambda}_+$. □

**Exemple.** $E = 142C^{(53)}$ :

| $E_0$ | équation | $D$ | $p$ | $\mu_0$ | $\mu_1, \lambda_1$ | $\mu_2, \lambda_2$ | $\tilde{\lambda}_-$ | $\tilde{\lambda}_+$ |
|---|---|---|---|---|---|---|---|---|
| 142C | [1,-1,0,-1,-3] | 53 | 3 | 3 | 0,2 | 0,6 | 2 | 4 |

La conjecture principale est vraie. Le nombre de Tamagawa est 12, ce qui est prévisible puisque $\mu_0$ est impair. Donc, l'ordre de $\mathrm{III}(E/\mathbb{Q})(3)$ est 9. Là encore $E(\mathbb{Q}_\infty) = 0$ car les polynômes $P_1$ et $P_2$ sont premiers respectivement à $\xi_1$ et à $\xi_2$. La variation du cardinal de $\mathrm{III}(E/\mathbb{Q}_n)(3)$ peut être obtenue facilement comme précédemment.

**Exemple.** $E = 142C^{(461)}$ :

| $E_0$ | équation | $D$ | $p$ | $\mu_0$ | $\mu_1, \lambda_1$ | $\mu_2, \lambda_2$ | $\mu_3, \lambda_3$ | $\tilde{\lambda}_-$ | $\tilde{\lambda}_+$ |
|---|---|---|---|---|---|---|---|---|---|
| 142C | [1,-1,0,-1,-3] | 461 | 3 | 3 | 2,2 | 0,4 | 0,12 | 6 | 2 |

La conjecture principale est vraie. Le nombre de Tamagawa est 12. Donc, l'ordre de $\text{III}(E/\mathbb{Q})(3)$ est 9. Ici, $P_1 = 18\xi_1$ et $P_2$ est premier avec $\xi_2$. Donc, le rang de $\check{S}_p(\mathbb{Q}(\mu_9)^\Delta)$ est $\geq 2$. S'il était égal à 4, on aurait alors $\hat{L}'_p(E/\mathbb{Q}) \equiv 0 \mod \xi_1$, ce qu'il faut vérifier numériquement. En tout cas $\check{S}_p(\mathbb{Q}_\infty) = \check{S}_p(\mathbb{Q}(\mu_9)^\Delta)$.

Les courbes suivantes vérifient la conjecture principale et l'ordre de $\text{III}(E/\mathbb{Q})(3)$ est égal à 9 :

| $E_0$ | équation | $D$ | $p$ | $\mu_0$ | $\mu_1, \lambda_1$ | $\mu_2, \lambda_2$ | $\mu_3, \lambda_3$ | $\tilde{\lambda}_-$ | $\tilde{\lambda}_+$ |
|---|---|---|---|---|---|---|---|---|---|
| 142C | [1,-1,0,-1,-3] | 461 | 3 | 3 | 2,2 | 0,4 | 0,12 | 6 | 2 |
| 142C | [1,-1,0,-1,-3] | -139 | 3 | 3 | 1,2 | 0,4 | 0,10 | 4 | 2 |
| 142C | [1,-1,0,-1,-3] | -467 | 3 | 3 | 2,2 | 0,4 | 0,10 | 4 | 2 |
| 142C | [1,-1,0,-1,-3] | 53 | 3 | 3 | 0,2 | 0,6 | | 2 | 4 |
| 52A | [0,0,0,1,-10] | -499 | 3 | 3 | 2,2 | 0,4 | 0,10 | 4 | 2 |
| 52A | [0,0,0,1,-10] | 469 | 3 | 3 | 1,2 | 0,4 | 0,10 | 4 | 2 |
| 142C | [1,-1,0,-1,-3] | -211 | 3 | 3 | 0,2 | 0,4 | | 2 | 2 |
| 52A | [0,0,0,1,-10] | -331 | 3 | 3 | 0,2 | 0,4 | | 2 | 2 |
| 124B | [0,0,0,-17,-27] | 109 | 3 | 3 | 0,2 | 0,4 | | 2 | 2 |

Mais pour beaucoup de courbes, les arguments précédents ne permettent pas de conclure : par exemple

| $E_0$ | équation | $D$ | $p$ | $\mu_0$ | $\mu_1, \lambda_1$ | $\mu_2, \lambda_2$ | $\mu_3, \lambda_3$ | $\tilde{\lambda}_-$ | $\tilde{\lambda}_+$ |
|---|---|---|---|---|---|---|---|---|---|
| 142C | [1,-1,0,-1,-3] | 485 | 3 | 3 | 1,2 | 0,6 | 0,12 | 6 | 4 |
| 142C | [1,-1,0,-1,-3] | 493 | 3 | 3 | 1,2 | 0,6 | 0,16 | 10 | 4 |
| 124B | [0,0,0,-17,-27] | 485 | 3 | 3 | 2,2 | 0,8 | 0,10 | 4 | 6 |
| 124B | [0,0,0,-17,-27] | 205 | 3 | 3 | 1,2 | 0,6 | 0,10 | 4 | 4 |

**6.3. Courbes de rang $\geq 1$.** *Rappelons que*

$$L'_{p,\omega_E}(E, \mathbf{1}) = [(1-\varphi)^{-1}(1-p^{-1}\varphi^{-1})\mathbf{1}(L'_p(E)), \omega_E]_{D_p(E)}$$
$$= [\mathbf{1}(L'_p(E)), \tilde{\omega}_E]_{D_p(E)} .$$

**6.3.1. Proposition.** *Supposons $\rho_p$ surjective, $L(E/\mathbb{Q}, 1) = 0$ et $L'_{p,\omega_E}(E, \mathbf{1}) \neq 0$, $\mu_+ = \mu_- = 0$. Si $\lambda_+ = 1$ ou $\lambda_- = 1$, la conjecture principale est vraie et le rang de $\check{S}_p(E)$ est égal à 1. Si $E(\mathbb{Q})$ est infini, il est de rang 1, l'ordre de $\text{III}(E/\mathbb{Q})(p)$ est celui prévu par la conjecture $p$-adique de Birch et Swinnerton-Dyer : ainsi si $P$ engendre un $\mathbb{Z}$-module de $E(\mathbb{Q})$ d'indice premier à $p$, $\text{III}(E/\mathbb{Q})(p)$ est fini et on a*

$$\text{ord}_p \text{III}(E/\mathbb{Q})(p) = \text{ord}_p \frac{L'_{p,\omega_E}(E, \mathbf{1})}{p} - \text{ord}_p \text{Tam}(E)(\frac{\log_{\omega_E} P}{p})^2$$

*Remarquons que la conjecture principale est ainsi prouvée sans calculer de points rationnels, mais uniquement avec des calculs de symboles modulaires (plus ou moins longs, j'en conviens !). Malheureusement, pour l'instant dans tous les calculs numériques, $\text{ord}_p \text{III}(E/\mathbb{Q})(p)$ est toujours trouvé nul.*

**Remarque.** On peut aussi obtenir des renseignements supplémentaires sur $E(\mathbb{Q}_\infty)$. Par exemple, si $\mu_1 = 0$, $\mu_2 = 0$ et $\lambda_+ < p(p-1)$ (pour assurer que $P_2$ est premier à $\xi_1$), $E(\mathbb{Q}_\infty) = E(\mathbb{Q})$ et $\text{III}(E/\mathbb{Q}_n)(p)$ est fini pour tout entier $n$. On obtient alors une formule pour la croissance du groupe de Shafarevich-Tate. De manière générale, si l'on veut connaître le rang de $E(\mathbb{Q}_\infty)$ et $\text{III}(E/\mathbb{Q}_n)(p)$, il faut de plus regarder si tous les polynômes $P_j$ avec les entiers $j$ tels que $\mu_j \neq 0$, ainsi que pour les deux premiers entiers $j$ et $j'$ respectivement pair et impair pour lesquelles $\mu_j = 0$, sont premiers à $\xi_j$.

*Démonstration.* La non-nullité de $L'_{p,\omega_E}(E, \mathbf{1})$ implique que le rang de $\check{S}_p(E/\mathbb{Q})$ est égal à 1 et donc que $I_{arith}$ est lui aussi nul sur le caractère trivial (3.4.1). Donc, $\tilde{\lambda}'_- \geq 1$ et $\tilde{\lambda}'_+ \geq 1$. Comme $\lambda_\pm = \lambda'_\pm + \lambda(g)$, $\lambda(g)$ est nul. Comme $\mu(g) = 0$, $g$ est une unité. □

**Exemple.** Prenons $E = 43A$ d'équation $y^2 + y = x^3 + x^2$ et $p = 7$. On a

| $E_0$ | équation | $D$ | $p$ | $\mu_0$ | $\mu_1, \lambda_1$ | $\mu_2, \lambda_2$ | $\tilde{\lambda}_-$ | $\tilde{\lambda}_+$ |
|---|---|---|---|---|---|---|---|---|
| 43A | [0,1,1,0,0] | 1 | 7 | $\infty$ | 0,1 | 0,9 | 1 | 3 |

Ainsi, $L(E, 1) = 0$ et ses invariants sont $\mu_1 = \mu_2 = 0$, $\lambda_1 = 1$, $\lambda_2 = 9$. On vérifie que $L'_{p,\omega_E}(E, \mathbf{1})$ est non nul : on a

$$L(E/\mathbb{Q}_7, 1)(1-\varphi)^{-1}(1 - 7^{-1}\varphi^{-1})\mathbf{1}(L'_7(E))$$
$$= (5{\times}7 + 6{\times}7^2 + O(7^3))\omega_E - 7(3{\times}7 + 4{\times}7^2 + O(7^3))\varphi\omega_E$$

Il avait déjà été vu que la conjecture principale 7-adique est vraie, que $\mathbf{III}(E/\mathbb{Q})(7)$ est trivial. Quand à $E(\mathbb{Q}_\infty)$, il est de rang 1 et $\mathbf{III}(E/\mathbb{Q}_n)(7)$ est fini pour tout entier $n$.

**Exemple.** Soit $E = 91A$ d'équation $y^2 + y = x^3 + x$ et $p = 3$. On a

| $E_0$ | équation | $D$ | $p$ | $\mu_0$ | $\mu_1, \lambda_1$ | $\mu_2$ | $\mu_4, \lambda_4$ | $\tilde{\lambda}_-$ | $\tilde{\lambda}_+$ |
|---|---|---|---|---|---|---|---|---|---|
| 91A | [0,0,1,1,0] | 1 | 3 | $\infty$ | 0, 1 | $\infty$ | 0, 27 | 1 | 7 |

On a $L(E, 1) = 0$, ses invariants sont $\mu_1 = 0$, $\mu_2 = \infty$, $\mu_4 = 0$ et $\lambda_1 = 1$, $\lambda_4 = 27$. On vérifie que $L'_{p,\omega_E}(E, \mathbf{1})$ est non nul : on a

$$L(E/\mathbb{Q}_3, 1)(1-\varphi)^{-1}(1 - 3^{-1}\varphi^{-1})\mathbf{1}(L'_3(E))$$
$$\sim (2{\times}3 + 2{\times}3^2 + 2{\times}3^4 + O(3^6))\omega_E$$
$$- 3(3 + 2{\times}3^2 + 3^3 + 2{\times}3^4 + O(3^6))\varphi\omega_E$$

On en déduit que $\check{S}_3(\mathbb{Q})$ est de rang 1, ce qui est bien sûr déjà connu ($E(\mathbb{Q})$ est engendré par le point $(0,0)$), que $\mathbf{III}(E/\mathbb{Q})(3)$ est trivial et que la conjecture principale est vraie. Comme $P_2$ est identiquement nul, $E(\mathbb{Q}(\mu_{27})^\Delta)$ ou au moins $\check{S}_3(E/\mathbb{Q}(\mu_{27})^\Delta)$ est de rang $\geq 1 + 6 = 7$. S'il était plus grand, on aurait $7 = \tilde{\lambda}_+ \geq 1 + 2{\times}7$ ce qui n'est pas possible. Un autre argument est le suivant. Si le rang était $\geq 15$, alors la série caractéristique de $H^2_{\infty, \{3\}}(\mathbb{Q}, T_3(E))^\Delta$ serait divisible par $\xi_2$ et $\lambda_-$ serait plus grand que $1+6$, ce qui n'est pas le cas. Donc, $\check{S}_3(\mathbb{Q}_\infty) = \check{S}_3(\mathbb{Q}(\mu_{27})^\Delta)$. Comme le régulateur de $E$ est

$$(2{\times}3 + 2{\times}3^2 + 2{\times}3^4 + O(3^6))\omega_E - 3(3 + 2{\times}3^2 + 3^3 + 2{\times}3^4 + O(3^6))\varphi\omega_E ,$$

que $\mathrm{Tam}(E/\mathbb{Q}) = 1$ et que $E(\mathbb{Q})$ n'a pas de torsion, tout va bien !

**Exemple.** Pour $37A$ d'équation $y^2 + y = x^3 - x$, les nombres premiers 17 et 19 sont supersinguliers. On trouve pour $L(E/\mathbb{Q}_p, 1)(1-\varphi)^{-1}(1-p^{-1}\varphi^{-1})\mathbf{1}(L'_p(E))$ en $\mathbf{1}$

$$(4{\times}17 + 11{\times}17^2 + O(17^3))\omega_E - 17(8{\times}17 + 12{\times}17^2 + O(17^3))\varphi\omega_E \text{ pour } p = 17$$
$$(13{\times}19 + 10{\times}19^2 + O(19^3))\omega_E - 19(18{\times}19 + 7{\times}19^2 + O(19^3))\varphi\omega_E \text{ pour } p = 19$$

Sans calculer de points, ces calculs impliquent que $\check{S}_p(E/\mathbb{Q})$ est de rang 1 pour ces deux nombres premiers. Bien sûr on sait que $E(\mathbb{Q})$ est infini et admet le point $(0, 0)$ comme générateur. On a

| $p$ | $\mu_0$ | $\mu_1, \lambda_1$ | $\mu_2, \lambda_2$ | $\tilde{\lambda}_-$ | $\tilde{\lambda}_+$ |
|---|---|---|---|---|---|
| 17 | $\infty$ | 0,1 | 0, 19 | 3 | 1 |
| 19 | $\infty$ | 0,1 | 0, 19 | 1 | 1 |

La conjecture principale est vraie pour $p = 17$ et pour $p = 19$. Sur la $\mathbb{Z}_{17}$-extension et sur la $\mathbb{Z}_{19}$, le rang reste égal à 1. La croissance du groupe de Shafarevich est contrôlée par les valeurs suivantes de $\lambda_+$ et $\lambda_-$ :

$$(53/18, 1/18) \text{ pour } p = 17$$
$$(19/20, 1/20) \text{ pour } p = 19$$

On trouve en fait pour le régulateur $p$-adique

$$(4 \times 17 + 11 \times 17^2 + O(17^3))\omega_E - 17(8 \times 17 + 12 \times 17^2 + O(17^3))\varphi\omega_E$$

$$(13 \times 19 + 10 \times 19^2 + O(19^3))\omega_E - 19(18 \times 19 + 7 \times 19^2 + O(19^3))\varphi\omega_E$$

Les exemples suivants se traitent de la même manière avec

$$L(E/\mathbb{Q}_p, 1)(1-\varphi)^{-1}(1-p^{-1}\varphi^{-1})\mathbf{1}(L'_p(E)) = \mathbb{L}_1\omega_E - p\mathbb{L}_2\varphi\omega_E$$

| $E$ | équation | $p$ | $P$ | $\mathbb{L}_1, \mathbb{L}_2$ | Tam | tors |
|---|---|---|---|---|---|---|
| *53A* | *1,-1,1,0,0* | *5* | *[0,0]* | $3 \cdot 5 + 2 \cdot 5^2 + O(5^5), 4 \cdot 5 + 5^2 + 2 \cdot 5^4 + O(5^5)$ | *1* | *1* |
| *91B* | *0,1,1,-7,5* | *11* | *[-1,3]* | $7 \cdot 11 + 8 \cdot 11^2 + O(11^3), 7 \cdot 11 + O(11^3)$ | *1* | *3* |
| *106B* | *1, 1, 0,-7,5* | *7* | *[2,1]* | $4 \cdot 7^2 + O(7^3), 6 \cdot 7 + 2 \cdot 7^2 + O(7^3)$ | *2* | *1* |
| *145A* | *1,-1,1,-3,2* | *3* | *[0,1]* | $3 + 2 \cdot 3^3 + 3^5 + O(3^6), 2(3 + 3^3 + 3^4 + 3^5 + O(3^6))$ | *1* | *2* |

*Pour ces courbes, la composante $p$-primaire du groupe de Shafarevich-Tate est triviale. La conjecture principale $p$-adique est vraie. J'ai d'autre part vérifié que le logarithme $p$-adique des points donnés est de valuation 1 et même que les calculs sont compatibles avec la conjecture de Birch et Swinnerton-Dyer $p$-adique.*

**Exemple.** Soit $E = 17A^{(-239)}$. On a $L(E,1) = 0$ et

| $E_0$ | équation | $D$ | $p$ | $\mu_0$ | $\mu_1, \lambda_1$ | $\mu_2, \lambda_2$ | $\mu_3, \lambda_3$ | $\tilde{d}a_-$ | $\tilde{\lambda}_+$ |
|---|---|---|---|---|---|---|---|---|---|
| 17A | [1,-1,1,-1,-14] | -239 | 3 | $\infty$ | 1,1 | 0,3 | 0,15 | 9 | 1 |

La conjecture principale est vraie.

**Exemple.** Pour les courbes suivantes, la conjecture principale est vraie, le rang de $\check{S}_p(E/\mathbb{Q})$ est 1 et $\mathbf{III}(E/\mathbb{Q})(p)$ est divisible, donc trivial si l'on peut exhiber un point d'ordre infini :

| $E_0$ | équation | $D$ | $p$ | $\mu_0$ | $\mu_1, \lambda_1$ | $\mu_2, \lambda_2$ | $\mu_3, \lambda_3$ | $\tilde{\lambda}_-$ | $\tilde{\lambda}_+$ |
|---|---|---|---|---|---|---|---|---|---|
| 73A | [1,-1,0,4,-3] | -211 | 3 | $\infty$ | 1,1 | 0,3 | 0,13 | 1 | 1 |
| 102C | [1,0,1,-256,1550] | -187 | 5 | $\infty$ | 1,1 | 0,5 | 0,25 | 5 | 1 |
| 102C | [1,0,1,-256,1550] | -187 | 5 | $\infty$ | 1,1 | 0,5 | 0,25 | 5 | 1 |
| 34A | [1,0,0,-3,1] | -11 | 5 | $\infty$ | 1,1 | 0,5 | 0,25 | 5 | 1 |
| 24A | [0,-1,0,-4,4] | -211 | 7 | $\infty$ | 1,1 | 0,7 | 0,49 | 1 | 7 |
| 98A | [1, 1, 0, -25, -111] | 269 | 5 | $\infty$ | 1,3 | 0,5 | 0,25 | 5 | 1 |
| 38A | [1, 0, 1, 9, 90] | 428 | 5 | $\infty$ | 2,1 | 0,5 | 0,25 | 5 | 1 |
| 14A | [1,0,1,4,-6] | 12 | 5 | $\infty$ | $\infty$ | 0,5 | 0,25 | 5 | 1 |
| 14A | [1,0,1,4,-6] | 185 | 11 | $\infty$ | 0,5 | 0,11 | | 5 | 1 |
| 11A | [0,-1,1,-10,-20] | 61 | 19 | $\infty$ | 0,3 | 0,19 | | 3 | 1 |
| 11A | [0,-1,1,-10,-20] | 65 | 19 | $\infty$ | 0,5 | 0,19 | | 5 | 1 |

**6.3.2. Proposition.** *Supposons $\rho_p$ surjective, $a_p = 0$, $L(E,1) = 0$ et $L^{(*)}_{p,\omega_E}(E, \mathbf{1}) \neq 0$. Supposons qu'il existe un entier $k$ tel que $\mu_k = 1$, $\mu_{k+2} = 0$, $\tilde{\lambda}_{-\epsilon} - \tilde{\lambda}_\epsilon \geq M_{k+1}$ et $\lambda_k = 1 + M_k$ avec $\epsilon = (-1)^{k+1}$. Alors, la conjecture principale est vraie et le rang de $\check{S}_p(E)$ est égal à 1. Mêmes conclusions que dans la proposition 6.3.1.*

On n'a pas forcément $\tilde{\lambda}_- = 1$ car $\mu_1$ est non nul.

*Démonstration.* On a comme précédemment $\lambda(g) = \tilde{\lambda}_{-\epsilon} - \tilde{\lambda}'_{-\epsilon} = \tilde{\lambda}_\epsilon - \tilde{\lambda}'_\epsilon$. Supposons $\mu'_k = 0$, on a nécessairement $\lambda'_k < p^k$, donc $\lambda'_{-\epsilon} < p^k - M_k$. D'autre part, $\tilde{\lambda}'_\epsilon \geq 1$, donc

$$\tilde{\lambda}_{-\epsilon} - p^k + M_k \leq \lambda(g) \leq \tilde{\lambda}_\epsilon - 1$$

et

$$\tilde{\lambda}_{-\epsilon} - \tilde{\lambda}_\epsilon < p^k - 1 - M_k = M_{k+1}$$

ce qui est contradictoire avec l'hypothèse faite. Donc $\mu'_k = 1$. On en déduit que $(P_k/p) \equiv g(Q_k/p) \mod \omega_k$ dans $\mathbb{Z}_p[x]/\omega_k(x)$ et donc que $\lambda(g) = \lambda_k - \lambda'_k \leq \lambda_k - 1 - M_k = 0$. Comme $\mu(g)$ est nul, $g$ est une unité et la conjecture principale est vraie. $\square$

**Exemple.** La proposition s'applique pour les courbes suivantes :

| $E_0$ | équation | $D$ | $p$ | $\mu_0$ | $\mu_1, \lambda_1$ | $\mu_2, \lambda_2$ | $\mu_3, \lambda_3$ | $\mu_4, \lambda_4$ | $\tilde{\lambda}_-$ | $\tilde{\lambda}_+$ |
|---|---|---|---|---|---|---|---|---|---|---|
| 62A | [1,-1,1,-1,1] | -40 | 3 | $\infty$ | 1,1 | 0,5 | 0,11 | | 5 | 3 |
| 46A | [1,-1,0,-10,-12] | 29 | 3 | $\infty$ | 1,1 | 0,5 | 0,13 | | 7 | 3 |
| 52A | [0,0,0,1,-10] | 293 | 3 | $\infty$ | 1,1 | 0,5 | 0,15 | | 9 | 3 |
| 94A | [1,-1,1,0,-1] | 137 | 3 | $\infty$ | 1,1 | 0,5 | 0,17 | | 11 | 3 |
| 62A | [1,-1,1,-1,1] | -59 | 3 | $\infty$ | 1,1 | 1,3 | 0,9 | 0,29 | 3 | 9 |
| 38A | [1,0,1,9,90] | -132 | 5 | $\infty$ | 1,1 | 0,7 | 0,27 | | 7 | 3 |
| 43A | [0,1,1,0,0] | -4 | 7 | $\infty$ | 1,1 | 0,7 | 0,49 | | 7 | 1 |
| 17A | [1,-1,1,-1,-14] | 493 | 3 | $\infty$ | 1,1 | 0,3 | 0,9 | | 3 | 1 |
| 17A | [1,-1,1,-1,-14] | -19 | 3 | $\infty$ | $\infty$ | 0,3 | 0,9 | | 3 | 1 |

Pour la dernière courbe, on a de plus un nouvel élément d'ordre infini dans $\check{S}_p(E/\mathbb{Q}(\mu_9)^\Delta)$. Ainsi, il est prévisible que $\tilde{\lambda}_- \geq 1+2 = 3$. Comme $\tilde{\lambda}_- = 3$, les rangs de $\check{S}_p(E/\mathbb{Q}(\mu_9)^\Delta)$ et de $\check{S}_p(E/\mathbb{Q}_\infty)$ sont exactement égaux à 3. Cela n'empêche pas $\tilde{\lambda}_+$ d'être égal à 1.

**Exemple.** Soit $E = 115A^{(-127)}$ la courbe d'équation $y^2 + y = x^3 + 112903x + 22020117$. Le point $P = (11049/484, 52817151/10648)$ appartient à $E(\mathbb{Q})$ et est d'ordre infini dont la puissance quatrième $Q$ a comme logarithme $2 \times 3^2 + 3^4$ mod $3^5$. Ce point n'est pas divisible par 3 dans $E(\mathbb{Q})$ car $P$ n'est pas dans le groupe formel $E_1(\mathbb{Q}_3)$ en 3 (l'indice de $E_1(\mathbb{Q}_3)$ dans $E(\mathbb{Q}_3)$ est premier à 3). D'autre part, $L'_{p,\omega_E}(E, \mathbf{1}) = 27u$ mod $3^4$ (avec $u$ une unité 3-adique). Le groupe de Shafarevich-Tate est donc d'ordre premier à 3. Les invariants sont $\mu_1 = 2$, $\mu_2 = 1$, $\mu_3 = \mu_4 = 0$, $\lambda_1 = 1$, $\lambda_2 = 3$, $\lambda_3 = 11$, $\lambda_4 = 27$. On a $\tilde{\lambda}_- = 5$, $\tilde{\lambda}_+ = 7$. La conjecture principale est donc montrée. On a enfin $E(\mathbb{Q}_\infty) = E(\mathbb{Q})$ et $\mathbf{III}(E/\mathbb{Q}_n)(3)$ est fini pour tout $n$.

6.4. **Courbes de rang $\geq 2$.** *Pour tout entier $n$, posons $\epsilon_n = 0$ ou $1$ selon que $\delta(L_p(E/\mathbb{Q})) = 0$ ou non pour $\delta$ caractère d'ordre $p^n$ (ou la composante) et soit l'ordre de multiplicité $\sigma_n$ de $L_p(E/\mathbb{Q})$ en $\delta$ et $\sigma'_n = \mathrm{rg}_{\mathbb{Z}_p} \check{S}_p(\mathbb{Q}_n) - \mathrm{rg}_{\mathbb{Z}_p} \check{S}_p(\mathbb{Q}_{n-1})$. Alors, on a $\sigma'_n \leq \sigma_n$ et*

$$\mathrm{rg}_{\mathbb{Z}} E(\mathbb{Q}_\infty) \leq \sum_n p^{n-1}(p-1)\sigma'_n \leq \sum_n p^{n-1}(p-1)\sigma_n$$

*et*

$$\sigma_0 + \sum_{n>0} \epsilon_n(\sigma_n - 1)p^{n-1}(p-1) + \sum_{\substack{n>0 \\ n\equiv 0 \bmod 2}} \epsilon_n p^{n-1}(p-1) \leq \tilde{\lambda}_+$$

$$\sigma_0 + \sum_{n>0} \epsilon_n(\sigma_n - 1)p^{n-1}(p-1) + \sum_{n\equiv 1 \bmod 2} \epsilon_n p^{n-1}(p-1) \leq \tilde{\lambda}_-$$

*ce qui peut aussi s'écrire*

$$\sigma_0 + \sum_{\substack{n>0 \\ n\equiv 0 \bmod 2}} \sigma_n p^{n-1}(p-1) + \sum_{n\equiv 1 \bmod 2} \epsilon_n(\sigma_n - 1)p^{n-1}(p-1) \leq \tilde{\lambda}_+$$

$$\sigma_0 + \sum_{n\equiv 1 \bmod 2} \sigma_n p^{n-1}(p-1) + \sum_{\substack{n>0 \\ n\equiv 0 \bmod 2}} \epsilon_n(\sigma_n - 1)p^{n-1}(p-1) \leq \tilde{\lambda}_-$$

*donc,*

$$2\sigma_0 + \sum_{n, \sigma_n \geq 1} (2\sigma_n - 1)p^{n-1}(p-1) \leq \tilde{\lambda}_- + \tilde{\lambda}_+$$

*et de même en remplaçant $\sigma_n$ par $\sigma'_n$. On peut en déduire une borne du rang de $E(\mathbb{Q}_\infty)$ en fonction de $\lambda_+$ et $\lambda_-$.*

**6.4.1. Proposition.** *Supposons $E(\mathbb{Q})$ de rang $\geq r$. Si $\tilde{\lambda}_+$ ou $\tilde{\lambda}_-$ est égal à $r$, alors la conjecture principale est vraie, $E(\mathbb{Q})$ est de rang $r$ et $\text{III}(E/\mathbb{Q})(p)$ est fini.*

**Exemple.** Reprenons la courbe $E = 1909A$ d'équation $y^2 + y = x^3 - 4x + 2$ et $p = 3$. La courbe $E$ est de rang $\geq 2$ sur $\mathbb{Q}$ et admet comme points indépendants les points $(-2, 1)$, $(0, 1)$. Les calculs donnent $\mu_1 = \mu_2 = 0$, $\lambda_1 = 2$, $\lambda_2 = 4$ et donc $\tilde{\lambda}_- = 2$, $\tilde{\lambda}_- = 2$. La conjecture principale est donc vraie. On a

$$L(E/\mathbb{Q}_3, 1)(1-\varphi)^{-1}(1-3^{-1}\varphi^{-1})\mathbf{1}(L_3''(E))$$
$$\sim 2((3^2 + O(3^3))\omega_E - 3(3^2 + O(3^3))\varphi\omega_E)$$

Comme elle est non nulle, la conjecture de Birch-Swinnerton 3-adique est vraie à une unité près. Le nombre de Tamagawa $\text{Tam}(E)$ est 1. Par le calcul, le régulateur $p$-adique des deux points est

$$(2\times 3^2 + 2\times 3^3 + 2\times 3^4 + 2\times 3^5 + O(3^7))\omega_E - 3(2\times 3^2 + 3^3 + O(3^5))\varphi\omega_E \ .$$

On en déduit que $\text{III}(E/\mathbb{Q})(3)$ est nul. Il n'est malheureusement pas possible pour l'instant de revenir à la conjecture complexe.

**Exemple.** Soit la courbe $E = 70A^{(-299)}$.

| $E_0$ | équation | $D$ | $p$ | $\mu_0$ | $\mu_1, \lambda_1$ | $\mu_2, \lambda_2$ | $\mu_3, \lambda_3$ | $\mu_4, \lambda_4$ | $\tilde{\lambda}_-$ | $\tilde{\lambda}_+$ |
|---|---|---|---|---|---|---|---|---|---|---|
| 70A | [1, -1, 1, 2, -3] | -299 | 3 | $\infty$ | 1,2 | 1,4 | 0,10 | 0,28 | 4 | 8 |

La courbe $E$ a deux points indépendants qui sont $(-132, 5356)$ et $(42625, 8779396)$. Leur régulateur est $(3^4 + O(3^6))\omega_E - 3(3^4 + O(3^5))\varphi\omega_E$. Il est facile de vérifier qu'ils engendrent un sous-groupe de $E(\mathbb{Q})$ d'indice premier à 3. De plus, le signe de l'équation fonctionnelle de $L_p(E)$ est $+1$, donc $\mathbf{1}(L_p'(E)) = 0$. Par contre,

$$L(E/\mathbb{Q}_3, 1)(1-\varphi)^{-1}(1-3^{-1}\varphi^{-1})\mathbf{1}(L_3''(E))$$
$$= (3^4 + O(3^5))\omega_E - 3(3^4 + O(3^5))\varphi(\omega_E) \ .$$

On en déduit que $\text{III}(E/\mathbb{Q})(3)$ est trivial, que le rang de $E(\mathbb{Q})$ est 2, que la conjecture principale est vraie, que $\text{III}(E/\mathbb{Q}_n)(3)$ est fini pour tout entier $n$ et que $E(\mathbb{Q}_\infty) = E(\mathbb{Q})$. En utilisant le fait que $\text{Tam}(E/\mathbb{Q}) = 32$, que le sous-groupe de torsion de $E(\mathbb{Q})$ est d'ordre 2, que les points trouvés engendrent $E(\mathbb{Q})$ modulo torsion et que $\mathbf{1}(L_p^*(E)) = \mathbf{1}(L_p''(E))/2$, ces calculs sont compatibles avec la conjecture $p$-adique de Birch et Swinnerton-Dyer.

**Exemple.** $E = 17A^{(-56)}$ :

| $E_0$ | équation | $D$ | $p$ | $\mu_0$ | $\mu_1, \lambda_1$ | $\mu_2, \lambda_2$ | $\tilde{\lambda}_-$ | $\tilde{\lambda}_+$ |
|---|---|---|---|---|---|---|---|---|
| 17A | [1, -1, 1, -1, -14] | -56 | 3 | $\infty$ | 0,2 | 0,4 | 2 | 2 |

Le régulateur $p$-adique des points $(-42, 1568)$ et $(4, 1560)$ est

$$(2\times 3^2 + 3^3 + 2\times 3^4 + 3^5 + O(3^6))\omega_E - 3(3^4 + O(3^5))\varphi\omega_E \ .$$

La dérivée première de $L_p(E)$ est nulle en $\mathbf{1}$ et on a

$$L(E/\mathbb{Q}_3, 1)(1-\varphi)^{-1}(1-3^{-1}\varphi^{-1})\mathbf{1}(L_3''(E))$$
$$= (3^2 + 3^3 + 2\times 3^4 + O(3^5))\omega_E - 3(2\times 3^4 + O(3^5))\varphi\omega_E$$

On en déduit que $\text{III}(E/\mathbb{Q})(3)$ est trivial, que le rang de $E(\mathbb{Q})$ est 2, que la conjecture principale est vraie, que $\text{III}(E/\mathbb{Q}_n)(3)$ est fini pour tout entier $n$ et que $E(\mathbb{Q}_\infty) = E(\mathbb{Q})$. En utilisant le fait que le sous-groupe de torsion de $E(\mathbb{Q})$ est d'ordre 2, que $\text{Tam}(E/\mathbb{Q}) = 16$ et que $\mathbf{1}(L_p^*(E)) = \mathbf{1}(L_p''(E))/2$, on vérifie la compatibilité avec la conjecture $p$-adique de Birch et Swinnerton-Dyer.

**Exemple.** $E = 145A^{(401)}$

| $E_0$ | équation | $D$ | $p$ | $\mu_0$ | $\mu_1, \lambda_1$ | $\mu_2, \lambda_2$ | $\mu_3, \lambda_3$ | $\tilde{\lambda}_-$ | $\tilde{\lambda}_+$ |
|---|---|---|---|---|---|---|---|---|---|
| 145A | [1, -1, 1, -3, 2] | 401 | 3 | $\infty$ | 1, 1 | 0, 5 | 0, 9 | 3 | 3 |

Le programme `mwrank` donne comme générateurs

$$P_1 = (3910, 239246) \ , \ P_2 = (\frac{7475822}{5041}, \frac{18407343402}{357911}) \ , P_3 = (\frac{694073}{169}, \frac{566652579}{2197}) \ .$$

Le régulateur $p$-adique de ces points est

$$(2 \times 3^3 + 2 \times 3^4 + 2 \times 3^5 + 2 \times 3^6 + O(3^8))\omega_E - 3(2 \times 3^3 + 2 \times 3^4 + 3^5 + O(3^7))\varphi\omega_E \ .$$

La conjecture principale est vraie et le rang de $E(\mathbb{Q})$ est bien 3 car $\tilde{\lambda}_+ = 3$. Quand à la dérivée troisième, on trouve que

$$L(E/\mathbb{Q}_3, 1)(1-\varphi)^{-1}(1-3^{-1}\varphi^{-1})\mathbf{1}(L_p^{(3)}(E))$$
$$= (3^4 + O(3^5))\omega_E - 3(3^4 + O(3^5))\varphi\omega_E$$

Le sous-groupe de torsion de $E(\mathbb{Q})$ est d'ordre 2, $\text{Tam}(E/\mathbb{Q}) = 4$ et $\mathbf{1}(L^{(*)}(E/\mathbb{Q})) = \mathbf{1}(L^{(3)}(E/\mathbb{Q}))/6$. On déduit de nouveau que $\text{III}(E/\mathbb{Q})(3)$ est trivial.

Remarquons que dans la base $(\omega_E, -3\varphi\omega_E)$, la pente du régulateur des deux points $(P_1, P_2)$ est $2 + 3 + 3^2 + 2 \times 3^3 + O(3^4)$, la pente du régulateur de $(P_2, P_3)$ est $2 \times 3^3 + 3^4 + 2 \times 3^6 + O(3^8)$ et la pente du régulateur des deux points $(P_1, P_3)$ est $3^3 + 2 \times 3^5 + 2 \times 3^7 + O(3^8)$, alors que la pente de $(1-\varphi)^{-1}(1-3^{-1}\varphi^{-1})\mathbf{1}(L_p^{(3)}(E))$ est $1 + O(3)$.

## Annexe A

**A.1. Critères de surjectivité de $\rho_p$.** *Nous donnons dans ce paragraphe des conditions pour que l'image de $G_\mathbb{Q}$ par la représentation $\rho_p = \rho_{p,E}$ donnant l'action de $G_\mathbb{Q}$ sur les points de $p$-torsion de $E$ soit égale à $GL_2(\mathbb{Z}/p\mathbb{Z})$ en supposant que $E$ a bonne réduction supersingulière et n'a pas multiplication complexe. Ces critères sont suffisants pour tous les exemples numériques que nous avons calculé. Ils sont tirés de* [23].

**A.1.1. Proposition.** *Soit $G$ un sous-groupe de $GL_2(\mathbb{Z}/p\mathbb{Z})$.*

*1 Si $p$ divise l'ordre de $G$, alors soit $G$ contient $SL_2(\mathbb{Z}/p\mathbb{Z})$, soit $G$ est contenu dans un sous-groupe de Borel de $GL_2(\mathbb{Z}/p\mathbb{Z})$.*

*2 Soit $G$ un sous-groupe de $GL_2(\mathbb{Z}/p\mathbb{Z})$ contenant un sous-groupe de Cartan $C$ (et $p \neq 5$ si $C$ est déployé). Alors soit $G = GL_2(\mathbb{Z}/p\mathbb{Z})$, soit $G$ est contenu dans un sous-groupe de Borel, soit $G$ est contenu dans le normalisateur d'un sous-groupe de Cartan.*

**A.1.2. Lemme.** *Supposons que $E$ a bonne réduction supersingulière en $p$. Alors, $\rho_p(G_\mathbb{Q})$ contient le normalisateur d'un sous-groupe de Cartan non déployé. En particulier, si $p \neq 2$, soit $\rho_p(G_\mathbb{Q}) = GL_2(\mathbb{Z}/p\mathbb{Z})$, soit $\rho_p(G_\mathbb{Q})$ est le normalisateur d'un sous-groupe de Cartan non déployé.*

[23, prop. 12]

Pour assurer que $\rho_p(G_\mathbb{Q}) = GL_2(\mathbb{Z}/p\mathbb{Z})$, il suffit de rajouter une condition pour que $\rho_p(G_\mathbb{Q})$ contienne un élément d'ordre $p$. Cela n'est bien sûr possible que si $E$ n'a pas multiplication complexe. Commençons par le cas $p = 3$.

**A.1.3. Lemme.** *Supposons que $E$ a réduction supersingulière en 3. Alors, $\rho_3(G_\mathbb{Q}) = GL_2(\mathbb{Z}/3\mathbb{Z})$ si et seulement si son discriminant $\Delta_E$ n'est pas un cube.*

Autrement dit, il existe un nombre premier divisant $N_E$ tel que $\text{ord}_\ell(\Delta_E) \not\equiv 0 \bmod 3$. En effet, l'ordre de $G = \rho_3(G_\mathbb{Q})$ est divisible par 3 si et seulement si $\Delta_E$ n'est pas un cube. Lorsque $j \neq 0$, cela est équivalent à ce que $j$ ne soit pas un cube.

**A.1.4. Lemme.** *Supposons que $E$ a bonne réduction supersingulière en $p$. S'il existe un nombre premier $\ell$ divisant strictement $N_E$ tel que $\text{ord}_\ell(j_E) \not\equiv 0 \bmod p$, alors $\rho_p(G_\mathbb{Q}) = GL_2(\mathbb{Z}/p\mathbb{Z})$.*

*Pour $p \neq 3$, cela impose que $E$ n'a pas potentiellement bonne réduction.*

**A.1.5. Proposition.** *Supposons $E$ semi-stable et ayant bonne réduction supersingulière en $p$. Alors, si $p \geq 7$, $\rho_p(G_{\mathbb{Q}}) = GL_2(\mathbb{Z}/p\mathbb{Z})$. Il en est de même pour tous les twists de $E$ par un discriminant premier à $pN_E$.*

**A.1.6. Proposition.** *Supposons $p \geq 5$ et $E$ ayant bonne réduction supersingulière en $p$. S'il existe un nombre premier $\ell$ tel que $a_l \not\equiv 0 \bmod p$, $a_l^2 \not\equiv 4\ell \bmod p$, $(\frac{a_l^2 - 4\ell}{p}) = 1$, alors $\rho_p(G_{\mathbb{Q}}) = GL_2(\mathbb{Z}/p\mathbb{Z})$.*

*Cela se déduit de la proposition 19 de [23] en utilisant le fait que $\rho_p(G_{\mathbb{Q}})$ contient le normalisateur d'un sous-groupe de Cartan non déployé.*

**A.2. Constante de Manin.** *Soit $E$ une courbe elliptique sur $\mathbb{Q}$ ayant bonne réduction supersingulière en $p$ et $\pi : X_0(N_E) \to E$ une paramétrisation minimale. Alors la constante $c_\pi$ est première à $p$. En effet, il est démontré que si $E$ est la courbe de Weil forte, cela est le cas. Comme $E$ a bonne réduction supersingulière, il n'existe pas d'isogénie sur $\mathbb{Q}$ de degré $p$. En effet, dans le cas contraire, l'image de $G_{\mathbb{Q}}$ dans $GL_2(\mathbb{Z}/p\mathbb{Z})$ serait contenue dans un sous-groupe de Borel, ce qui est impossible car elle contient le normalisateur d'un sous-groupe de Cartan. On démontre donc en même temps que le sous-groupe de $E(\mathbb{Q})$ est d'ordre premier à $p$.*

## Annexe B. Du côté des calculs

**B.1. Précisions numériques.** *Pour contrôler la précision numérique dans les calculs, on utilise les affirmations suivantes de [15]. Soit $F$ une fonction localement analytique. Si $F(x) = \sum_i c_i(x-a)_p^i$ est le développement en $a$, soit $I(F, a, n)$ l'idéal engendré par les $c_j p^{(j-1)n}$ pour $j \geq 1$. Posons $I(F, n) = \sum_{a \bmod p^n} I(F, a, n)$. Alors avec $M_E' = \mathbb{Z}_p \omega_E - p\varphi \omega_E \mathbb{Z}_p$,*

$$\int_{\mathbb{Z}_p^*} F d\mu - \sum_{a \bmod p^n} F(a) \mu(a + p^n \mathbb{Z}_p) \in C^{-1} p^{\lfloor n/2 \rfloor} I(F, n) M_E' \ .$$

*Numériquement, pour décider de la précision à utiliser, nous avons besoin de ces congruences. Donnons quelques exemples :*

**B.1.1.** *Prenons les fonctions puissances $x^k$. On a alors $I(x^k, n) \subset \mathbb{Z}_p$. En effet $x^k = \sum_j \binom{k}{j} a^j (x-a)^j$. On a donc en appliquant la formule pour $n = 1$*

$$\int_{\mathbb{Z}_p^*} x^k d\mu - \sum_{a \bmod p} a^k \mu(a + p\mathbb{Z}_p) \in M_E' \ .$$

*Comme $\mu(a + p\mathbb{Z}_p) \in p^{-1} M_E'$ et que $a^k$ est un entier (et même une unité), on en déduit que $L_p(E)(\chi^k) \in p^{-1} M_E'$. Ensuite en prenant $n$ quelconque,*

$$\int_{\mathbb{Z}_p^*} x^k d\mu - \sum_{a \bmod p^n} a^k \mu(a + p^n \mathbb{Z}_p) \in p^{\lfloor n/2 \rfloor} M_E' \ .$$

*Si $k$ est congru à $0 \bmod p^{s-1}(p-1)$, on peut être plus précis sur l'idéal $I(x^k, n)$. Par définition, c'est l'idéal engendré par les $\binom{k}{j} p^{(j-1)n}$. Or on a $\binom{k}{j} = \frac{k}{j}\binom{k-1}{j-1} \in \frac{k}{j}\mathbb{Z}_p$ pour $j \geq 1$. Or pour $j \geq 2$ et $n \geq 1$,*

$$\mathrm{ord}_p((j-1)^{-1} p^{(j-1)n}) \geq j - 1 - \mathrm{ord}_p(j-1) \geq 0 \ .$$

*Donc pour $k \equiv 0 \bmod p^{s-1}(p-1)$, $I(x^k, n) \subset p^{s-1} \mathbb{Z}_p$,*

$$\int_{\mathbb{Z}_p^*} x^k d\mu - \sum_{a \bmod p^n} a^k \mu(a + p^n \mathbb{Z}_p) \in p^{s-1+\lfloor n/2 \rfloor} M_E' \ .$$

*Ainsi, lorsque $\int_{\mathbb{Z}_p^*} d\mu = 0$, la valeur de $\int_{\mathbb{Z}_p^*} x^k d\mu$ est de plus en plus divisible par 0.*

B.1.2. *Comparons $x^k$ et $x^{k'}$ pour $k \equiv k' \bmod p^{s-1}(p-1)$, ce qui revient à calculer l'idéal $I(x^k(x^{p^s(p-1)} - 1), n) \subset I(x^j - 1, n)$ avec $j = p^{s-1}(p-1)$. Le $p-1$ sert (peut-être) par $a^{p-1} \equiv 1 \bmod p$. Je ne pense pas qu'on puisse faire mieux que $I(*, n) = p^{s-1}\mathbb{Z}_p$. On en déduit que*

$$\int_{\mathbb{Z}_p^*} x^k d\mu - \int_{\mathbb{Z}_p^*} x^{k'} d\mu \equiv \sum_{a \bmod p^n} (a^k - a^{k'})\mu(a + p^n\mathbb{Z}_p) \bmod p^{s-1+\lfloor n/2 \rfloor} M'_E$$

$$\in p^{s-(\lfloor n/2 \rfloor+1)}M'_E + p^{s-1+\lfloor n/2 \rfloor}M'_E = p^{s-1+\lfloor n/2 \rfloor}M'_E$$

*car $a^k - a^{k'} \equiv 0 \bmod p^s$. En prenant $n = 1$, on obtient*

$$\int_{\mathbb{Z}_p^*} x^k d\mu \equiv \int_{\mathbb{Z}_p^*} x^{k'} d\mu \bmod p^{s-1}M'_E .$$

*Autrement dit, si $k \equiv k' \bmod p^s(p-1)$, on a*

$$L_p(E)(\chi^k) \equiv L_p(E)(\chi^{k'}) \bmod p^s M'_E .$$

B.1.3. *Passons au logarithme. Du développement $\log_p x = \log_p a + \sum_j (-1)^{j-1} \frac{(x-a)^j}{j}$, on déduit que $I(\log_p x, n) \subset \sum_{j \geq 1} j^{-1} p^{n(j-1)} \mathbb{Z}_p = \mathbb{Z}_p$. D'où la congruence que l'on avait utilisé*

$$\int_{\mathbb{Z}_p^*} \log_p x \, d\mu - \sum_{a \bmod p^n} \log_p a \, \mu(a + p^n\mathbb{Z}_p) \in p^{\lfloor n/2 \rfloor} M'_E .$$

*Comme de plus $\log_p a$ est à valeurs dans $p\mathbb{Z}_p$, en appliquant cette formule pour $n = 1$, on en déduit que $\int_{\mathbb{Z}_p^*} \log_p x \, d\mu \in \mathbb{Z}_p$.*

*On a de même $I(\log_p^k x, n) \subset \mathbb{Z}$. Comme*

$$\sum_{a \bmod p^n} \log_p^k a \, \mu(a + p^n\mathbb{Z}_p) \in p^{k-(\lfloor n/2 \rfloor+1)} M'_E .$$

*en prenant $n = 1$, on obtient que $\int_{\mathbb{Z}_p^*} \log_p^k x \, d\mu \in M'_E$. Explicitement,*

$$\int_{\mathbb{Z}_p^*} \log_p^k x \, d\mu \equiv X_m \omega_E - Y_m(p\varphi\omega_E) \bmod p^m M'_E$$

*avec*

(B.1.1)
$$X_m = (-p)^{-m} \sum_{\substack{a=1 \\ (a,p)=1}}^{p^{2m}} \log_p^k a \, x^+(\frac{a}{p^{2m}})$$

$$Y_m = (-p)^{-m-1} \sum_{\substack{a=1 \\ (a,p)=1}}^{p^{2m}} \log_p^k a \, x^+(\frac{a}{p^{2m-1}})$$

B.1.4. *Formule de Taylor. Prenons comme fonction analytique la fonction analytique (dépendant de $k$ et de $r$)*

$$F(x) = \frac{1}{k^{r+1}}(x^k - \sum_{j=0}^{r} \frac{k^j \log_p^j x}{j!}) .$$

*Dit rapidement, $F(x)$ est une combinaison linéaire infinie de $\log_p^s x$ avec $s > r$ à coefficients entiers. On en déduit que $I(F,n) \subset \mathbb{Z}_p$. De plus les valeurs de $F$ sont dans $\mathbb{Z}_p$. On en déduit que*

$$k^{-(r+1)} \left( \int_{\mathbb{Z}_p^*} x^k d\mu - \sum_{j=0}^{r} k^j \frac{\int_{\mathbb{Z}_p^*} \log_p^j x \, d\mu}{j!} \right) \in p^{-1} M_E' \ .$$

*Ainsi,*

$$\int_{\mathbb{Z}_p^*} x^k d\mu - \sum_{j=0}^{r-1} k^j \frac{\int_{\mathbb{Z}_p^*} \log_p^j x \, d\mu}{j!} \equiv k^r \frac{\int_{\mathbb{Z}_p^*} \log_p^j x \, d\mu}{r!} \bmod p^{-1} k^{r+1} M_E'$$

*ou encore*

$$L_p(E)(\chi^k) - \sum_{j=0}^{r-1} k^j \frac{\mathbf{1}(L_p^{(j)}(E))}{j!} \equiv k^r \frac{\mathbf{1}(L_p^{(r)}(E))}{r!} \bmod p^{-1} k^{r+1} M_E'$$

*ou encore lorsque $\mathbf{1}(L_p^{(j)}(E)) = 0$ pour tout $j < r$,*

$$k^{-r} L_p(E)(\chi^k) \equiv \frac{\mathbf{1}(L_p^{(r)}(E))}{r!} \bmod p^{-1} k M_E' \ .$$

*Ainsi, la limite de $k^{-r} L_p(E)(\chi^k)$ lorsque $k$ tend vers $0$ $p$-adiquement est $\frac{\mathbf{1}(L_p^{(r)}(E))}{r!}$. On en déduit la proposition anecdotique suivante :*

**B.1.1. Proposition.** *Dans une base de $D_p(E)$, la limite de la pente de $L_p(E)(\chi^k)$ lorsque $k$ tend vers $0$ $p$-adiquement est égale à la pente de la première dérivée non nulle de $L_p(E)$ en $1$.*

**B.1.5.** *Version tordue par un caractère quadratique. Pour rentabiliser certains calculs et avoir accès à plus de courbes, il est commode de faire les expérimentations relativement aux courbes elliptiques $E^{(d)}$ tordues par un caractère quadratique de discriminant $d$ d'une courbe fixée $E$. On renvoie à [15] et à [3] pour les formules. Par exemple, les polynômes d'interpolation à calculer sont essentiellement (à une unité $p$-adique près) les polynômes*

$$P_n^{d,(j)} = \sum_{\substack{a=0 \\ (a,dp)=1}}^{|d|p^{n+1}} \binom{d}{a} \text{Teich}^j(a) x^{s(d)} (\frac{a}{p^{n+1}|d|}) (1+x)^{r_n(a)}$$

*où $s(d)$ est le signe de $d$.*

## Annexe C. Tableaux

*Dans les tableaux suivants, nous donnons une liste des "types" de $\lambda_i, \mu_i$ que nous avons rencontré dans les calculs. Plus précisément, on donne*

(1) *le nombre premier $p$ (supersingulier pour les courbes étudiées)*
(2) *dans la colonne $i$, la liste $(\lambda_i, \mu_i)$ pour $i = 1, \ldots, 4$ lorsque ces nombres ne se déduisent pas des précédents, c'est-à-dire jusqu'au premier $i$ pair et le premier $i$ impair tel que $\mu_i = 0$.*
(3) *dans la colonne CP,*

   (a) *CP signifie que la conjecture principale est montrée (simplement avec ces nombres)*
   (b) *\*CP signifie qu'il faut de plus vérifier que $L'_{p,\omega_E}(E, \mathbf{1}) \neq 0$*
   (c) *\*CP-rang signifie qu'il faut de plus vérifier $L^{(*)}_{p,\omega_E}(E, \mathbf{1}) \neq 0$ et que le rang de $E(\mathbb{Q})$ est supérieur à $2$.*

(d) **CP si on vérifie de plus que $L_{p,\omega_E}^{(*)}(E,\mathbf{1}) \neq 0$ et que $[L'_p(E),\tilde{\omega}_E]_{D_p(E)} \not\equiv 0 \bmod \xi_1$.

(e) *CP-tam signifie qu'il faut vérifier en plus que* Tam *n'est pas premier à p, ce qui avec les données en question ne peut être que vrai.*

(4) *dans les colonnes $\tilde{\lambda}_-$ et $\tilde{\lambda}_+$, les valeurs de ces invariants, qui se calculent simplement à l'aide des données précédentes.*

*Ainsi, CP, *CP et **CP ne font intervenir que des conditions de type modulaire.*

| $p$ | $\mu_0$ | $\mu_1,\lambda_1$ | $\mu_2,\lambda_2$ | $\mu_3,\lambda_3$ | $\mu_4,\lambda_4$ | $CP$ | $\tilde{\lambda}_-$ | $\tilde{\lambda}_+$ |
|---|---|---|---|---|---|---|---|---|
| 3 | 0 | 0,0 | | | | $CP$ | 0 | 0 |
| 3 | 1 | 0,2 | 0,4 | | | $CP$-tam | 2 | 2 |
| 3 | 1 | 0,2 | 0,6 | | | $CP$ | 2 | 4 |
| 3 | 1 | 0,2 | 0,8 | | | $CP$ | 2 | 6 |
| 3 | 1 | 0,2 | 1,2 | | 0,28 | $CP$ | 2 | 8 |
| 3 | 1 | 0,2 | 1,2 | | 0,32 | $CP$ | 2 | 12 |
| 3 | 1 | 1,0 | 0,4 | 0,10 | | $CP$ | 4 | 2 |
| 3 | 1 | 1,0 | 0,4 | 0,12 | | $CP$ | 6 | 2 |
| 3 | 1 | 1,0 | 0,4 | 0,16 | | $CP$ | 10 | 2 |
| 3 | 1 | 1,0 | 0,4 | 0,18 | | $CP$ | 12 | 2 |
| 3 | 1 | 1,0 | 0,6 | 0,10 | | | 4 | 4 |
| 3 | 1 | 1,0 | 0,6 | 0,12 | | $CP$ | 6 | 4 |
| 3 | 1 | 1,0 | 0,6 | 0,16 | | $CP$ | 10 | 4 |
| 1 | 1 | 1,0 | 1,2 | 0,10 | 0,30 | $CP$ | 4 | 10 |
| 3 | 2 | 0,2 | 0,4 | | | | 2 | 2 |
| 3 | 2 | 0,2 | 0,6 | | | $CP$ | 2 | 4 |
| 3 | 2 | 0,2 | 0,8 | | | $CP$ | 2 | 6 |
| 3 | 2 | 0,2 | 1,4 | | 0,28 | $CP$ | 2 | 8 |
| 3 | 2 | 0,2 | 1,4 | | 0,30 | $CP$ | 2 | 10 |
| 3 | 2 | 0,2 | 1,4 | | 0,32 | $CP$ | 2 | 12 |
| 3 | 2 | 0,2 | 1,6 | | 0,28 | $CP$ | 2 | 8 |
| 3 | 2 | 0,2 | 1,6 | | 0,30 | $CP$ | 2 | 10 |
| 3 | 2 | 0,2 | 1,6 | | 0,32 | $CP$ | 2 | 12 |
| 3 | 2 | 1,2 | 0,4 | 0,10 | | $CP$ | 4 | 2 |
| 3 | 2 | 1,2 | 0,4 | 0,12 | | $CP$ | 6 | 2 |
| 3 | 2 | 1,2 | 0,4 | 0,14 | | $CP$ | 8 | 2 |
| 3 | 2 | 1,2 | 0,6 | 0,10 | | | 4 | 4 |
| 3 | 2 | 1,2 | 0,6 | 0,12 | | $CP$ | 6 | 4 |
| 3 | 2 | 1,2 | 0,6 | 0,16 | | $CP$ | 10 | 4 |
| 3 | 2 | 1,2 | 0,8 | 0,12 | | | 6 | 6 |
| 3 | 2 | 1,2 | 1,4 | 0,10 | 0,28 | $CP$ | 4 | 8 |
| 3 | 2 | 1,2 | 1,4 | 0,12 | 0,28 | $CP$ | 6 | 8 |
| 3 | 2 | 1,2 | 1,6 | 0,10 | 0,28 | $CP$ | 4 | 8 |
| 3 | 2 | 2,0 | 0,4 | 0,10 | | $CP$ | 4 | 2 |
| 3 | 2 | 2,0 | 0,4 | 0,12 | | $CP$ | 6 | 2 |
| 3 | 2 | 2,0 | 0,4 | 0,14 | | $CP$ | 8 | 2 |
| 3 | 2 | 2,0 | 0,4 | 0,18 | | $CP$ | 12 | 2 |
| 3 | 2 | 2,0 | 0,6 | 0,10 | | | 4 | 4 |
| 3 | 2 | 2,0 | 0,6 | 0,14 | | $CP$ | 8 | 4 |
| 3 | 2 | 2,0 | 1,4 | 0,10 | 0,28 | $CP$ | 4 | 8 |
| 3 | 2 | 2,0 | 1,8 | 0,10 | 0,28 | $CP$ | 4 | 8 |
| 3 | 3 | 0,2 | 0,4 | | | $CP$-tam | 2 | 2 |
| 3 | 3 | 0,2 | 0,6 | | | $CP$ | 2 | 4 |
| 3 | 3 | 1,2 | 0,4 | 0,10 | | $CP$ | 2 | 4 |
| 3 | 3 | 1,2 | 0,6 | 0,10 | | | 4 | 4 |
| 3 | 3 | 1,2 | 0,6 | 0,12 | | | 6 | 4 |
| 3 | 3 | 1,2 | 0,6 | 0,16 | | | 10 | 4 |
| 3 | 3 | 2,2 | 0,4 | 0,10 | | $CP$ | 4 | 2 |
| 3 | 3 | 2,2 | 0,4 | 0,12 | | $CP$ | 6 | 2 |
| 3 | 3 | 2,2 | 0,8 | 0,10 | | | 4 | 6 |
| 3 | 4 | 0,2 | 0,4 | | | | 2 | 2 |

| $p$ | $\mu_0$ | $\mu_1,\lambda_1$ | $\mu_2,\lambda_2$ | $\mu_3,\lambda_3$ | $\mu_4,\lambda_4$ | $CP$ | $\tilde{\lambda}_-$ | $\tilde{\lambda}_+$ |
|---|---|---|---|---|---|---|---|---|
| 3 | $\infty$ | 0,1 | 0,3 | | | *CP | 1 | 1 |
| 3 | $\infty$ | 0,1 | 0,5 | | | *CP | 1 | 3 |
| 3 | $\infty$ | 0,1 | 0,7 | | | *CP | 1 | 5 |
| 3 | $\infty$ | 0,1 | 1,3 | | 0,27 | *CP | 1 | 7 |
| 3 | $\infty$ | 0,1 | 1,3 | | 0,29 | *CP | 1 | 9 |
| 3 | $\infty$ | 0,1 | 1,5 | | 0,27 | *CP | 1 | 7 |
| 3 | $\infty$ | 0,1 | 1,5 | | 0,29 | *CP | 1 | 9 |
| 3 | $\infty$ | 0,1 | 1,7 | | 0,27 | *CP | 1 | 7 |
| 3 | $\infty$ | 0,1 | 2,5 | | 0,29 | *CP | 1 | 9 |
| 3 | $\infty$ | 0,1 | $\infty$ | | 0,31 | *CP | 1 | 11 |
| 3 | $\infty$ | 0,1 | $\infty$ | | 0,27 | *CP | 1 | 7 |
| 3 | $\infty$ | 0,2 | 0,4 | | | *CP-rang | 2 | 2 |
| 3 | $\infty$ | 0,2 | 0,6 | | | *CP-rang | 2 | 4 |
| 3 | $\infty$ | 0,2 | 0,8 | | | *CP-rang | 2 | 6 |
| 3 | $\infty$ | 0,2 | 1,4 | | 0,28 | *CP-rang | 2 | 8 |
| 3 | $\infty$ | 0,2 | 1,4 | | 0,30 | *CP-rang | 2 | 10 |
| 3 | $\infty$ | 0,2 | 1,6 | | 0,28 | *CP-rang | 2 | 8 |
| 3 | $\infty$ | 0,2 | 1,6 | | 0,30 | *CP-rang | 2 | 10 |
| 3 | $\infty$ | 0,2 | 1,8 | | 0,30 | *CP-rang | 2 | 10 |
| 3 | $\infty$ | 0,2 | 2,4 | | 0,28 | *CP-rang | 2 | 8 |
| 3 | $\infty$ | 1,1 | 0,3 | 0,9 | | *CP | 3 | 1 |
| 3 | $\infty$ | 1,1 | 0,3 | 0,11 | | *CP | 5 | 1 |
| 3 | $\infty$ | 1,1 | 0,3 | 0,13 | | *CP | 7 | 1 |
| 3 | $\infty$ | 1,1 | 0,3 | 0,15 | | *CP | 9 | 1 |
| 3 | $\infty$ | 1,1 | 0,3 | 0,17 | | *CP | 11 | 1 |
| 3 | $\infty$ | 1,1 | 0,3 | 0,19 | | *CP | 13 | 1 |
| 3 | $\infty$ | 1,1 | 0,5 | 0,9 | | | 3 | 3 |
| 3 | $\infty$ | 1,1 | 0,5 | 0,11 | | *CP | 5 | 3 |
| 3 | $\infty$ | 1,1 | 0,5 | 0,13 | | *CP | 7 | 3 |
| 3 | $\infty$ | 1,1 | 0,5 | 0,15 | | *CP | 9 | 3 |
| 3 | $\infty$ | 1,1 | 0,5 | 0,17 | | *CP | 11 | 3 |
| 3 | $\infty$ | 1,1 | 0,7 | 0,9 | | | 3 | 5 |
| 3 | $\infty$ | 1,1 | 0,7 | 0,11 | | | 5 | 5 |
| 3 | $\infty$ | 1,1 | 0,7 | 0,13 | | *CP | 7 | 5 |
| 3 | $\infty$ | 1,1 | 0,7 | 0,15 | | *CP | 9 | 5 |
| 3 | $\infty$ | 1,1 | 1,3 | 0,9 | 0,27 | *CP | 3 | 7 |
| 3 | $\infty$ | 1,1 | 1,3 | 0,9 | 0,29 | *CP | 3 | 9 |
| 3 | $\infty$ | 1,1 | 1,3 | 0,9 | 0,33 | *CP | 3 | 13 |
| 3 | $\infty$ | 1,1 | 1,3 | 0,11 | 0,27 | | 5 | 7 |
| 3 | $\infty$ | 1,1 | 1,3 | 0,11 | 0,29 | | 5 | 9 |
| 3 | $\infty$ | 1,1 | 1,3 | 0,11 | 0,33 | *CP | 5 | 13 |
| 3 | $\infty$ | 1,1 | 1,3 | 0,13 | 0,27 | | 7 | 7 |
| 3 | $\infty$ | 1,1 | 1,3 | 0,13 | 0,29 | | 7 | 9 |
| 3 | $\infty$ | 1,1 | 1,5 | 0,9 | 0,27 | | 3 | 7 |
| 3 | $\infty$ | 1,1 | 1,5 | 0,9 | 0,29 | | 3 | 9 |
| 3 | $\infty$ | 1,1 | 1,5 | 0,9 | 0,31 | | 3 | 11 |
| 3 | $\infty$ | 1,1 | 1,7 | 0,9 | 0,31 | | 3 | 11 |

| $p$ | $\mu_0$ | $\mu_1,\lambda_1$ | $\mu_2,\lambda_2$ | $\mu_3,\lambda_3$ | $\mu_4,\lambda_4$ | $CP$ | $\tilde{\lambda}_-$ | $\tilde{\lambda}_+$ |
|---|---|---|---|---|---|---|---|---|
| *3* | $\infty$ | *1,2* | *0,4* | *0,10* | | *\*CP-rang* | *4* | *2* |
| *3* | $\infty$ | *1,2* | *0,4* | *0,12* | | *\*CP-rang* | *6* | *2* |
| *3* | $\infty$ | *1,2* | *0,4* | *0,14* | | *\*CP-rang* | *8* | *2* |
| *3* | $\infty$ | *1,2* | *0,4* | *0,16* | | *\*CP-rang* | *10* | *2* |
| *3* | $\infty$ | *1,2* | *0,6* | *0,10* | | | *4* | *4* |
| *3* | $\infty$ | *1,2* | *0,6* | *0,14* | | *\*CP-rang* | *8* | *4* |
| *3* | $\infty$ | *1,2* | *0,8* | *0,10* | | | *4* | *6* |
| *3* | $\infty$ | *1,2* | *0,8* | *0,14* | | *\*CP-rang* | *8* | *6* |
| *3* | $\infty$ | *1,2* | *1,4* | *0,10* | *0,28* | *\*CP-rang* | *4* | *8* |
| *3* | $\infty$ | *2,1* | *0,3* | *0,9* | | *\*CP* | *3* | *1* |
| *3* | $\infty$ | *2,1* | *0,3* | *0,11* | | *\*CP* | *5* | *1* |
| *3* | $\infty$ | *2,1* | *0,3* | *0,13* | | *\*CP* | *7* | *1* |
| *3* | $\infty$ | *2,1* | *0,3* | *0,15* | | *\*CP* | *9* | *1* |
| *3* | $\infty$ | *2,1* | *0,5* | *0,9* | | | *3* | *3* |
| *3* | $\infty$ | *2,1* | *0,5* | *0,11* | | | *5* | *3* |
| *3* | $\infty$ | *2,1* | *0,5* | *0,13* | | | *7* | *3* |
| *3* | $\infty$ | *2,1* | *0,7* | *0,9* | | | *3* | *5* |
| *3* | $\infty$ | *2,1* | *0,7* | *0,11* | | | *5* | *5* |
| *3* | $\infty$ | *2,1* | *0,7* | *0,19* | | | *13* | *5* |
| *3* | $\infty$ | *2,1* | *1,3* | *0,9* | *0,27* | | *3* | *7* |
| *3* | $\infty$ | *2,1* | *1,3* | *0,9* | *0,29* | *\*CP* | *3* | *9* |
| *3* | $\infty$ | *2,1* | *1,3* | *0,11* | *0,27* | | *5* | *7* |
| *3* | $\infty$ | *2,1* | *1,3* | *0,13* | *0,27* | | *7* | *7* |
| *3* | $\infty$ | *2,2* | *0,4* | *0,10* | | *\*CP-rang* | *4* | *2* |
| *3* | $\infty$ | *2,2* | *0,4* | *0,12* | | *\*CP-rang* | *6* | *2* |
| *3* | $\infty$ | *2,2* | *0,6* | *0,12* | | | *6* | *4* |
| *3* | $\infty$ | *2,2* | *0,6* | *0,14* | | | *8* | *4* |
| *3* | $\infty$ | *2,2* | *0,8* | *0,12* | | | *6* | *6* |
| *3* | $\infty$ | *3,1* | *0,3* | *0,9* | | *\*CP* | *3* | *1* |
| *3* | $\infty$ | *3,1* | *0,3* | *0,11* | | *\*CP* | *5* | *1* |
| *3* | $\infty$ | *3,1* | *0,3* | *0,15* | | *\*CP* | *7* | *1* |
| *3* | $\infty$ | *3,1* | *0,5* | *0,9* | | | *3* | *3* |
| *3* | $\infty$ | *3,1* | *0,5* | *0,11* | | | *5* | *3* |
| *3* | $\infty$ | *3,1* | *0,5* | *0,13* | | | *7* | *3* |
| *3* | $\infty$ | *3,1* | *0,5* | *0,15* | | | *9* | *3* |
| *3* | $\infty$ | *3,1* | *1,3* | *0,13* | *0,27* | | *7* | *7* |
| *3* | $\infty$ | *4,1* | *0,3* | *0,9* | | *\*CP* | *3* | *1* |
| *3* | $\infty$ | $\infty$ | *0,3* | *0,11* | | *\*CP* | *5* | *1* |
| *3* | $\infty$ | $\infty$ | *0,3* | *0,13* | | *\*CP* | *7* | *1* |
| *3* | $\infty$ | $\infty$ | *0,3* | *0,17* | | *\*CP* | *11* | *1* |
| *3* | $\infty$ | $\infty$ | *0,4* | *0,10* | | *\*CP-rang* | *4* | *2* |
| *3* | $\infty$ | $\infty$ | *0,4* | *0,12* | | *\*CP-rang* | *6* | *2* |
| *3* | $\infty$ | $\infty$ | *0,5* | *0,9* | | *\*\*CP* | *3* | *3* |
| *3* | $\infty$ | $\infty$ | *0,5* | *0,11* | | | *5* | *3* |
| *3* | $\infty$ | $\infty$ | *0,5* | *0,13* | | | *7* | *3* |
| *3* | $\infty$ | $\infty$ | *0,5* | *0,15* | | | *9* | *3* |
| *3* | $\infty$ | $\infty$ | *0,6* | *0,10* | | | *4* | *4* |
| *3* | $\infty$ | $\infty$ | *0,6* | *0,12* | | | *6* | *4* |
| *3* | $\infty$ | $\infty$ | *0,7* | *0,9* | | *\*\*CP* | *3* | *5* |
| *3* | $\infty$ | $\infty$ | *0,7* | *0,11* | | | *5* | *5* |
| *3* | $\infty$ | $\infty$ | *0,7* | *0,13* | | | *7* | *5* |
| *3* | $\infty$ | $\infty$ | *1,3* | *0,9* | *0,27* | *\*\*CP* | *3* | *7* |
| *3* | $\infty$ | $\infty$ | *1,3* | *0,9* | *0,29* | *\*\*CP* | *3* | *9* |
| *3* | $\infty$ | $\infty$ | *1,3* | *0,9* | *0,31* | *\*\*CP* | *3* | *11* |
| *3* | $\infty$ | $\infty$ | *1,3* | *0,11* | *0,27* | | *5* | *7* |
| *3* | $\infty$ | $\infty$ | *1,4* | *0,10* | *0,28* | | *4* | *8* |
| *3* | $\infty$ | $\infty$ | *1,4* | *0,12* | *0,28* | | *4* | *8* |
| *3* | $\infty$ | $\infty$ | *1,5* | *0,9* | *0,27* | *\*\*CP* | *3* | *7* |

| $p$ | $\mu_0$ | $\mu_1,\lambda_1$ | $\mu_2,\lambda_2$ | $\mu_3,\lambda_3$ | $\mu_4,\lambda_4$ | $CP$ | $\tilde{\lambda}_-$ | $\tilde{\lambda}_+$ |
|---|---|---|---|---|---|---|---|---|
| 5 | 0 | 0,0 | | | | $CP$ | 0 | 0 |
| 5 | 1 | 0,2 | 0,6 | | | $CP$-tam | 2 | 2 |
| 5 | 1 | 0,2 | 0,8 | | | $CP$ | 2 | 4 |
| 5 | 1 | 0,2 | 0,12 | | | $CP$ | 2 | 8 |
| 5 | 1 | 0,4 | 0,6 | | | $CP$ | 4 | 2 |
| 5 | 1 | 1,0 | 0,6 | 0,26 | | $CP$ | 6 | 2 |
| 5 | 2 | 0,2 | 0,6 | | | | 2 | 2 |
| 5 | 2 | 0,2 | 0,8 | | | $CP$ | 2 | 4 |
| 5 | 2 | 0,2 | 0,10 | | | $CP$ | 2 | 6 |
| 5 | 2 | 0,4 | 0,6 | | | $CP$ | 4 | 2 |
| 5 | 2 | 1,2 | 0,6 | 0,26 | | $CP$ | 6 | 2 |
| 5 | 2 | 1,2 | 0,6 | 0,28 | | $CP$ | 8 | 2 |
| 5 | $\infty$ | 0,1 | 0,5 | | | *$CP$ | 1 | 1 |
| 5 | $\infty$ | 0,1 | 0,7 | | | *$CP$ | 1 | 3 |
| 5 | $\infty$ | 0,1 | 0,9 | | | *$CP$ | 1 | 5 |
| 5 | $\infty$ | 0,1 | 0,11 | | | *$CP$ | 1 | 7 |
| 5 | $\infty$ | 0,1 | 0,13 | | | *$CP$ | 1 | 9 |
| 5 | $\infty$ | 0,1 | 0,15 | | | *$CP$ | 1 | 11 |
| 5 | $\infty$ | 0,2 | 0,6 | | | | 2 | 2 |
| 5 | $\infty$ | 0,2 | 0,8 | | | *$CP$-rang | 2 | 4 |
| 5 | $\infty$ | 0,2 | 0,10 | | | *$CP$-rang | 2 | 6 |
| 5 | $\infty$ | 0,2 | 0,12 | | | *$CP$-rang | 2 | 8 |
| 5 | $\infty$ | 0,3 | 0,5 | | | *$CP$ | 3 | 1 |
| 5 | $\infty$ | 0,3 | 0,7 | | | | 3 | 3 |
| 5 | $\infty$ | 0,3 | 0,9 | | | | 3 | 5 |
| 5 | $\infty$ | 0,3 | 0,11 | | | | 3 | 7 |
| 5 | $\infty$ | 0,3 | 0,13 | | | | 3 | 9 |
| 5 | $\infty$ | 0,4 | 0,6 | | | *$CP$-rang | 4 | 2 |
| 5 | $\infty$ | 0,4 | 0,8 | | | | 4 | 4 |
| 5 | $\infty$ | 0,4 | 0,10 | | | | 4 | 6 |
| 5 | $\infty$ | 0,4 | 0,12 | | | | 4 | 8 |
| 5 | $\infty$ | 1,1 | 0,5 | 0,25 | | *$CP$ | 5 | 1 |
| 5 | $\infty$ | 1,1 | 0,5 | 0,27 | | *$CP$ | 7 | 1 |
| 5 | $\infty$ | 1,1 | 0,5 | 0,29 | | *$CP$ | 9 | 1 |
| 5 | $\infty$ | 1,1 | 0,5 | 0,31 | | *$CP$ | 11 | 1 |
| 5 | $\infty$ | 1,1 | 0,7 | 0,25 | | | 5 | 3 |
| 5 | $\infty$ | 1,1 | 0,7 | 0,27 | | *$CP$ | 7 | 3 |
| 5 | $\infty$ | 1,1 | 0,9 | 0,27 | | | 7 | 5 |
| 5 | $\infty$ | 1,1 | 0,9 | 0,25 | | | 5 | 5 |
| 5 | $\infty$ | 1,1 | 0,11 | 0,25 | | | 5 | 7 |
| 5 | $\infty$ | 1,2 | 0,6 | 0,26 | | *$CP$-rang | 6 | 2 |
| 5 | $\infty$ | 1,2 | 0,6 | 0,28 | | *$CP$-rang | 8 | 2 |
| 5 | $\infty$ | 1,2 | 0,8 | 0,26 | | | 6 | 3 |
| 5 | $\infty$ | 1,3 | 0,5 | 0,25 | | *$CP$ | 5 | 1 |
| 5 | $\infty$ | 1,3 | 0,7 | 0,27 | | | 7 | 3 |
| 5 | $\infty$ | 1,3 | 0,5 | 0,27 | | *$CP$ | 7 | 1 |
| 5 | $\infty$ | 1,3 | 0,15 | 0,25 | | | 5 | 11 |
| 5 | $\infty$ | 2,1 | 0,5 | 0,25 | | *$CP$ | 5 | 1 |
| 5 | $\infty$ | $\infty$ | 0,5 | 0,25 | | *$CP$ | 5 | 1 |
| 5 | $\infty$ | $\infty$ | 0,7 | 0,25 | | **$CP$ | 5 | 3 |
| 5 | $\infty$ | $\infty$ | 0,8 | 0,26 | | | 6 | 4 |

| $p$ | $\mu_0$ | $\mu_1,\lambda_1$ | $\mu_2,\lambda_2$ | $\mu_3,\lambda_3$ | $\mu_4,\lambda_4$ | $CP$ | $\tilde{\lambda}_-$ | $\tilde{\lambda}_+$ |
|---|---|---|---|---|---|---|---|---|
| 7 | 0 | 0,0 | | | | CP | 0 | 0 |
| 7 | 2 | 0,2 | 0,8 | | | | 2 | 2 |
| 7 | 2 | 0,2 | 0,10 | | | CP | 2 | 4 |
| 7 | 2 | 0,4 | 0,8 | | | CP | 4 | 2 |
| 7 | ∞ | 0,1 | 0,7 | | | *CP | 1 | 1 |
| 7 | ∞ | 0,1 | 0,9 | | | *CP | 1 | 3 |
| 7 | ∞ | 0,1 | 0,11 | | | *CP | 1 | 5 |
| 7 | ∞ | 0,1 | 0,13 | | | *CP | 1 | 7 |
| 7 | ∞ | 0,1 | 0,15 | | | *CP | 1 | 9 |
| 7 | ∞ | 0,2 | 0,8 | | | *CP-rang | 2 | 2 |
| 7 | ∞ | 0,2 | 0,10 | | | *CP-rang | 2 | 4 |
| 7 | ∞ | 0,3 | 0,7 | | | CP | 3 | 1 |
| 7 | ∞ | 0,3 | 0,9 | | | | 3 | 3 |
| 7 | ∞ | 0,3 | 0,11 | | | | 3 | 5 |
| 7 | ∞ | 0,3 | 0,13 | | | | 3 | 7 |
| 7 | ∞ | 0,4 | 0,8 | | | *CP-rang | 4 | 2 |
| 7 | ∞ | 0,4 | 0,10 | | | | 4 | 4 |
| 7 | ∞ | 0,4 | 0,12 | | | | 4 | 6 |
| 7 | ∞ | 0,5 | 0,7 | | | *CP | 5 | 1 |
| 7 | ∞ | 0,5 | 0,9 | | | | 5 | 3 |
| 7 | ∞ | 0,5 | 0,11 | | | | 5 | 5 |
| 7 | ∞ | 0,5 | 0,13 | | | | 5 | 7 |
| 7 | ∞ | 0,6 | 0,8 | | | *CP-rang | 6 | 2 |
| 7 | ∞ | 0,6 | 0,10 | | | | 6 | 4 |
| 7 | ∞ | 1,1 | 0,7 | 0,49 | | *CP | 7 | 1 |
| 7 | ∞ | 1,1 | 0,9 | 0,49 | | **CP | 7 | 3 |
| 7 | ∞ | 1,4 | 0,10 | 0,50 | | | 8 | 4 |
| 7 | ∞ | ∞ | 0,9 | 0,49 | | **CP | 7 | 3 |
| 11 | 0 | 0,0 | | | | CP | | |
| 11 | ∞ | 0,1 | 0,11 | | | *CP | 1 | 1 |
| 11 | ∞ | 0,1 | 0,13 | | | *CP | 1 | 3 |
| 11 | ∞ | 0,1 | 0,15 | | | *CP | 1 | 5 |
| 11 | ∞ | 0,2 | 0,12 | | | *CP-rang | 2 | 2 |
| 11 | ∞ | 0,3 | 0,11 | | | *CP | 3 | 1 |
| 11 | ∞ | 0,4 | 0,14 | | | | 4 | 4 |
| 11 | ∞ | 0,5 | 0,11 | | | *CP | 5 | 1 |
| 17 | 0 | 0,0 | | | | CP | 0 | 0 |
| 17 | 2 | 0,2 | 0,18 | | | | 2 | 2 |
| 19 | 0 | 0,0 | | | | CP | 0 | 0 |
| 19 | ∞ | 0,3 | 0,19 | | | *CP | 3 | 1 |
| 19 | ∞ | 0,5 | 0,19 | | | *CP | 5 | 1 |

| $p$ | $\mu_0$ | $\mu_1,\lambda_1$ | $\mu_2,\lambda_2$ | $\mu_3,\lambda_3$ | $\mu_4,\lambda_4$ | $CP$ | $\tilde{\lambda}_-$ | $\tilde{\lambda}_+$ |
|---|---|---|---|---|---|---|---|---|
| 3 | 0 | 0,0 | 0,2 | 0,6 | 0,20 | $CP$ | 0 | 0 |
| 3 | 2 | 0,2 | 0,4 | 0,8 | 0,22 | | 2 | 2 |
| 3 | 2 | 0,2 | 0,6 | 0,8 | 0,24 | $CP$ | 2 | 4 |
| 3 | 2 | 0,2 | 0,8 | 0,8 | 0,26 | $CP$ | 2 | 6 |
| 3 | 2 | 0,2 | 1,2 | 0,8 | 0,28 | $CP$ | 2 | 8 |
| 3 | 2 | 1,2 | 0,4 | 0,12 | 0,22 | $CP$ | 6 | 2 |
| 3 | 2 | 1,2 | 0,4 | 0,10 | 0,22 | $CP$ | 4 | 2 |
| 3 | 2 | 1,2 | 0,4 | 0,14 | 0,22 | $CP$ | 8 | 2 |
| 3 | 2 | 1,2 | 0,4 | 0,16 | 0,22 | $CP$ | 10 | 2 |
| 3 | 2 | 1,2 | 0,6 | 0,10 | 0,24 | | 4 | 4 |
| 3 | 2 | 2,0 | 0,4 | 0,10 | 0,22 | $CP$ | 4 | 2 |
| 3 | 2 | 2,0 | 0,4 | 0,12 | 0,22 | $CP$ | 6 | 2 |
| 3 | 2 | 2,0 | 0,8 | 0,10 | 0,26 | $CP$ | 4 | 6 |
| 3 | 4 | 0,2 | 0,4 | 0,8 | 0,22 | | 2 | 2 |
| 3 | $\infty$ | 0,1 | 0,3 | 0,7 | 0,21 | $*CP$ | 1 | 1 |
| 3 | $\infty$ | 0,1 | 0,5 | 0,7 | 0,23 | $*CP$ | 1 | 3 |
| 3 | $\infty$ | 0,1 | 0,7 | 0,7 | 0,25 | $*CP$ | 1 | 5 |
| 3 | $\infty$ | 0,1 | 1,1 | 0,7 | 0,27 | $*CP$ | 1 | 7 |
| 3 | $\infty$ | 0,1 | 1,1 | 0,7 | 0,29 | $*CP$ | 1 | 9 |
| 3 | $\infty$ | 0,2 | 0,4 | 0,8 | 0,22 | | 2 | 2 |
| 3 | $\infty$ | 0,2 | 0,6 | 0,8 | 0,24 | | 2 | 4 |
| 3 | $\infty$ | 0,2 | 0,8 | 0,8 | 0,26 | | 2 | 6 |
| 3 | $\infty$ | 0,2 | 1,2 | 0,8 | 0,28 | | 2 | 8 |
| 3 | $\infty$ | 1,1 | 0,3 | 0,9 | 0,21 | $*CP$ | 3 | 1 |
| 3 | $\infty$ | 1,1 | 0,3 | 0,11 | 0,21 | | 5 | 1 |
| 3 | $\infty$ | 1,1 | 0,3 | 0,15 | 0,21 | $*CP$ | 9 | 1 |
| 3 | $\infty$ | 1,1 | 0,3 | 0,19 | 0,21 | $*CP$ | 13 | 1 |
| 3 | $\infty$ | 1,1 | 0,3 | 0,13 | 0,21 | $*CP$ | 7 | 1 |
| 3 | $\infty$ | 1,1 | 0,5 | 0,9 | 0,23 | | 3 | 3 |
| 3 | $\infty$ | 1,1 | 0,5 | 0,11 | 0,23 | | 5 | 3 |
| 3 | $\infty$ | 1,1 | 0,5 | 0,13 | 0,23 | | 7 | 3 |
| 3 | $\infty$ | 1,1 | 0,5 | 0,17 | 0,23 | | 11 | 3 |
| 3 | $\infty$ | 1,1 | 0,7 | 0,9 | 0,25 | | 3 | 5 |
| 3 | $\infty$ | 1,1 | 0,7 | 0,11 | 0,25 | | 5 | 5 |
| 3 | $\infty$ | 1,2 | 0,4 | 0,10 | 0,22 | | 4 | 2 |
| 3 | $\infty$ | 1,2 | 0,4 | 0,14 | 0,22 | | 8 | 2 |
| 3 | $\infty$ | 1,2 | 0,6 | 0,10 | 0,24 | | 4 | 4 |
| 3 | $\infty$ | 2,1 | 0,3 | 0,9 | 0,21 | | $*CP$ 3 | 1 |
| 3 | $\infty$ | 2,1 | 0,5 | 0,9 | 0,23 | | 3 | 3 |
| 3 | $\infty$ | 2,1 | 0,7 | 0,11 | 0,25 | | 5 | 5 |
| 3 | $\infty$ | 2,2 | 0,4 | 0,10 | 0,22 | | 4 | 2 |
| 3 | $\infty$ | 3,1 | 0,3 | 0,9 | 0,21 | | $*CP$ 3 | 1 |
| 3 | $\infty$ | 3,1 | 0,3 | 0,11 | 0,21 | $*CP$ | 5 | 1 |
| 3 | $\infty$ | $\infty$ | 0,3 | 0,9 | 0,21 | $*CP$ | 3 | 1 |
| 3 | $\infty$ | $\infty$ | 0,3 | 0,11 | 0,21 | $*CP$ | 5 | 1 |
| 3 | $\infty$ | $\infty$ | 0,3 | 0,13 | 0,21 | $*CP$ | 7 | 1 |
| 3 | $\infty$ | $\infty$ | 0,4 | 0,10 | 0,22 | $*CP$ | 4 | 2 |
| 3 | $\infty$ | $\infty$ | 0,5 | 0,9 | 0,23 | $**CP$ | 3 | 3 |
| 3 | $\infty$ | $\infty$ | 0,5 | 0,11 | 0,23 | | 5 | 3 |
| 3 | $\infty$ | $\infty$ | 0,5 | 0,13 | 0,23 | | 7 | 3 |

Table pour $a_3 \neq 0$

## Références

Mathématiques, Bat 425, Université Paris-Sud, F-91 405 Orsay, France
*E-mail address*: `bpr@math.u-psud.fr`